\newcommand{\R}{\mathds R}

\newcommand{\Z}{\mathds Z}
\newcommand{\N}{\mathds N}
\newcommand{\Spl}{\mathrm{Sp}}
\newcommand{\spl}{\mathfrak{sp}}
\newcommand{\GL}{\mathrm{GL}}

\newcommand{\indice}{\mathrm n_-}
\newcommand{\coindice}{\mathrm n_+}
\newcommand{\segnatura}{\mathrm{sign}}

\newcommand{\Ddt}{\tfrac{\mathrm D}{\mathrm dt}}
\newcommand{\ddt}{\tfrac{\mathrm d}{\mathrm dt}}
\newcommand{\DdtR}{\tfrac{\mathrm D_{\scriptscriptstyle{\mathrm R}}}{\mathrm dt}}
\newcommand{\DdttR}{\tfrac{\mathrm D_{\scriptscriptstyle{\mathrm R}}^2}{\mathrm dt^2}}
\newcommand{\Dim}{\mathrm{dim}}
\newcommand{\Ker}{\mathrm{Ker}}
\newcommand{\dd}{\mathrm d}

\newcommand{\pint}{\langle\cdot,\cdot\rangle}
\newcommand{\anul}{{\mathrm o}}
\newcommand{\Img}{\mathrm{Im}}
\newcommand{\Id}{\mathrm{Id}}

\newcommand{\iCZ}{\mathfrak i_{\mathrm{CZ}}}
\newcommand{\iMaslov}{\mathfrak i_{\mathrm{M}}}
\newcommand{\Gr}{\mathrm{Gr}}

\newcommand{\perpB}{{\perp_{\!B}}}

\newcommand{\Ver}{\mathrm{Ver}}
\newcommand{\Hor}{\mathrm{Hor}}

\newcommand{\gR}{g_{\mathrm{R}}}
\newcommand{\llangle}{\langle\!\!\langle}
\newcommand{\rrangle}{\rangle\!\!\rangle}

\newcommand{\LLambda}{\mathbf\Lambda}

\renewcommand{\marginpar}[1]{\null}

\documentclass[twoside,final,a4paper,10pt]{amsart}

\usepackage{times}
\usepackage{dsfont}
\usepackage[all]{xy}
\usepackage{hyperref}

\numberwithin{equation}{section}

\title[Morse theory for Lorentzian closed geodesics]{On a Gromoll-Meyer type theorem\\ in globally hyperbolic stationary spacetimes}
\author[L. Biliotti]{Leonardo Biliotti}
\author[F. Mercuri]{Francesco Mercuri}
\author[P.\ Piccione]{Paolo Piccione}
\address{Departamento de Matem\'atica,
Universidade de S\~ao Paulo, \hfill\break\indent Rua do Mat\~ao
1010, CEP 05508-900, S\~ao Paulo, SP, Brazil}
\email{biliotti@dipmat.univpm.it, mercuri@ime.unicamp.br,\hfill\break\phantom{\textit{Email address:}}piccione@ime.usp.br}
\urladdr{http://www.ime.usp.br/\~{}piccione}
\thanks{F. M. and P. P. are partially sponsored by CNPq and Fapesp}


\date{January 23rd, 2007}

\begin{document}


\theoremstyle{plain}\newtheorem*{teon}{Theorem}
\theoremstyle{definition}\newtheorem*{defin*}{Definition}
\theoremstyle{plain}\newtheorem{teo}{Theorem}[section]
\theoremstyle{plain}\newtheorem{prop}[teo]{Proposition}
\theoremstyle{plain}\newtheorem{lem}[teo]{Lemma}
\theoremstyle{plain}\newtheorem{cor}[teo]{Corollary}
\theoremstyle{definition}\newtheorem{defin}[teo]{Definition}
\theoremstyle{remark}\newtheorem{rem}[teo]{Remark}
\theoremstyle{plain} \newtheorem{assum}[teo]{Assumption}
\swapnumbers
\theoremstyle{definition}\newtheorem{example}{Example}[section]
\theoremstyle{plain} \newtheorem*{acknowledgement}{Acknowledgements}
\theoremstyle{definition}\newtheorem*{notation}{Notation}


\begin{abstract}
Following the lines of the celebrated Riemannian result of Gromoll and Meyer \cite{GroMey2},
we use infinite dimensional equivariant Morse theory to establish the existence
of infinitely many \emph{geometrically distinct} closed geodesics in a class
of globally hyperbolic stationary Lorentzian manifolds.
\end{abstract}

\maketitle
\tableofcontents

\begin{section}{Introduction}
The question of existence of closed geodesics is one of the most classical theme of Riemannian
\marginpar{\textbf{Franco, vuoi\\ rivedere  l'introduzione?}}
geometry (see \cite{Kli}); spectacular contributions to the theory of global geometry have
been given in this area by very many authors, including Hadamard, Cartan, Poincar\'e, Birkhoff, Morse
and many others. Variational techniques for establishing the existence and the multiplicity of
closed geodesics have been developed and employed by many authors, including among others Ljusternik,
Schnirelman, Fet,  Klingenberg, Gromoll and Meyer.
Closed geodesics are critical points of the geodesic action functional
in the space of closed paths, and existence results may be obtained by applying global
variational techniques to this variational problem. In particular, Morse theory has been used
by Gromoll and Meyer (see \cite{GroMey2}) to establish the celebrated result on the existence of
infinitely many geometrically distinct closed geodesics in simply connected Riemannian manifolds, whose space
of free closed curves has unbounded rational Betti numbers.

As to the existence of closed geodesics in manifolds endowed with a non positive definite
metric, very few results are available in the literature, and basically nothing is known on
their multiplicity. An earlier result by Tipler
(see \cite{Tipler}) gives the existence of one closed timelike geodesic in compact Lorentzian
manifolds that admit a regular covering which has a compact Cauchy surface. More recently, Guediri
(see \cite{Gue1, Gue2}) has extended Tipler's result to the case that the Cauchy surface in
the covering is not necessarily compact. In this situation, a closed geodesic is proven to exist
in each free timelike homotopy class which is determined by a central deck transformation.
It is also proved in \cite{Gue1} that compact flat spacetimes contain a causal (i.e., nonspacelike)
closed geodesic, and in \cite{Gue2} the author proves that such spacetimes contain a closed timelike
geodesic if and only if the fundamental group of the underlying manifold contains a nontrivial timelike
translation.
The existence of closed timelike geodesic has been established also by Galloway in \cite{Gal1}, where
the author proves the existence of a longest closed timelike curve, which is necessarily a geodesic,
in each \emph{stable} free timelike homotopy class. Also non existence results for
nonspacelike geodesics are available, see \cite{Gal2, Gue3}.

All these results are based on the notion of Lorentzian distance function (see \cite[Ch.~4]{BEE}).
Recall that in Lorentzian geometry \emph{only} nonspacelike geodesics have length extremizing properties,
while for spacelike geodesics usual geometrical constructions (curve shortening methods) do not work.
The question of existence of closed geodesics of arbitrary causal character
has to be studied using the quadratic geodesic action functional in the Hilbert manifold of closed paths
of Sobolev regularity $H^1$; its critical points are typically saddle points.
In the Lorentzian (or semi-Riemannian) case, the variational theory associated to the study of the critical points
of this quadratic functional is complicated by the fact that, unlike the Riemannian counterpart, the
condition of Palais and Smale is never satisfied. Moreover, this functional is not bounded from
below, and its critical points always have infinite Morse index.
In \cite{AntSam} the authors use an approximation scheme in the theory of Ljusternik and Schnirel\-man
to determine the existence of a critical point of the geodesic action functional in the space of closed $H^1$
curves in a class of product Lorentzian manifolds, whose metric is of splitting type. Such critical
point corresponds to a spacelike closed geodesic; in this situation, thanks to the result of Galloway
\cite{Gal1}, one has two geometrically distinct closed geodesic, one is timelike and the other is spacelike.
A. Masiello has proved the existence of one (spacelike) closed geodesic in standard stationary
Lorentzian manifolds $M=M_0\times\R$ whose spatial component $M_0$ is compact.
More recently (see \cite{CapMasPic}), using variational methods the authors have established the existence of a
closed geodesic in each free homotopy class corresponding to an element of the fundamental group having
finite conjugacy class, in the case of \emph{static} compact Lorentzian manifold. In \cite{San06},
the author shows that one closed timelike geodesic exists in compact Lorentzian manifolds
that are conformally static, provided that that the group of deck transformations of some globally
hyperbolic regular covering of the manifold admits a finite conjugacy class containing a closed timelike curve.

In this paper, we develop a Morse theory for closed geodesics in a class of stationary Lorentzian manifolds,
obtaining a result of existence of infinitely many distinct closed geodesics analogous to the corresponding
result of Gromoll and Meyer in the Riemannian case. More precisely, the result will hold
for simply connected stationary manifolds whose free loop space has
unbounded Betti numbers (relatively to any coefficient field), and that admit a compact Cauchy surface.\footnote{%
Recall that any two
Cauchy surfaces of a globally hyperbolic spacetime are homeomorphic.}
Let us recall briefly the essential ingredients required  in Gromoll and Meyer's theory.
One considers the geodesic action functional $f$ on the Hilbert manifold $\Lambda M$ of
all closed paths of Sobolev class $H^1$ on a compact and simply connected Riemannian manifold
$(M,g)$; this functional is bounded from below, it
satisfies the Palais--Smale condition and its critical points are exactly the closed geodesics.
The compact group $\mathrm O(2)$ acts equivariantly on $\Lambda M$ via the operation of
$\mathrm O(2)$ on the parameter circle $\mathds S^1$; the orbits of this action are smooth
(compact) submanifolds of $\Lambda M$. In particular, the critical points of $f$ are never isolated;
nevertheless, using a generalized Morse Lemma for possibly degenerate isolated critical point (see \cite{GroMey1}),
generalized Morse inequalities can be applied to obtain estimates on the number
of critical orbits, provided that these orbits are isolated.  Finally, one has to distinguish between
critical orbits that correspond to iterates of the same closed geodesic. This is done
using an iteration formula for the Morse index (and the nullity) of the $n$-fold covering of a given closed
geodesic, which is obtained from a celebrated result due to Bott (see \cite{Bo4}).
Using a Morse index theorem for closed geodesics,
the Morse index of the $n$-th iterated of a closed geodesic is proven to be either bounded,
or to have linear growth in $n$. Using this fact one proves that if $(M,g)$ has only a finite
number of geometrically distinct closed geodesics, then the rational Betti numbers of $\Lambda M$ must
form a bounded sequence. Restriction to the case of a field of characteristic zero was used by the authors
to prove an estimate on the dimension of relative homology spaces of certain fiber bundles;
an elementary argument based on the Mayer--Vietoris sequence, discussed in Appendix~\ref{sec:homfiberS1},
shows that such restriction is not necessary.\footnote{%
Extension of the Gromoll and Meyer result to non zero characteristic seems to have been
established in the subsequent literature.}
Thus, if $\Lambda M$ has an unbounded sequence of Betti numbers,
$(M,g)$ must contain infinitely many geometrically distinct closed geodesics.
It is known (see \cite{VigSul}), that the existence of an unbounded sequence of
rational Betti numbers of the free loop space of $M$ is equivalent to the fact that the rational
\marginpar{\textbf{Franco: citare metri\-che Finsler e lavori\\ relativi.}}
cohomology algebra of $M$ is generated by at least two elements. In particular,
if $M$ has the same homotopy type of the product of two simply connected
compact manifold, then $\Lambda M$ has unbounded rational Betti numbers.
Ziller \cite{Zil} has proved that any compact symmetric space of rank greater than $1$
has unbounded $\Z_2$-Betti numbers; McCleary and Ziller \cite{McCleZil}
have later proved that the same conclusion holds for compact, simply connected
homogeneous spaces which are not diffeomorphic to a symmetric space of rank one.

Several extensions
of the theory have been developed in the context of Riemannian manifolds
(see \cite{BalThoZil, BanHin, BanKli, GroHalVig, GroHal, GroTan, Tan}), Finsler manifolds (\cite{Mat})
and, recently, of Riemannian orbifolds (see \cite{GurHae}). Reference \cite{Tan2} is
a good survey paper on the classical results of Gromoll--Meyer type.

When passing to the case of Lorentzian metrics, none of the arguments above works.
First, the geodesic action functional $f$ is not bounded from below and it does
not satisfy the PS condition; besides, the Morse index of each critical point
is infinite.
In this paper we consider the case of
stationary Lorentzian manifolds that admit a complete timelike Killing vector field.
Timelike invariance of the metric tensor allows
to determine a smooth embedded submanifold $\mathcal N$ of $\Lambda M$ with the following
properties:
\begin{itemize}
\item $f\vert_{\mathcal N}$ has the same critical points of $f$;
\item $\mathcal N$ has the same homotopy type of $\Lambda M$;
\item $f$ is bounded from below and it satisfies the PS condition
on each connected component of $\mathcal N$;
\item each critical point of $f\vert_{\mathcal N}$ has finite Morse index;
\item if a critical point is degenerate for $f\vert_{\mathcal N}$, then it is also
degenerate for $f$.
\end{itemize}
The abelian group $G=\mathrm O(2)\times\R$ acts (isometrically) on $\mathcal N$,
and $f$ is $G$-invariant. The group $\mathrm O(2)$ acts on the parameter space
$\mathds S^1$ of the curves, and as in the Riemannian case, this action is not smooth,
but only continuous. Nevertheless, if $\gamma$ is a smooth curve, then the orbit
$\mathrm O(2)\gamma$ is smooth, and it is diffeomorphic to $\mathrm O(2)$ if $\gamma$
is not constant. In particular, critical orbits are always smooth.
The group $\R$ acts by translation along the flow lines of the timelike
Killing vector field; obviously, the actions of $\mathrm O(2)$ and of $\R$ commute.
In this situation, we define \emph{geometrically distinct} two closed geodesics
that belong to different $G$-orbits, and that cannot be obtained one from another by iteration.
The action of $\R$ is free,  the orbit space given by
the quotient $\widetilde{\mathcal N}=\mathcal N/\R$ is a smooth manifold
and $\mathcal N$ is diffeomorphic to the product $\widetilde{\mathcal N}\times\R$.
Thus, in order to study multiplicity of distinct closed geodesics, it suffices
to study geometrically distinct critical $\mathrm O(2)$-orbits for the constrained functional
$f\vert_{\widetilde{\mathcal N}}$. The central result of this paper, which gives
the existence of infinitely many distinct closed geodesics in a class of stationary Lorentzian
manifolds, is obtained applying equivariant Morse theory to this setup.
Essential tools for the development of the theory are
a calculation of the Morse index for each critical point of $f\vert_{\mathcal N}$,
and a formula that describes its growth under iterations. The Morse index is given
in terms of symplectic invariants of the geodesic, such as the Conley--Zehnder and
the Maslov index, and it is computed explicitly in Theorem~\ref{thm:MorseIndexTheorem},
which is a Morse index theorem for possibly degenerate closed Lorentzian geodesics.
This result is obtained by purely functional analytical techniques, proving a
preliminary result (Theorem~\ref{thm:magicbilin}) that gives a method for computing
the index of essentially positive symmetric bilinear forms, possibly degenerate, in
terms of restrictions to possibly degenerate subspaces.
We believe that this result has interest in its own, and that its applicability
should go beyond the purposes of the present paper.
Using this method, one reduces the computation of the Morse index for
periodic geodesics to the Morse index of the corresponding fixed endpoint
geodesic, avoiding the usual assumption of \emph{orientability} of
the closed geodesic (see \cite{MarPicTau}).
The Morse index theorem is given in terms of a symplectic invariant of the geodesic,
called the Maslov index; in order to estimate its growth by iteration, we use a
recent formula that gives an estimate on the growth of another symplectic
invariant, called the Conley--Zehnder index (Proposition~\ref{eq:itformCZ}).
For orientation preserving geodesics, the two indices are related by a simple formula,
involving a four-fold index, which is called the H\"ormander index (Proposition~\ref{thm:compind2}).
Using the growth formula for the Conley--Zehnder index and the (non trivial) fact
that the Morse index is, up to a bounded perturbation, nondecreasing by iteration
(Lemma~\ref{thm:nondecreasing}), we then obtain a superlinear estimate on the
growth of the Maslov index of an iterate of a closed geodesic (Proposition~\ref{thm:lineargrowth} and
Corollary~\ref{thm:corcrescitalinearearb}).
As to the nullity of an iterate, the result is totally analogous to the Riemannian case
using the linearized Poincar\'e map (Lemma~\ref{thm:nuliteration}).
This setup paves the path to an application of infinite dimensional equivariant Morse theory,
in the same spirit as Gromoll and Meyer's celebrated result, that gives
the existence of infinitely many critical points for the functional
$f\vert_{\widetilde{\mathcal N}}$.

We will now give a formal statement of the main result of the paper.
Let $(M,g)$ be a globally hyperbolic stationary Lorentzian manifold, and let
us assume that $M$ admits a \emph{complete} timelike Killing vector field
$\mathcal Y$.
Denote by $\mathcal F_t$, $t\in\R$, the flow of $\mathcal Y$; clearly, if
$\gamma$ is a (closed) geodesic in $M$, then also $\mathcal F_t\circ\gamma$ is a (closed)
geodesic for all $t\in\R$.

In order to state our main result, we need to give an appropriate notion of geometric
equivalence of closed geodesics.

\begin{defin*}
Given closed geodesics $\gamma_i:[a_i,b_i]\to M$, $i=1,2$, in a stationary Lorentzian manifold $(M,g)$,
we will say that they are \emph{geometrically distinct}, if there exists no $t\in\R$ such that the sets
$\mathcal F_t\circ\gamma_1\big([a_1,b_1]\big)$ and $\gamma_2\big([a_2,b_2]\big)$ coincide.
\end{defin*}

The main result of this paper is the following:

\begin{teon}
Let $(M,g)$ be a simply connected globally hyperbolic stationary Lorentzian manifold having
a complete timelike Killing vector field, and having a compact Cauchy surface.
Assume that the free loop space $\Lambda M$ has unbounded Betti numbers with respect
to some coefficient field. Then, there are infinitely many geometrically distinct non
trivial (i.e., non constant) closed geodesics in $M$.
\end{teon}
Note that, by causality, every closed geodesic in $(M,g)$ is \emph{spacelike}.
It should be observed here that, although the notion of geometric equivalence given
above depends on the choice of a complete timelike Killing vector field, the property of
existence of infinitely many geometrically distinct closed geodesics is intrinsic to
$(M,g)$  (see Remark~\ref{thm:remindependence}). It is also interesting to observe that
the statement of the Theorem admits a mild
generalization to the case of non simply connected manifolds (see Remark~\ref{thm:remnonsimplyconnected}).

\smallskip

The paper is organized as follows.
Section~\ref{sec:preliminaries} contains a few basic facts concerning bilinear forms
and their index; in Section~\ref{sec:indexesspos} we prove the main result concerning the
computation of the index of an essentially positive symmetric bilinear form on a
real Hilbert space (Theorem~\ref{thm:magicbilin}). In Section~\ref{sec:iteformMaslov}
we recall the notions of Conley--Zehnder index for a continuous symplectic path,
and of Maslov index for a continuous Lagrangian path. The central result is an inequality
(Corollary~\ref{thm:coriteformMaslovindex}) that provides an estimate on the growth
of the Maslov index of the iterate of a periodic solution of a Hamiltonian system.
The definition of such index depends on the choice of a periodic symplectic trivialization
along the solution of the Hamiltonian.
When applied to the case of periodic geodesics on a semi-Riemannian manifold $M$, under
a certain orientability assumption we have a canonical choice of a class of periodic
symplectic trivializations along the corresponding periodic solution of the geodesic
Hamiltonian in the cotangent bundle $TM^*$ (Subsection~\ref{sub:Maslovgeo}), and
we therefore obtain estimates on the growth of the Maslov index of orientation preserving
periodic geodesics. The results in Section~\ref{sec:iteformMaslov} are valid
for closed geodesics in arbitrary semi-Riemannian manifolds.
Section~\ref{sec:varprobPS} contains some material on the geodesic variational
problem in stationary Lorentzian manifolds and on the Palais--Smale condition of
the relative action functional. In Section~\ref{sec:indexthm} we prove a general
version of the Morse index theorem for closed geodesics in stationary Lorentzian
manifold, that holds in the general case of possibly degenerate and non orientation
preserving geodesics (Theorem~\ref{thm:MorseIndexTheorem}). This is obtained as
an application of Theorem~\ref{thm:magicbilin}, which reduces the periodic case
to the case of fixed endpoints geodesics. In subsection~\ref{sub:iterMorse} we
first show that the Morse index of an $N$-th iterate is nondecreasing on $N$, up to adding a bounded
sequence. Then, we use the Morse index theorem and the estimates on the growth of the Maslov index
to get an estimate on the growth of the Morse index under iteration. The central
result (Proposition~\ref{thm:lineargrowth}, Corollary~\ref{thm:corcrescitalinearearb}),
which provides an alternative approach to the iteration theory of Bott \cite{Bo4} also
for the Riemannian case, says that the index of an $N$-th iterate is either bounded or
it has linear growth in $N$, up to adding a bounded sequence. The nullity of an iterate
is studied in Subsection~\ref{sub:nullity}, and the result is totally analogous to the
Riemannian case. Finally, in Section~\ref{sec:proof}, we use equivariant Morse theory
for isolated critical $\mathrm O(2)$-orbits of the action functional $f$ in $\widetilde{\mathcal N}$
to prove our main result. We follow closely the original paper by Gromoll and Meyer,
but we take advantage of a more recent approach to equivariant Morse theory (\cite{Cha, Wang}),
that simplifies some of the constructions in \cite{GroMey2}. The local homological invariant
at an isolated critical orbit is defined as the relative homology of the critical
sublevel, modulo the sublevel minus the critical orbit. Using excision, this invariant
is computed as the relative homology of a fiber bundle over the circle modulo a subbundle;
these bundles can be described as associated bundles to the principal fiber bundle $\mathrm O(2)\to\mathrm O(2)/\Gamma$,
where $\Gamma\subset\mathrm{SO}(2)$ is the stabilizer of the orbit. One of the crucial
steps in Gromoll and Meyer construction is an estimate on the dimension of this relative
homology (see \eqref{eq:disugHkorbHkpto}); this estimate is proven in Appendix~\ref{sec:homfiberS1}
in the case of homology with coefficients in arbitrary fields using the Mayer--Vietoris sequence
in relative homology. This allows a slight generalization of the original result in \cite{GroMey2}, in
that no restriction is posed on the characteristic of the coefficient field.

In order to make the paper essentially self-contained, and to facilitate its reading,
we have opted to include in the present version of the manuscript the full statement
of some results already appearing in the literature and needed in our theory.
Quotations of the original authors and complete bibliographical references are
given for the proof of these results.

\end{section}

\begin{section}{Preliminaries about bilinear forms on normed vector spaces}
\label{sec:preliminaries}
\marginpar{Nella versione\\ da spedire alla\\ rivista questa\\ e la prossima sezione\\
si possono\\ compattare in\\ una sola sezione.}
In this section we collect a few elementary facts on bilinear forms on vector spaces.
All vector spaces considered in the entire text are assumed to be
real. Given a (normed) vector space $X$, we will denote by $X^*$ its (topological) dual;
let $B:X\times X\to\R$ be a bilinear
form on $X$. The $B$-orthogonal complement of a subspace $S\subset X$ is
defined by: \[S^\perpB=\big\{x\in X:\text{$B(x,y)=0$, for all $y\in
S$}\big\};\] the {\em kernel\/} of $B$ is defined by:
\[\Ker(B)=X^\perpB=\big\{x\in X:\text{$B(x,y)=0$, for all $y\in
X$}\big\}.\]
 Given a subspace $S\subset X$, then $\Ker(B\vert_{S\times S})=S\cap S^\perpB$; if $(S_i)_{i\in I}$ is a family of subspaces of $X$, then
$\Big(\sum_{i\in I}S_i\Big)^\perpB=\bigcap_{i\in I}S_i^\perpB$.
We say that $B$ is {\em nondegenerate\/} if
$\Ker(B)=\{0\}$. A subspace $S\subset X$ is called {\em
isotropic\/} if $B\vert_{S\times S}=0$. Assume now that $B$ is
symmetric. We say that $B$ is {\em positive definite\/} (resp.,
{\em positive semi-definite}) if $B(x,x)>0$ for all nonzero $x\in
X$ (resp., $B(x,x)\ge0$, for all $x\in X$). Similarly, we say that
$B$ is {\em negative definite\/} (resp., {\em negative
semi-definite}) if $B(x,x)<0$ for all nonzero $x\in X$ (resp.,
$B(x,x)\le0$, for all $x\in X$). A subspace $S\subset X$ is called
{\em positive\/} (resp., {\em negative}) for $B$ if
$B\vert_{S\times S}$ is positive definite (resp., negative
definite). The {\em index\/} of $B$ is the (possibly infinite)
natural number defined by:
\[\mathrm n_-(B)=\sup\big\{\Dim(W):\text{$W\subset X$ is a negative subspace for
$B$}\big\},\] and the {\em co-index\/} of $B$ is defined by:
\[\mathrm n_+(B)=\mathrm n_-(-B).\]
When not both $\mathrm n_-(B)$ and $\mathrm n_+(B)$ are infinite, the \emph{signature}  of $B$ is defined as the
difference $\segnatura(B)=\indice(B)-\coindice(B)$.

We collect in the following lemma all the elementary results concerning bilinear forms
on vector spaces that will be used throughout.
\begin{lem}\label{thm:KerBSSperp}
Let $X$ be a vector space and let  $B$ be a symmetric bilinear form on $X$.
\begin{enumerate}
\item\label{thm:ortogkernels} If $X=X_1\oplus X_2$ is a {\em $B$-orthogonal direct sum
decomposition}, i.e., $B(x_1,x_2)=0$ for all $x_1\in X_1$, $x_2\in
X_2$, then
$\Ker(B)=\Ker(B\vert_{X_1\times X_1})\oplus\Ker(B\vert_{X_2\times
X_2})$. In particular, $B$ is nondegenerate if and only if
$B\vert_{X_1\times X_1}$ and $B\vert_{X_2\times X_2}$ are both
nondegenerate.
\item\label{thm:BSSnondeg} If $S\subset X$ is a subspace with $X=\Ker(B)\oplus S$, then
$B\vert_{S\times S}$ is nondegenerate.
\item\label{thm:SSperp} If $X=S\oplus S^\perpB$ then
$B\vert_{S\times S}$ is nondegenerate. Conversely, if $S$ is
finite dimensional and $B\vert_{S\times S}$ is nondegenerate then
$X=S\oplus S^\perpB$.
\item\label{thm:indcodindmonot} Given any subspace $Y\subset X$, then:
\begin{gather*}
\mathrm n_+(B\vert_{Y\times Y})\le \mathrm n_+(B)\le \mathrm n_+(B\vert_{Y\times Y})+\mathrm{codim}_X(Y),\\
\mathrm n_-(B\vert_{Y\times
Y})\le \mathrm n_-(B)\le \mathrm n_-(B\vert_{Y\times Y})+\mathrm{codim}_X(Y).
\end{gather*}
\item\label{thm:calcind} Let $X=X_1\oplus X_2$ be a direct sum decomposition such that $B$
is positive definite on $X_1$ and negative semi-definite on $X_2$
(resp., negative definite on $X_1$ and positive semi-definite on
$X_2$). Then $\mathrm n_+(B)=\Dim(X_1)$ (resp., $\mathrm n_-(B)=\Dim(X_1)$).
\item\label{thm:twoposnegs} Assume that $X=X_1\oplus X_2$ is a $B$-orthogonal
direct sum decomposition such that $B$ is positive definite (resp.,
negative definite) on both $X_1$ and $X_2$. Then $B$ is positive
definite (resp., negative definite) on $X$. Similarly, if $B$ is
positive semi-definite (resp., negative semi-definite) on both
$X_1$ and $X_2$ then $B$ is positive semi-definite (resp.,
negative semi-definite) on $X$.
\item \label{thm:maximalneg}
Let $S\subset X$ be a maximal subspace on which $B$ is positive
definite (resp., negative definite). Then $B$ is negative
semi-definite (resp., positive semi-definite) on $S^\perpB$.
\item \label{thm:ortogindcoind}
Let  $X=X_1\oplus X_2$ be a $B$-orthogonal direct sum
decomposition. Then:
\begin{gather*} \mathrm n_+(B)=\mathrm n_+(B\vert_{X_1\times
X_1})+\mathrm n_+(B\vert_{X_2\times X_2}),\\
\mathrm n_-(B)=\mathrm n_-(B\vert_{X_1\times X_1})+\mathrm n_-(B\vert_{X_2\times X_2}).
\end{gather*}
\item \label{thm:BSSindcoind}
Let $S\subset X$, $N\subset\Ker(B)$ be subspaces with $X=N\oplus
S$. Then:
\[\mathrm n_+(B\vert_{S\times S})=\mathrm n_+(B),\quad \mathrm n_-(B\vert_{S\times
S})=\mathrm n_-(B).\]
\item\label{thm:Bbar}
Let  $Y$ be a vector space and let $q:X\to Y$ be surjective linear map with
$\Ker(q)\subset\Ker(B)$. Then there exists a unique map $\overline
B:Y\times Y\to\R$ such that:
\[\overline B\big(q(x_1),q(x_2)\big)=B(x_1,x_2),\quad\text{for all $x_1,x_2\in X$};\]
the map $\overline B$ is a symmetric bilinear form on $Y$.
Moreover:
\begin{gather} \label{eq:kerbarB}
\Ker(\overline B)=q\big(\Ker(B)\big)\cong\Ker(B)/\Ker(q),\\
\label{eq:n+-barB} \mathrm n_+(\overline B)=\mathrm n_+(B),\quad \mathrm n_-(\overline
B)=\mathrm n_-(B).
\end{gather}
In particular, if $\Ker(q)=\Ker(B)$ then $\overline B$ is
nondegenerate.
\item \label{thm:n+n-deg}
The following formula holds:
\begin{equation}\label{eq:Dim+-deg}
\Dim(X)=\mathrm n_+(B)+\mathrm n_-(B)+\Dim\big(\Ker(B)\big).
\end{equation}
\item \label{thm:alittlebetter}
Let   $X_1,X_2\subset X$ be {\em $B$-orthogonal subspaces}, i.e.,
$B(x_1,x_2)=0$ for all $x_1\in X_1$, $x_2\in X_2$. If $X=X_1+X_2$
(not necessarily a direct sum) then:
\begin{gather*} \mathrm n_+(B)=\mathrm n_+(B\vert_{X_1\times
X_1})+\mathrm n_+(B\vert_{X_2\times X_2}),\\
\mathrm n_-(B)=\mathrm n_-(B\vert_{X_1\times X_1})+\mathrm n_-(B\vert_{X_2\times X_2}).
\end{gather*}
\item \label{thm:seminondeg}
Let $B$ be a symmetric bilinear form on a vector space $X$. If $B$
is positive semi-definite or negative semi-definite then:
\[\Ker(B)=\big\{x\in X:B(x,x)=0\big\}.\]
In particular if $B$ is positive semi-definite (resp., negative
semi-definite) and nondegenerate then $B$ is positive definite
(resp., negative definite).
\item\label{thm:upbdisotrop}
If $B$ is nondegenerate and symmetric, and  $S\subset X$ is an isotropic subspace, then:
\[\Dim(S)\le \mathrm n_-(B),\quad\Dim(S)\le \mathrm n_+(B).\qed\]
\end{enumerate}
\end{lem}

Let us now consider a (real) normed space $X$. If $T:X\to Y$ is a continuous
linear map between normed spaces then $T^*:Y^*\to X^*$ denotes the
transpose map defined by $T^*(\alpha)=\alpha\circ T$. If $S\subset
X$ is a subspace we denote by $S^\anul\subset X^*$ the {\em
annihilator\/} of $S$, i.e.:
\[S^\anul=\big\{\alpha\in X^*:\alpha\vert_S=0\big\}.\]

If $X$, $Y$ are Banach spaces and $T:X\to Y$ is a continuous
linear map then $\Ker(T^*)=\Img(T)^\anul$ and
$\Img(T^*)\subset\Ker(T)^\anul$. Moreover, if $\Img(T)$ is closed
in $Y$ then $\Img(T^*)=\Ker(T)^\anul$. Given a closed subspace
$S\subset X$, denote by $q:X\to X/S$ the quotient map; then
$q^*:(X/S)^*\to X^*$ is injective and its image equals $S^\anul$.
Moreover, if $X$ is reflexive, by identifying $X$ with $X^{**}$ in
the usual way, the bi-annihilator $(S^\anul)^\anul$ equals the
closure of $S$.

\begin{rem}\label{thm:saturated}
If $X$, $Y$ are sets and $f:X\to Y$ is a map then a subset
$S\subset X$ is called {\em saturated for $f$\/} (or {\em
$f$-saturated}) if $x_1\in S$, $x_2\in X$ and $f(x_1)=f(x_2)$
imply $x_2\in S$. If $X$, $Y$ are vector spaces, $f$ is linear and
$S\subset X$ is a subspace then $S$ is $f$-saturated if and only
if $\Ker(f)\subset S$. Observe that if $X$, $Y$ are Banach spaces
and $f:X\to Y$ is a surjective continuous linear map then, by the
open mapping theorem, $f$ is a quotient map in the topological
sense; hence a saturated subset $S\subset X$ is open (resp.,
closed) in $X$ if and only if $f(S)$ is open (resp., closed) in
$Y$. Similarly, a subset $U\subset Y$ is open (resp., closed) in
$Y$ if and only if $f^{-1}(U)$ is open (resp., closed) in $X$.
\end{rem}

\begin{lem}\label{thm:finitecodplussomething}
Let $X$ be a Banach space and let $Y\subset X$ be a finite
co-dimensional closed subspace. If $Z\subset X$ is a subspace with
$Y\subset Z$ then $Z$ is also closed in $X$.
\end{lem}
\begin{proof}
If $q:X\to X/Y$ denotes the quotient map then $q(Z)$ is closed in
$X/Y$ because $X/Y$ is finite dimensional. But since $Y\subset Z$
we have that $Z$ is $q$-saturated and thus $Z$ is closed in $X$
(see Remark~\ref{thm:saturated}).
\end{proof}

\begin{lem}\label{thm:closedplusfinite}
Let $X$ be a Banach space, $Y_1\subset X$ be a closed subspace and
$Y_2\subset X$ be a finite dimensional subspace. Then $Y_1+Y_2$ is
closed in $X$.
\end{lem}
\begin{proof}
If $q:X\to X/Y_1$ denotes the quotient map then $q(Y_2)$ is closed
in $X/Y_1$ because it is finite dimensional. But then
$Y_1+Y_2=q^{-1}\big(q(Y_2)\big)$ is also closed in $X$.
\end{proof}

Finally, let us prove the following:

\begin{lem}\label{thm:Sperpperp}
Let $X$ be a reflexive Banach space, $B$ be a continuous symmetric
bilinear form on $X$ and $S\subset X$ be a subspace. If
$\Img(B)+S^\anul$ is closed in $X^*$ then the bi-orthogonal
complement of $S$ is given by:
\[(S^\perpB)^\perpB=\overline S+\Ker(B).\]
\end{lem}
\begin{proof}
Clearly, $S^\perpB=B^{-1}(S^\anul)$
and:
\[(S^\perpB)^\perpB=B^{-1}\big[\big(B^{-1}(S^\anul)\big)^\anul\big].\]
Denote by $q:X^*\to X^*/S^\anul$ the quotient map (observe that an
annihilator is always closed). Obviously
$B^{-1}(S^\anul)=\Ker(q\circ B)$. Since $\Img(B)+S^\anul$ is
$q$-saturated and closed in $X^*$, we have that $\Img(q\circ B)$
is closed in $X^*/S^\anul$ (recall Remark~\ref{thm:saturated});
then,  we have:
\[\big(B^{-1}(S^\anul)\big)^\anul=\Ker(q\circ B)^\anul=\Img\big((q\circ B)^*\big).\]
Now, using
the fact that $(q\circ B)^*=B^*\circ
q^*$, we obtain:
\[\Img\big((q\circ B)^*\big)=B(\overline S).\]
Finally:
\[(S^\perpB)^\perpB=B^{-1}\big[\big(B^{-1}(S^\anul)\big)^\anul\big]=B^{-1}\big(B(\overline
S)\big)=\overline S+\Ker(B).\qedhere\]
\end{proof}
\end{section}

\begin{section}{On the index of essentially positive bilinear forms}
\label{sec:indexesspos}
In this section we will discuss some functional analytical
preliminaries needed for the index theorem. The central result is Theorem~\ref{thm:magicbilin},
that gives a result concerning the computation of the index of
symmetric bilinear forms, possibly degenerate, using restrictions to
possibly degenerate subspaces.
Throughout this section we will always identify continuous
bilinear forms $B:X\times X\to\R$ on a normed space $X$ with the
continuous linear map $B:X\to X^*$ given by $x\mapsto B(x,\cdot)$.

Let $(X,\pint)$ be a Hilbert space and let $B:X\times X\to\R$ be a
continuous bilinear form. We say that a continuous linear operator
$T:X\to X$ {\em represents\/} $B$ with respect to $\pint$ if:
\[B(x,y)=\big\langle T(x),y\big\rangle,\]
for all $x,y\in X$.
Observe that if $i:X\to X^*$ denotes the isometry $x\mapsto\langle
x,\cdot\rangle$ given by Riesz representation theorem, then $T$
represents $B$ if and only if $i\circ T=B$. In particular,
$\Ker(T)=\Ker(B)$.
A continuous bilinear form $B:X\times X\to\R$ on a Banach space
$X$ is called {\em strongly nondegenerate\/} if the linear map
$B:X\to X^*$ is an isomorphism.
Obviously if $B$ is strongly nondegenerate then $B$ is
nondegenerate. The converse holds if we know that the linear map
$B:X\to X^*$ is a Fredholm operator of index zero (for instance, a
compact perturbation of an isomorphism).

Let $(X,\pint)$ be a Hilbert space. A continuous linear operator
$P:X\to X$ is called {\em positive\/} if the bilinear form
$\langle P\cdot,\cdot\rangle$ represented by $P$ is symmetric and
positive semi-definite. We recall the following standard result:
\begin{lem}\label{thm:carprodint}
Let $(X,\pint)$ be a Hilbert space and let $P:X\to X$ be a
continuous linear operator. Then the bilinear form $\langle
P\cdot,\cdot\rangle$ represented by $P$ is an inner product on $X$
that defines the same topology as $\pint$ if and only if $P$ is a
positive isomorphism of $X$.\qed
\end{lem}

One proves easily the following:
\begin{lem}\label{thm:BSSstrongly}
Let $X$ be a Banach space, $B$ be a continuous bilinear form on
$X$ and $S\subset X$ be a closed subspace. If $X=S\oplus S^\perpB$
then $B\vert_{S\times S}$ is nondegenerate. Conversely, if
$B\vert_{S\times S}$ is strongly nondegenerate then $X=S\oplus
S^\perpB$.\qed
\end{lem}

\begin{defin}\label{thm:defRCPPI}
Let $(X,\pint)$ be a Hilbert space and let $B$ be a continuous
symmetric bilinear form on $X$. We say that $B$ is {\em
essentially positive\/}
if the operator $T:X\to X$ that represents $B$ is of
the form $T=P+K$, with $P:X\to X$ a positive isomorphism and
$K:X\to X$ a (symmetric) compact operator.
\end{defin}

The following summarizes the main properties of essentially positive bilinear forms.
\begin{lem}\label{thm:propRCPPI}
Let $X$ and $Y$ be (real) Hilbert spaces. The following results hold:
\begin{enumerate}
\item\label{itm:RCPPI1} if
$B:X\times X\to\R$ is a continuous symmetric bilinear form then
$B$ is essentially positive if and only if there exists an inner product
$\pint_1$ on $X$ and a compact operator $K:X\to X$ such that
$B=\pint_1+\langle K\cdot,\cdot\rangle$ and such that $\pint_1$
defines the same topology on $X$ as $\pint$.

\item\label{itm:RCPPI2} If  $T:X\to Y$
is a continuous isomorphism and $B:Y\times Y\to\R$ is an essentially positive
symmetric bilinear form, then the {\em pull-back\/}
$T^*B=B(T\cdot,T\cdot):X\times X\to\R$ is also essentially positive.

\item\label{itm:RCPPI3} If
$B:X\times X\to\R$ is an essentially positive
symmetric bilinear form, then $B$ is also essentially positive with respect to any other inner
product $\pint_1$ on $X$ which defines the same topology as
$\pint$.

\item\label{itm:RCPPI4} If
$B:X\times X\to\R$ is an essentially positive
symmetric bilinear form, then for every closed subspace $S\subset X$, the
restriction $B\vert_{S\times S}$ is again essentially positive.

\item\label{itm:RCPPI4half} if  $B$ is an essentially positive
symmetric bilinear form on $X$,
$q:X\to Y$ is a surjective continuous linear map with
$\Ker(q)\subset\Ker(B)$ and $\overline B$ is defined as in
part~\eqref{thm:Bbar} in Lemma~\ref{thm:KerBSSperp} then $\overline B$ is also (bilinear,
symmetric, continuous and) essentially positive.

\item\label{itm:RCPPI5} If
$B:X\times X\to\R$ is an essentially positive
symmetric bilinear form, then there exists an inner product $\pint_1$ on $X$ that
defines the same topology as $\pint$ and such that $B$ is
represented with respect to $\pint_1$ by an operator of the form
identity plus compact.

\item\label{itm:RCPPI6} If
$B:X\times X\to\R$ is an essentially positive
symmetric bilinear form, then $\Ker(B)$ is finite dimensional and $\mathrm n_-(B)$ is
finite.
\end{enumerate}
\end{lem}
\begin{proof}\
\begin{enumerate}
\item It is a straightforward consequence of Lemma~\ref{thm:carprodint}.

\item If $B$ is represented by $P+K$ with $P:Y\to Y$ a positive
isomorphism and $K:Y\to Y$ compact then the pull-back $T^*B$ is
represented by $T^*(P+K)T$, where $T^*$ is identified with a
linear map from $Y$ to $X$, using Riesz representation theorem.
Now it is easy to see that $T^*PT$ is a positive isomorphism of
$X$ and that $T^*KT$ is a compact operator on $X$.

\item Follows from item~\eqref{itm:RCPPI2}, setting
$T=\Id:(X,\pint_1)\to(X,\pint)$.

\item Write $B=\pint_1+\langle K\cdot,\cdot\rangle$, as in
item~\eqref{itm:RCPPI1}. Denote by $\pi:X\to S$ the orthogonal
projection onto $S$ and set $K'=\pi\circ K\vert_S$. Then
$\pint_1\vert_{S\times S}$ is an inner product on $S$ that defines
the same topology as $\pint\vert_{S\times S}$ and $B\vert_{S\times
S}-\pint_1\vert_{S\times S}$ is represented by $K':S\to S$, which
is a compact operator.

\item Let $S$ be the orthogonal complement of $\Ker(q)$ with
respect to the Hilbert space inner product of $X$. Then
$q\vert_S:S\to Y$ is an isomorphism and $\overline B$ is the
pull-back of $B\vert_{S\times S}$ by $(q\vert_S)^{-1}$. The
conclusion follows from items~\eqref{itm:RCPPI2} and
\eqref{itm:RCPPI4}.

\item Take $\pint_1$ to be the inner product whose existence is given in item~\eqref{itm:RCPPI1}.

\item By item~\eqref{itm:RCPPI5}, we can choose the inner product
on $X$ such that $B$ is represented by $\Id+K$, with $K$ compact.
Then $\Ker(B)=\Ker(\Id+K)$ is finite-dimensional, because $\Id+K$
is a Fredholm operator. If $X_\lambda\subset X$ denotes the
$\lambda$-eigenspace of $\Id+K$ then by the spectral theorem for
compact self-adjoint operators on real Hilbert spaces, $X$ is the
closure of the algebraic ($\pint$-orthogonal) direct sum
$\bigoplus_{\lambda\in\sigma(\Id+K)}X_\lambda$, where
$\sigma(\Id+K)\subset\R$ denotes the spectrum of the operator
$\Id+K$; moreover, for $\lambda\ne1$, $X_\lambda$ is
finite-dimensional. Set:
\[X_-=\bigoplus_{\substack{\lambda\in\sigma(\Id+K)\\\lambda<0}}X_\lambda,\quad
X_+=\overline{\bigoplus_{\substack{\lambda\in\sigma(\Id+K)\\\lambda\ge0}}X_\lambda};\]
then it is easy to see that $B$ is negative definite on $X_-$,
positive semi-definite on $X_+$ and that $X=X_-\oplus X_+$.
Finally, from  part~\eqref{thm:calcind} in Lemma~\ref{thm:KerBSSperp}, $\mathrm n_-(B)=\Dim(X_-)<+\infty$,
recalling that $\big\{\lambda\in\sigma(\Id+K):\lambda<0\big\}$ is
finite, because $\sigma(\Id+K)$ is bounded and has $1$ as its only
limit point.\qedhere
\end{enumerate}
\end{proof}

\begin{rem}\label{thm:remRCPPInondeg}
If a continuous symmetric bilinear form $B$ on a Hilbert space $X$
is essentially positive then $B$ is nondegenerate if and only if $B$ is strongly
nondegenerate. Namely, $B$ is represented by a Fredholm operator
of index zero.
\end{rem}

\begin{rem}\label{thm:remperperRCPPI}
If $B$ is a continuous symmetric bilinear form on a Hilbert space
$X$ which is essentially positive and if $S\subset X$ is a subspace then:
\[(S^\perpB)^\perpB=\overline S+\Ker(B);\]
namely, we only have to check the hypotheses of
Lemma~\ref{thm:Sperpperp}. Obviously $X$ is reflexive, being a
Hilbert space. Moreover, $\Img(B)$ is closed and finite
co-dimensional in $X^*$, because $B$ is a Fredholm operator; thus,
$\Img(B)+S^\anul$ is closed in $X^*$, by
Lemma~\ref{thm:finitecodplussomething}.
\end{rem}

\begin{lem}\label{thm:WplusWperpclosed}
Let $X$ be a Hilbert space and let $B$ be a continuous symmetric
bilinear form on $X$ which is essentially positive. If $W\subset X$ is a closed
subspace then $W+W^\perpB$ is also closed in $X$.
\end{lem}
\begin{proof}
We can choose the inner product $\pint$ on $X$ such that $B$ is
represented by an operator of the form $T=\Id+K$, with $K$ compact
(see Lemma~\ref{thm:propRCPPI}, item~\eqref{itm:RCPPI5}). If $W'$
denotes the orthogonal complement of $W$ with respect to $\pint$
then $W^\perpB=T^{-1}(W')$. We then have to show that
$W+T^{-1}(W')$ is closed in $X$. Since $T$ is a Fredholm operator,
its image is closed in $X$ and so $T:X\to\Img(T)$ is a surjective
continuous linear operator between Banach spaces. We have
$\Ker(T)\subset T^{-1}(W')\subset W+T^{-1}(W')$, so that
$W+T^{-1}(W')$ is $T$-saturated (see Remark~\ref{thm:saturated});
thus $W+T^{-1}(W')$ is closed in $X$ if and only if
$T\big[W+T^{-1}(W')]=T(W)+\big(W'\cap\Img(T)\big)$ is closed in
$\Img(T)$. But:
\[T(W)+\big(W'\cap\Img(T)\big)=\big(T(W)+W'\big)\cap\Img(T),\]
and
therefore the proof will be completed once we show that $T(W)+W'$
is closed in $X$. We have:
\[T(W)+W'=\big\{x+y+K(x):x\in W,\ y\in
W'\big\}=\Img(\Id+K\circ\pi),\] where $\pi$ denotes the orthogonal
projection (with respect to $\pint$) onto $W$. Since $K\circ\pi$
is compact, $\Id+K\circ\pi$ is a Fredholm operator and hence its
image is closed in $X$.
\end{proof}

\begin{lem}\label{thm:ZcaseBnondeg}
Let $X$ be a Hilbert space and let $B$ be a nondegenerate
continuous symmetric bilinear form on $X$ which is essentially positive. If
$Z\subset X$ is an isotropic subspace then:
\begin{equation}\label{eq:n-BZ}
\mathrm n_-(B)=\mathrm n_-(B\vert_{Z^\perpB\times Z^\perpB})+\Dim(Z),
\end{equation}
all the terms on the equality above being finite
natural numbers.
\end{lem}
\begin{proof}
By~\eqref{itm:RCPPI6} in Lemma~\ref{thm:propRCPPI},
$\mathrm n_-(B)<+\infty$ and thus $\Dim(Z)<+\infty$, by
part~\eqref{thm:upbdisotrop} in Lemma~\ref{thm:KerBSSperp}. This proves that all terms on the
equality \eqref{eq:n-BZ} are finite natural numbers. Since $Z$ is
isotropic, we have $Z\subset Z^\perpB$ and thus we can find a
closed subspace $U\subset Z^\perpB$ with $Z^\perpB=Z\oplus U$ (for
instance, take $U$ to be the orthogonal complement of $Z$ in
$Z^\perpB$ with respect to any Hilbert space inner product). By
Lemma~\ref{thm:KerBSSperp}, we have:
\[\Ker(B\vert_{Z^\perpB\times
Z^\perpB})=Z^\perpB\cap(Z^\perpB)^\perpB.\] Now
Remark~\ref{thm:remperperRCPPI} implies that $(Z^\perpB)^\perpB=Z$.
We have thus proven that:
\[\Ker(B\vert_{Z^\perpB\times Z^\perpB})=Z,\]
and from part~\eqref{thm:BSSnondeg} in Lemma~\ref{thm:KerBSSperp} we obtain that
$B\vert_{U\times U}$ is nondegenerate. By item~\eqref{itm:RCPPI4}
of Lemma~\ref{thm:propRCPPI} and by
Remark~\ref{thm:remRCPPInondeg}, we obtain that $B\vert_{U\times
U}$ is actually strongly nondegenerate; thus, by
Lemma~\ref{thm:BSSstrongly}:
\[X=U\oplus U^\perpB.\]
From part~\eqref{thm:ortogindcoind} in Lemma~\ref{thm:KerBSSperp}, we then obtain:
\[\mathrm n_-(B)=\mathrm n_-(B\vert_{U\times U})+\mathrm n_-(B\vert_{U^\perpB\times
U^\perpB});\] again using part~\eqref{thm:ortogindcoind} in Lemma~\ref{thm:KerBSSperp}, the
$B$-orthogonal decomposition $Z^\perpB=Z\oplus U$ and the fact that
$Z$ is isotropic implies that:
\[\mathrm n_-(B\vert_{U\times U})=\mathrm n_-(B\vert_{Z^\perpB\times Z^\perpB}).\]
To complete the proof, it suffices to show that:
\[\mathrm n_-(B\vert_{U^\perpB\times U^\perpB})=\Dim(Z).\]
First, we claim that $\Dim(U^\perpB)=2\,\Dim(Z)$. Namely, since
$X=U\oplus U^\perpB$, the dimension of $U^\perpB$ equals the
co-dimension of $U$ in $X$. We have: \[U\subset Z^\perpB\subset
X;\] since $Z^\perpB=Z\oplus U$, the co-dimension of $U$ in
$Z^\perpB$ equals the dimension of $Z$. Since $B:X\to X^*$ is an
isomorphism (see Remark~\ref{thm:remRCPPInondeg}) and
$Z^\perpB=B^{-1}(Z^\anul)$, $B$
induces an isomorphism:
\[X/Z^\perpB\xrightarrow{\;\;\scriptscriptstyle\cong\;\;}X^*/Z^\anul;\]
moreover, $X^*/Z^\anul\cong Z^*\cong Z$. Thus the co-dimension of $Z^\perpB$
in $X$ is equal to the dimension of $Z$, which proves that
$\Dim(U^\perpB)=2\,\Dim(Z)$. To complete the proof of the lemma,
observe that $B$ is nondegenerate on $U^\perpB$ by
part~\eqref{thm:ortogkernels} in Lemma~\ref{thm:KerBSSperp}, and thus, by
part~\eqref{thm:n+n-deg} in Lemma~\ref{thm:KerBSSperp}:
\[\mathrm n_+(B\vert_{U^\perpB\times U^\perpB})+\mathrm n_-(B\vert_{U^\perpB\times
U^\perpB})=\Dim(U^\perpB)=2\,\Dim(Z),\] and by
part~\eqref{thm:upbdisotrop} in Lemma~\ref{thm:KerBSSperp}:
\[\mathrm n_+(B\vert_{U^\perpB\times U^\perpB})\ge\Dim(Z),\quad \mathrm n_-(B\vert_{U^\perpB\times
U^\perpB})\ge\Dim(Z).\] This proves that both
$\mathrm n_+(B\vert_{U^\perpB\times U^\perpB})$ and
$\mathrm n_-(B\vert_{U^\perpB\times U^\perpB})$ are equal to $\Dim(Z)$.
\end{proof}

\begin{lem}\label{thm:WcaseBnondeg}
Let $X$ be a Hilbert space and let $B$ be a nondegenerate
continuous symmetric bilinear form on $X$ which is essentially positive. If
$W\subset X$ is a closed subspace then:
\[\mathrm n_-(B)=\mathrm n_-(B\vert_{W\times W})+\mathrm n_-(B\vert_{W^\perpB\times
W^\perpB})+\Dim(W\cap W^\perpB),\] all the terms on the equality
above being finite natural numbers.
\end{lem}
\begin{proof}
Obviously $Z=W\cap W^\perpB$ is an isotropic subspace, and  we
can apply Lemma~\ref{thm:ZcaseBnondeg} to obtain:
\[\mathrm n_-(B)=\mathrm n_-(B\vert_{Z^\perpB\times Z^\perpB})+\Dim(W\cap W^\perpB).\]
The conclusion will follow from part~\eqref{thm:alittlebetter} in Lemma~\ref{thm:KerBSSperp},
once we show that $Z^\perpB=W+W^\perpB$. Using
Lemma~\ref{thm:KerBSSperp} and Remark~\ref{thm:remperperRCPPI}, we
compute:
\[(W+W^\perpB)^\perpB=W^\perpB\cap(W^\perpB)^\perpB=W^\perpB\cap W=Z.\]
Now using Lemma~\ref{thm:WplusWperpclosed} and
Remark~\ref{thm:remperperRCPPI} we obtain:
\[Z^\perpB=\big[(W+W^\perpB)^\perpB\big]^\perpB=W+W^\perpB.\qedhere\]
\end{proof}

Finally, the central result we aimed at:
\begin{teo}\label{thm:magicbilin}
Let $X$ be a Hilbert space and let $B$ be a continuous symmetric
essentially positive bilinear form on $X$. If $W\subset X$ is a
closed subspace and $S$ denotes the $B$-orthogonal space to $W$,
then:
\[\mathrm n_-(B)=\mathrm n_-\big(B\vert_{W\times W}\big)+\mathrm n_-\big(B\vert_{S\times S}\big)+\Dim(W\cap S)-\Dim\big(W\cap\Ker(B)\big),\] all
the terms on the equality above being finite natural numbers.
\end{teo}
\begin{proof}
Set $N=\Ker(B)$, $Y=X/N$ and denote by $q:X\to Y$ the quotient
map. Define $\overline B$ as in part~\eqref{thm:Bbar} in Lemma~\ref{thm:KerBSSperp}; then
$\overline B$ is a nondegenerate continuous symmetric bilinear
form on $Y$ and $\overline B$ is essentially positive, by
item~\eqref{itm:RCPPI4half} in Lemma~\ref{thm:propRCPPI}. We will
apply Lemma~\ref{thm:WcaseBnondeg} to $\overline B$ and to the
subspace $q(W)$ of $Y$; we first check that $q(W)$ is closed in
$Y$. By Remark~\ref{thm:saturated}, it suffices to observe that
$q^{-1}\big(q(W)\big)=W+N$ is closed in $X$ and this follows from
Lemma~\ref{thm:closedplusfinite} (recall that $\Dim(N)<+\infty$,
by item~\eqref{itm:RCPPI6} in Lemma~\ref{thm:propRCPPI}). Now:
\[\mathrm n_-(\overline B)=\mathrm n_-(\overline B\vert_{q(W)\times
q(W)})+\mathrm n_-(\overline B\vert_{q(W)^\perpB\times
q(W)^\perpB})+\Dim\big(q(W)\cap q(W)^\perpB\big).\] It is
straightforward to verify that:
\[q(W)^\perpB=q(W^\perpB).\]
Now using part~\eqref{thm:Bbar} in Lemma~\ref{thm:KerBSSperp} and considering the surjective
linear maps $q$, $q\vert_W:W\to q(W)$ and
$q\vert_{W^\perpB}:W^\perpB\to q(W)^\perpB$, we obtain:
\begin{gather*}
\mathrm n_-(\overline B)=\mathrm n_-(B),\\
\mathrm n_-(\overline B\vert_{q(W)\times q(W)})=\mathrm n_-(B\vert_{W\times
W}),\\
\mathrm n_-(\overline B\vert_{q(W)^\perpB\times
q(W)^\perpB})=\mathrm n_-(B\vert_{W^\perpB\times W^\perpB}).
\end{gather*}
To complete the proof, we have to show that:
\[\Dim\big(q(W)\cap q(W^\perpB)\big)=\Dim(W\cap W^\perpB)-\Dim(W\cap
N).\] Keeping in mind that $N\subset W^\perpB$, we compute:
\[q^{-1}\big(q(W)\cap q(W^\perpB)\big)=(W+N)\cap W^\perpB=(W\cap
W^\perpB)+N,\] so that $q(W)\cap q(W^\perpB)=q\big((W\cap
W^\perpB)+N\big)\cong\big[(W\cap W^\perpB)+N\big]/N$. Then:
\begin{multline*}
\Dim\big(q(W)\cap q(W^\perpB)\big)=\Dim\big[(W\cap
W^\perpB)+N\big]-\Dim(N)\\
=\Dim(W\cap W^\perpB)-\Dim(W\cap N). \end{multline*} This concludes
the proof.
\end{proof}

\end{section}

\begin{section}{On the Maslov index and iteration formulas}
\label{sec:iteformMaslov}
In this section we will prove an iteration formula for the Maslov index of a periodic solution
of a Hamiltonian system, using a similar formula proved in \cite{GosGosPic} for the
Conley--Zehnder index, and a formula relating the two indices via the H\"ormander index.
The reader should note that in the literature there are several definitions
Maslov index for a continuous Lagrangian path; these definitions
differ by a boundary term when the path has endpoints in the Maslov cycle.
In Robbin and Salamon \cite{RobSal}, the Maslov index is a half integer, obtained as half
of the variation of the signature function of certain bilinear forms, whereas in our case
we replace half of the signature with the \emph{extended coindex} (i.e., index plus nullity),
see formula~\eqref{eq:muL0tilde}. Obviously, the two definitions are totally equivalent; however,
the reader should observe that, using our definition, the Maslov index takes integer values, but
it fails to have the property that, when one changes the sign of the symplectic form
the absolute value of the Maslov index remains constant.
The reader should also be aware of the fact that the definition of Maslov index for a semi-Riemannian
geodesic adopted here differs slightly from previous definitions, originated from Helfer \cite{Hel},
in that here we also consider the contribution given by the initial endpoint, which is always conjugate.
The results in this section (more specifically, in Subsection~\ref{sub:Maslovgeo}) are valid
in any semi-Riemannian manifold $(M,g)$.

\subsection{Maslov and Conley--Zehnder index}
Let us recall a few definitions from the theory of Maslov index.
Let $V$ be a finite dimensional real vector space endowed with a symplectic form $\omega$,
and let $\Spl(V,\omega)$ denote the symplectic group of $(V,\omega)$; set $\Dim(V)=2n$.
Denote by $\LLambda=\LLambda(V,\omega)$ the Grassmannian
of all $n$-dimensional subspaces of $V$, which is a
$\frac12n(n+1)$-dimensional real-analytic compact manifold.
For $L_0\in\LLambda$, one has a smooth fibration $\beta_{L_0}:\Spl(V,\omega)\to\LLambda$
defined by:
\[\beta_{L_0}(\Phi)=\Phi[L_0].\]
Let $L_0,L_1\in\LLambda$ be transverse Lagrangians; any other Lagrangian $L$ which is transverse
to $L_1$  is the graph
of a unique linear map $T:L_0\to L_1$; we will denote by $\varphi_{L_0,L_1}(L)$
is defined to be the restriction of the bilinear map
$\omega(T\cdot,\cdot)$ to $L_0\times L_0$, which is a symmetric bilinear form on $L_0$.
For $L\in\LLambda$, we will denote by $\LLambda_0(L)$ the set of all
lagrangians $L'\in\LLambda$ that are transverse to $L$; this is a dense open subset of $\LLambda$.

Denote by $\pi(\LLambda)$ the \emph{fundamental groupoid} of $\LLambda$,  endowed with the partial operation
of concatenation $\diamond$.
For all $L_0\in\LLambda$, there exists a unique $\Z$-valued
groupoid homomorphism
$\mu_{L_0}$ on $\pi(\LLambda)$ such that:
\begin{equation}\label{eq:muL0tilde}
\mu_{L_0}\big([\gamma]\big)=\mathrm n_+\big(\varphi_{L_0,L_1}(\gamma(1))\big)+\Dim\big(\gamma(1)\cap L_0\big)-
\mathrm n_+\big(\varphi_{L_0,L_1}(\gamma(0))\big)-\Dim\big(\gamma(0)\cap L_0\big)
\end{equation}
for all continuous curve $\gamma:[0,1]\to\LLambda_0(L_1)$ and for all
$L_1\in\LLambda_0(L_0)$.
The map $\mu_{L_0}:\pi(\LLambda)\to\Z$ is called the \emph{$L_0$-Maslov index}.

Given four Lagrangians $L_0,L_1,L_0',L_1'\in\LLambda$ and
any continuous curve $\gamma:[a,b]\to\LLambda$ such that
$\gamma(a)=L_0'$ and $\gamma(b)=L_1'$, then the value of
the quantity $\mathfrak q(L_0,L_1;L_0',L_1')=\mu_{L_1}(\gamma)-\mu_{L_0}(\gamma)$ does
{\em not\/} depend on the choice of $\gamma$, and it is called the
{\em H\"ormander index\/} of the quadruple $(L_0,L_1;L_0',L_1')$.
Consider the direct
sum $V^2=V\oplus V$, endowed with the symplectic form $\omega^2=\omega\oplus(-\omega)$,
defined by:
\[\omega^2\big((v_1,v_2),(w_1,w_2)\big)=\omega(v_1,v_2)-\omega(w_1,w_2),\quad v_1,v_2,w_1,w_2\in V,\]
and let $\Delta\subset V^2$ denote the diagonal subspace.
If $\Phi\in\Spl(V,\omega)$, then the graph of $\Phi$, denoted by $\Gr(\Phi)$, is given by $(\mathrm{Id}\oplus\Phi)[\Delta]\in\LLambda(V^2,\omega^2)$;
in particular $\Delta=\Gr(\mathrm{Id})$ and $\Delta^o=\{(v,-v):v\in V\}=\Gr(-\mathrm{Id})$ are Lagrangian subspaces
of $V^2$.
Given a continuous curve $\Phi$ in
$\Spl(V,\omega)$, the \emph{Conley--Zehnder index $\iCZ(\Phi)$ of $\Phi$}
 is the $\Delta$-Maslov
index of the curve $t\mapsto\Gr\big(\Phi(t)\big)\in\LLambda(V^2,\omega^2)$:
\[\iCZ(\Phi):=\mu_\Delta\big(t\mapsto\Gr(\Phi(t))\big).\]
We have the following relation between the Maslov index and the H\"ormander index:
\begin{lem}\label{thm:comind}
Let $\Phi:[a,b]\to\Spl(V,\omega)$ be a continuous curve and
let $L_0,L_1,L_1'\in\LLambda(V,\omega)$ be fixed. Then:
\[\mu_{L_0}\big(\beta_{L_1}\circ\Phi\big)-\mu_{L_0}\big(\beta_{L_1'}\circ\Phi\big)=
\mathfrak q\big(L_1,L_1';\Phi(a)^{-1}(L_0),\Phi(b)^{-1}(L_0)\big).\]
\end{lem}
\begin{proof}
Using the Maslov index for pairs and the symplectic
invariance, we compute as follows:
\[\mu_{L_0}\big(\beta_{L_1}\circ\Phi\big)=\mu\big(\beta_{L_1}\circ\Phi,L_0\big)
=\mu\big(L_1,t\mapsto\Phi(t)^{-1}(L_0)\big)=-\mu_{L_1}\big(t\mapsto\Phi(t)^{-1}(L_0)
\big).
\]
Similarly,
\[\mu_{L_0}\big(\beta_{L_1'}\circ\Phi\big)=
-\mu_{L_1'}\big(t\mapsto\Phi(t)^{-1}(L_0)\big).\]
The conclusion follows easily from the  definition
of $\mathfrak q$.
\end{proof}
The following relation between the notions of Maslov, Conley--Zehnder and H\"ormander index holds:
\begin{prop}
\label{thm:compind2}
Let $\Phi:[a,b]\to\Spl(V,\omega)$ be a continuous curve and $L_0,\ell_0\in\LLambda(V,\omega)$ be
fixed. Then:
\[\iCZ(\Phi)+\mu_{L_0}\big(\beta_{\ell_0}\circ\Phi\big)=
\mathfrak q\big(\Delta,L_0\oplus\ell_0;\Gr\big(\Phi(a)^{-1}\big),
\Gr\big(\Phi(b)^{-1}\big)\big).\]
In particular, if $\Phi$ is a loop, then $\iCZ(\Phi)=-\mu_{L_0}(\beta_{\ell_0}\circ\Phi)$.
\end{prop}
\begin{proof}
We compute:
\[\iCZ(\Phi)=\mu_\Delta\big(t\mapsto(\mathrm{Id}\oplus\Phi(t))(\Delta)\big)
\]
and, using the properties of the Maslov index for pairs of curves,
\[\mu_{L_0}\big(\beta_{\ell_0}\circ\Phi\big)=-\mu_{\Delta}\big(t\mapsto L_0\oplus
\beta_{\ell_0}\circ\Phi(t)\big)=-\mu_\Delta\big(t\mapsto(\mathrm{Id}\oplus
\Phi(t))(L_0\oplus\ell_0)\big).\]
The result follows now easily applying Lemma~\ref{thm:comind} to
the curve $t\mapsto\mathrm{Id}\oplus\Phi(t)$ in the symplectic group $\Spl(V^2,\omega^2)$ and
to the Lagrangians $\Delta,L_0\oplus\ell_0\in\LLambda(V^2,\omega^2)$.
\end{proof}
\subsection{Periodic solutions of Hamiltonian systems}\label{sub:persolham}
The notion of Conley--Zehnder index is used in the theory of periodic solutions
for Hamiltonian systems. Let us recall a few basic facts; let $(\mathcal M,\varpi)$ be a
$2n$-dimensional symplectic manifold, and let $H:\mathcal M\times\R\to\R$
be a (possibly time-dependent) smooth Hamiltonian. Assume that $H$ is $T$-periodic in time,
and that $z:[0,T]\to M$ is a solution of $H$ (i.e., $\dot z=\vec H(z)$
such that $z(0)=z(T)$,
where $\vec H$ is the time-dependent Hamiltonian vector field, defined by $\varpi(\vec H,\cdot)=\mathrm dH$).
Then, the iterates $z^{(N)}$ of $z$, defined as the concatenation:
\[z^{(N)}=\underbrace{z\diamond\cdots\diamond z}_{\text{$N$-times}}:[0,NT]\longrightarrow\mathcal M\]
are also solutions of $H$. Assume that the following objects are given:
\begin{itemize}
\item a \emph{periodic symplectic trivialization} of the tangent bundle of $\mathcal M$ along $z$
(i.e., of the pull-back $z^*T\mathcal M$), which consists of a smooth family $\Psi=\{\psi_t\}_{t\in[0,T]}$
of symplectomorphisms $\psi_t:T_{z(0)}\mathcal M\to T_{z(t)}\mathcal M$ with $\psi_0=\psi_T=\mathrm{Id}$;
\item a Lagrangian subspace $L_0\subset T_{z(0)}\mathcal M$.
\end{itemize}
By a simple orientability argument, periodic symplectic
trivializations along periodic solutions always exist.
By the periodicity assumption, we have a smooth extension $\R\ni t\mapsto\psi_t$ by setting
$\psi_{t+NT}=\psi_t$ for all $t\in[0,T]$.
Denote by $\mathcal F^H_{t,t'}:\mathcal M\to\mathcal M$ the flow of $\vec H$,\footnote{%
For our purposes, we will not be interested in questions of global existence of the flow
$\mathcal F^H$.}
i.e., $\mathcal F^H_{t,t'}(p)=\gamma(t')$, where $\gamma$ is the unique
integral curve of the time-dependent vector field $\vec H$ on $\mathcal M$
satisfying $\gamma(t)=p$. It is well known that for all $t,t'$, the $\mathcal F^H_{t,t'}$ is
a symplectomorphism between open subsets of $\mathcal M$. Left composition with $\psi_t^{-1}$ gives
a smooth map $\R\ni t\mapsto \Psi(t)=\psi_t^{-1}\circ\mathcal F^H_{0,t}\big(z(0)\big)$
of linear symplectomorphisms of $T_{z(0)}\mathcal M$; clearly $X(t)=\Psi'(t)\Psi(t)^{-1}$ lies in the Lie algebra
$\spl\big(T_{z(0)}\mathcal M,\varpi_{z(0)}\big)$ of the symplectic group $\Spl\big(T_{z(0)}\mathcal M,\varpi_{z(0)}\big)$.

The \emph{linearized Hamilton equation along $z$} is the linear system
\begin{equation}\label{eq:linhamequations}
v'(t)=X(t)v(t),
\end{equation}
in $T_{z(0)}\mathcal M$; the fundamental solution of this linear system
is a smooth symplectic path $\Phi:\R\to\Spl\big(T_{z(0)}\mathcal M,\varpi_{z(0)}\big)$ that satisfies $\Phi(0)=\mathrm{Id}$
and $\Phi'=X\Phi$.

\begin{defin}
The \emph{Conley--Zehnder index} of the solution $z=z^{(1)}$ associated to the symplectic trivialization
$\Psi$, denoted by $\iCZ(z,\Psi)$, is
the Conley--Zehnder of the path in $\Spl\big(T_{z(0)}\mathcal M,\varpi_{z(0)}\big)$ obtained by restriction
of the fundamental solution $\Phi$ to the interval $[0,T]$. Similarly, the \emph{$L_0$-Maslov index}
of the solution $z$ associated to the symplectic trivialization
$\Psi$, denoted by $\mu_{L_0}(z,\Psi)$, is
the $L_0$-Maslov index  of the path in $\Spl\big(T_{z(0)}\mathcal M,\varpi_{z(0)}\big)$
given by $[0,T]\ni t\mapsto\Phi(t)[L_0]\in\LLambda\big(T_{z(0)}\mathcal M,\varpi_{z(0)}\big)$.
\end{defin}

\begin{rem}
We observe here that both the notions of Conley--Zehnder index and of Maslov index for a periodic
solution $z$ of a Hamiltonian system depend on the choice of a symplectic trivialization. More precisely,
given two periodic symplectic trivializations $\Psi=\{\psi_t\}_t$, $\widetilde\Psi=\{\tilde\psi_t\}_t$
and setting $G_t=\psi_t^{-1}\circ\tilde\psi_t\in\Spl\big(T_{z(0)}\mathcal M,\varpi_{z(0)}\big)$,
the corresponding paths $\Phi$ and $\widetilde\Phi$ in $\Spl\big(T_{z(0)}\mathcal M,\varpi_{z(0)}\big)$ are
related by:
\[\Phi(t)=G_t\circ\widetilde\Phi(t),\quad\forall\, t\in[0,T].\]
Clearly, $[0,T]\ni t\mapsto G_t$ is a closed path in $\Spl\big(T_{z(0)}\mathcal M,\varpi_{z(0)}\big)$ with
endpoint in the identity; in this situation, one proves easily\footnote{For instance, using the product formula in
\cite[Lemma~3.3]{GosGosPic}.} that $\iCZ(z,\Psi)=\iCZ(G)+\iCZ(z,\widetilde\Psi)$.
In particular, if the loop $G$ is homotopically trivial, then $\iCZ(z,\Psi)=\iCZ(z,\widetilde\Psi)$.
Similarly, if $G_t[L_0]=L_0$ for all $t$, and if $G$ is homotopically trivial, then $\mu_{L_0}(z,\Psi)=\mu_{L_0}(z,\widetilde\Psi)$.
\end{rem}
This observation will be used in a situation described in the following Lemma:
\begin{lem}\label{thm:trivialloop}
Let $\mathcal V$ be a finite dimensional vector space and set $V=\mathcal V\oplus\mathcal V^*$; $V$ is a symplectic space, endowed
with its \emph{canonical symplectic form} $\omega\big((v,\alpha),(w,\beta)\big)=\beta(v)-\alpha(w)$, $v,w\in\mathcal V$, $\alpha,\beta\in\mathcal V^*$.
Given any $\eta\in\GL(\mathcal V)$, then the linear map
\[G=\begin{pmatrix}\eta&0\cr0&{\eta^*}^{-1}\end{pmatrix}:V\to V\] is a symplectomorphism of $(V,\omega)$.
If $[a,b]\ni t\mapsto G_t\in\Spl(V,\omega)$ is a continuous map of symplectomorphisms of this type with $G_a=G_b=\mathrm{Id}$,
then $G$ is homotopically trivial in $\Spl(V,\omega)$.
\end{lem}
\begin{proof}
The first statement is immediate. In order to prove that $G$ is homotopically trivial, it is not restrictive to
assume $\mathcal V=\R^n$; identifying ${\R^n}^*$ with $\R^n$ via the Euclidean inner product, we will consider the
canonical complex structure on $V\cong\R^{2n}$. The thesis is obtained if we prove that, denoting by $G_t=u_tp_t$ the polar decomposition
of $G_t$, with $u_t$ unitary and $p_t$ positive definite, then $t\mapsto u_t$ is homotopically trivial
in $\mathrm U(n)$. This is equivalent to the fact that the closed in loop $t\mapsto\mathrm{det}(u_t)\in\mathds S^1$
is homotopically trivial in $\mathds S^1$. If $\eta_t=o_tq_t$ is the polar decomposition of $\eta_t$, with $o_t\in\mathrm O(n)$ and
$q_t$ positive definite, then the unitary $u_t$ is given by $\begin{pmatrix} o_t&0\cr0& o_t\end{pmatrix}\in\mathrm U(n)$, which has
constant determinant equal to $1$. The conclusion follows easily, recalling that the determinant map $\mathrm{det}:\mathrm U(n)\to\mathds S^1$
induces an isomorphism between the fundamental groups.
\end{proof}

\subsection{An iteration formula for the Maslov index}
Let us recall the following iteration formula for the Conley--Zehnder index, proved in \cite{GosGosPic}:
\begin{prop}
In the notations of Subsection~\ref{sub:persolham}, the following inequality holds:
\begin{equation}\label{eq:itformCZ}
\left\vert\iCZ\big(z^{(N)},\Psi\big)-N\cdot\iCZ\big(z,\Psi\big)\right\vert\le n(N-1).
\end{equation}
In particular, $\left\vert\iCZ\big(z^{(N)},\Psi\big)\right\vert$ has sublinear growth in $N$.
Moreover, if $\big\vert\iCZ\big(z,\Psi\big)\big\vert>n$, then $\iCZ\big(z^{(N)},\Psi\big)$
 has superlinear growth in $N$.
\end{prop}
\begin{proof}
See \cite[Corollary~4.4]{GosGosPic}. Observe that we are using here a slightly different definition of
Conley--Zehnder index, and the inequality \eqref{eq:itformCZ} differs by a factor $2$ from the corresponding
inequality in \cite[Corollary~4.4]{GosGosPic}.
\end{proof}
Let us prove that a similar iteration formula holds for the Maslov index:
\begin{cor}\label{thm:coriteformMaslovindex}
 The following inequality holds:
\[
\left\vert\,\mu_{L_0}\big(z^{(N)},\Psi\big)-N\cdot\mu_{L_0}\big(z,\Psi\big)\right\vert\le n(7N+5).
\]
In particular, $\left\vert\mu_{L_0}\big(z^{(N)},\Psi\big)\right\vert$ has sublinear growth in $N$;
moreover,  if $\mu_{L_0}\big(z,\Psi\big)>7n$, then $\mu_{L_0}\big(z^{(N)},\Psi\big)$
has superlinear growth in $N$.
\end{cor}
\begin{proof}
The inequality is obtained easily from \eqref{eq:itformCZ}, using Proposition~\ref{thm:compind2}:
\begin{multline*}\left\vert\mu_{L_0}\big(z^{(N)},\Psi\big)-N\cdot\mu_{L_0}\big(z,\Psi\big)\right\vert\le
\left\vert\iCZ\big(z^{(N)},\Psi\big)-N\cdot\iCZ\big(z,\Psi\big)\right\vert\\+\left\vert\mathfrak q\big(\Delta,L_0\oplus L_0;\Delta,\Gr\big(\Phi(NT)\big)\big)
-N\cdot \mathfrak q\big(\Delta,L_0\oplus L_0;\Delta,\Gr\big(\Phi(T)\big)\big)\right\vert\\ \le n(N-1)+6n(N+1)=n(7N+5).\qedhere\end{multline*}
\end{proof}

\subsection{Maslov index of a geodesic and of the corresponding Hamiltonian solution}
\label{sub:Maslovgeo}
Let us now define the notion of Maslov index for a closed geodesic $\gamma$ in a semi-Riemannian
manifold $(M,g)$; we will show that when $\gamma$ is orientation preserving, then its Maslov
index coincides with the Maslov index of the corresponding periodic solution of the
geodesic Hamiltonian in the cotangent bundle $TM^*$.

Let us recall
the notion of Maslov index for a \emph{fixed endpoint geodesic}. If $\gamma:[0,1]\to M$ is any geodesic,
consider a continuous trivialization of $TM$ along $\gamma$, i.e., a continuous family
of isomorphisms $h_t:T_{\gamma(0)}M\to T_{\gamma(t)}M$, $t\in[0,1]$. Consider the symplectic space
$V=T_{\gamma(0)}M\oplus T_{\gamma(0)}M^*$ endowed with its canonical symplectic structure
(recall Lemma~\ref{thm:trivialloop}), the Lagrangian subspace $L_0=\{0\}\oplus T_{\gamma(0)}M^*$,
and the continuous curve of Lagrangians $\ell(t)\in\LLambda(V,\omega)$ given by:
\[\ell(t)=\Big\{\big(h_t^{-1}[J(t)],h_t^*[g(\Ddt J(t))]\big):\text{$J$ Jacobi field along $\gamma$, with $J(0)=0$}\Big\}.\]
In the above formula, the metric tensor $g$ is seen as a map $g:T_{\gamma(t)}M\to T_{\gamma(t)}M^*$.
The Maslov index of $\gamma$, denoted by $\iMaslov(\gamma)$ is defined as the $L_0$-Maslov index
of the continuous path $[0,1]\ni t\mapsto\ell(t)$.\footnote{A different convention was originally
adopted by Helfer \cite{Hel} in the definition of Maslov index of a semi-Riemannian geodesic.
In Helfer's original definition, given a geodesic $\gamma:[0,1]\to M$ with non conjugate endpoints,
$\iMaslov(\gamma)$ was given by the $L_0$-Maslov index of the continuous path $[\varepsilon,1]\ni t\mapsto\ell(t)$,
where $\varepsilon>0$ is small enough so that there are no conjugate instants in $\left]0,\varepsilon\right]$.
This convention was motivated by the necessity of avoiding dealing with curves in the Lagrangian Grassmannian with endpoints
in the Maslov cycle. An immediate calculation using  \eqref{eq:muL0tilde}, shows that, if $g$ is Lorentzian,
the following simple relation holds: $\iMaslov\big(\gamma\vert_{[\varepsilon,1]}\big)=\iMaslov(\gamma)+1$.}
This quantity does not depend on the choice of the trivialization of $TM$ along $\gamma$.
Let us now consider the case of a closed geodesic, in which case one may study the existence
of \emph{periodic} trivializations of $TM$ along $\gamma$.

Recall that a closed curve $\gamma:[a,b]\rightarrow M$ is said to be {\em
orientation preserving\/} if for some (and hence for any) continuous
trivialization $h_t:T_{\gamma(a)}M\to T_{\gamma(t)}M$, $t\in[a,b]$, of $TM$ along $\gamma$, the
isomorphism $h_b^{-1}\circ h_a:T_{\gamma(a)}M\to T_{\gamma(a)}M$ is orientation preserving.
It is easy to prove that if $\gamma$ is orientation preserving then there exists a smooth
trivialization $h_t:T_{\gamma(a)}M\to T_{\gamma(t)}M$, $t\in[a,b]$, of $TM$ along $\gamma$
with $h_b^{-1}\circ h_a$ the identity of $T_{\gamma(a)}M$.

Assume that $\gamma:[0,1]\to M$ is a closed geodesic in $M$, which is orientation
preserving. Let $\Gamma:[0,1]\to TM^*$ be the corresponding periodic solution
of the geodesic Hamiltonian:
\[H(q,p)=g^{-1}(p,p).\]
Given a smooth periodic trivialization of $TM$ along $\gamma$, $h_t:T_{\gamma(0)}M\to T_{\gamma(t)}M$, $t\in[0,1]$,
$h_0=h_1$, then one can define a smooth periodic symplectic trivialization of the tangent bundle $T(TM^*)$ along
$\Gamma$ as follows. Denote by $\pi:TM^*\to M$ the canonical projection; for $p\in TM^*$, denote by $\Ver_p=\Ker(\dd\pi_p)$
the vertical subspace of $T_p(TM^*)$ and by $\Hor_p$ the horizontal subspace of $T_p(TM^*)$ relatively to the Levi--Civita
connection $\nabla$. One has a canonical identification $\Ver_p=T_p(T_xM^*)\cong (T_xM)^*$, while
the restriction of the differential $\dd\pi_p$ to $\Hor_p$ gives an identification $\Hor_p\cong T_xM$, where $x=\pi(p)$.
Since $\nabla$ is torsion free, with these identifications, the \emph{canonical symplectic form} $\varpi$ of $TM^*$ at $p\in TM^*$ becomes the canonical symplectic
form of $T_xM\oplus (T_xM)^*$; moreover,  for all $t\in[0,1]$ we define an isomorphism:
\begin{multline*}\psi_t:T_{\Gamma(0)}(TM^*)=\Hor_{\Gamma(0)}\oplus\Ver_{\Gamma(0)}\cong T_{\gamma(0)}M\oplus (T_{\gamma(0)}M)^*
\\ \longrightarrow T_{\gamma(t)}M\oplus (T_{\gamma(t)}M)^*\cong\Hor_{\Gamma(t)}\oplus\Ver_{\Gamma(t)}=T_{\Gamma(t)}(TM^*)
\end{multline*}
by setting:
\[\psi_t(v,\alpha)=\big(h_t(v),{h_t^*}^{-1}(\alpha)\big),\]
for all $v\in T_{\gamma(0)}M$ and $\alpha\in(T_{\gamma(0)}M)^*$. This is obviously a symplectomorphism for all $t$,
hence we obtain a smooth periodic symplectic trivialization $\Psi=\{\psi_t\}_{t\in[0,1]}$ of $T(TM^*)$ along
$\Gamma$. It is immediate to observe that the Maslov index $\iMaslov(\gamma)$ of the geodesic $\gamma$ coincides
with the $L_0$-Maslov index $\mu_{L_0}(\Gamma,\Psi)$ of the solution $\Gamma$ associated to the symplectic trivialization
$\Psi$, where $L_0$ is the Lagrangian subspace $\{0\}\oplus (T_{\gamma(0)}M)^*$ of $T_{\gamma(0)}M\oplus (T_{\gamma(0)}M)^*$.
\begin{lem} Let $\gamma$ be an orientation preserving closed geodesic in $(M,g)$, and let $\Gamma$ be the corresponding
periodic solution of the geodesic Hamiltonian in $TM^*$.
The $L_0$-Maslov index $\mu_{L_0}(\Gamma,\Psi)$, where $\Psi$ is the smooth periodic trivialization of $T(TM^*)$ along
$\Gamma$ constructed from a smooth periodic trivialization $\{h_t\}_{t\in[0,1]}$ of $TM$ along $\gamma$, as described
above, does not depend on the choice of $\{h_t\}_{t\in[0,1]}$.
\end{lem}
\begin{proof}
This is an immediate consequence of Lemma~\ref{thm:trivialloop}, observing that two distinct trivializations
$\{h_t\}$ and $\{\tilde h_t\}$ of $TM$ along $\gamma$, with $\eta_t=\tilde h_t\circ h_t\in\GL\big(T_{\gamma(0)}\big)M$,
yield periodic symplectic trivializations $\{\psi_t\}$ and $\{\tilde\psi_t\}$ of $T(TM^*)$ along
$\Gamma$ that \emph{differ} by a loop
$\{G_t\}$ in $\GL\big(T_{\Gamma(0)}\big)$ of the form:
\[G_t=\begin{pmatrix}\eta_t&0\\0&{\eta_t^*}^{-1}\end{pmatrix}.\]
By Lemma~\ref{thm:trivialloop}, this loop is contractible in $\Spl\big(T_{\Gamma(0)}(TM^*),\varpi_{\Gamma(0)}\big)$, and
clearly $G_t[L_0]=L_0$ for all $t$, which concludes the proof.
\end{proof}

Using the construction above and Corollary~\ref{thm:coriteformMaslovindex} we obtain immediately:

\begin{cor}\label{thm:iterateMaslovgeo}
Let $\gamma$ be an orientation preserving closed geodesic in $(M,g)$.
Then, denoting by $\gamma^{(N)}$ the $N$-th iterated of $\gamma$, $N\ge1$, the following inequality holds:
\[ \left\vert\,\iMaslov\big(\gamma^{(N)}\big)-N\cdot \iMaslov(\gamma)\right\vert\le \Dim(M)(7N+5).
\]
In particular, $\left\vert\iMaslov\big(\gamma^{(N)}\big)\right\vert$ has sublinear growth in $N$;
moreover,  if $\iMaslov(\gamma)>7\,\Dim(M)$, then $\iMaslov\big(\gamma^{(N)}\big)$
has superlinear growth in $N$.\qed
\end{cor}

\end{section}

\begin{section}{The variational setup}
\label{sec:varprobPS}
\noindent

Let $(M,g)$ be a stationary Lorentzian manifold, and let $\mathcal Y\in\mathfrak X(M)$ be a timelike Killing vector field
in $M$. Consider the auxiliary Riemannian metric $\gR$ on $M$, defined by
\begin{equation}\label{eq:defgr}
\gR(v,w)=g(v,w)-2\frac{g(v,\mathcal Y) g(w,\mathcal Y)}{g(\mathcal Y,\mathcal Y)};
\end{equation}
observe that $\mathcal Y$ is Killing also relatively to $\gR$.
Let $\mathds S^1$ be the unit circle, viewed as the quotient $[0,1]/\{0,1\}$, and denote by $\Lambda M=H^1(\mathds S^1,M)$ the infinite dimensional Hilbert manifold of
all loops $\gamma:[0,1]\to M$, i.e., $\gamma(0)=\gamma(1)$, of Sobolev class $H^1$;
if $\Lambda^0M$ is the set of continuous loops in $M$ endowed with the compact-open topology, the inclusion $\Lambda M\hookrightarrow\Lambda^0M$ is
a homotopy equivalence (this can be proved, for instance, using the results in \cite{Palais66}).
Set:
\[\mathcal N=\Big\{\gamma\in\Lambda M:g(\dot\gamma,\mathcal Y)=c_\gamma\ \text{(constant)}\ \text{a.e.\ on $\mathds S^1$}\Big\}.\]
For all $\gamma\in\Lambda M$, the tangent space $T_\gamma\Lambda M$ is identified with the space
of all sections of the pull-back $\gamma^*(TM)$ (i.e., periodic vector fields along $\gamma$)
of Sobolev class $H^1$; this space will be endowed with the
Hilbert space inner product:
\begin{equation}\label{eq:metricaHilbert}
\llangle V,W\rrangle=\int_0^1\Big[\gR(V,W)+\gR\big(\DdtR V,\DdtR W\big)\Big]\,\mathrm dt,
\end{equation}
where $\DdtR$ denotes the covariant differentiation along $\gamma$ relatively to the Levi--Civita
connection of the metric $\gR$.

Recall the definition of the classical geodesic energy functional on $\Lambda M$:
\[f(\gamma)=\tfrac12\int_{0}^1g(\dot\gamma,\dot\gamma)\,\mathrm dt.\]
\subsection{The constrained variational problem}
It is well known that the critical points of $f$ in $\Lambda M$ are exactly the closed geodesics in
$M$; it is also clear that the set $\mathcal N$ contains the closed geodesics in $M$.
It is proven that the equality $g(\dot\gamma,\mathcal Y)=c_\gamma$ provides a \emph{natural constraint}
for the critical points of the geodesic action functional in a stationary Lorentzian manifold; more precisely:
\begin{prop}\label{thm:constrainedvarprob}
The following statements hold:
\begin{enumerate}
\item $\mathcal N$ is a smooth embedded closed submanifold of $\Lambda M$, and for $\gamma\in\mathcal N$, the
tangent space $T_\gamma\mathcal N$ is given by the space of sections $V$ of the pull-back $\gamma^*(TM)$ of Sobolev class $H^1$,
satisfying:
\begin{equation}\label{eq:defCV}
g\big(\Ddt V,\mathcal Y\big)-g\big(V,\Ddt\mathcal Y\big)=C_V\ \ \text{(constant)  a.e.\ on $[0,1]$};
\end{equation}
\item\label{itm:Ycomplete} if $\mathcal Y$ is complete, then $\mathcal N$ is a strong deformation retract of $\Lambda M$ (hence it is
homotopy equivalent to $\Lambda M$);
\item a curve $\gamma\in\mathcal N$ is a critical point of the restriction of $f$ to $\mathcal N$ if and only if
$\gamma$ is a critical point of $f$ in $\Lambda M$, i.e., if and only if $\gamma$ is a closed geodesic
in $(M,g)$;
\item if $\gamma$ is a critical point of $f$, then the Hessian $\mathrm H^{f\vert_{\mathcal N}}$ of the restriction $f\vert_{\mathcal N}$
at $\gamma$ is given by the restriction of the index form:
\[I_\gamma(V,W)=\int_0^1g\big(\Ddt V,\Ddt W\big)+g\big(R_{\gamma(t)}(\dot\gamma,V)\,\dot\gamma,W\big)\,\mathrm dt\]
to the tangent space $T_\gamma\mathcal N$;
\item\label{thm:constrainedvarprobitm5} if $\gamma$ is a critical point of $f$, then the index form
$I_\gamma$ is essentially positive on $T_\gamma\mathcal N$,
and in particular the Morse index of $f\vert_{\mathcal N}$ at $\gamma$ is finite.
\end{enumerate}
\end{prop}
\begin{proof}
See \cite{asian, CAG1, Masiello, MarPicTau}.
\end{proof}

It is clear that $f$ does \emph{not} satisfy the Palais--Smale condition in $\mathcal N$;
namely, all its critical orbits are non compact.
\smallskip

Given a closed geodesic $\gamma$ in $(M,g)$, let us denote by $\mu(\gamma)$ the Morse index
of $f\vert_{\mathcal N}$ at $\gamma$, i.e., the index of the restriction of $I_\gamma$ to $T_\gamma\mathcal N$.
This index will be computed explicitly using the Morse index theorem (Theorem~\ref{thm:MorseIndexTheorem})
in Section~\ref{sec:indexthm}.
Moreover, let us denote by $\mu_0(\gamma)$ the \emph{extended index} if $f\vert_{\mathcal N}$ at $\gamma$,
which is the sum of the index $\mu(\gamma)$ and the \emph{nullity} $\mathrm n(\gamma)$:
\[\mathrm n(\gamma)=\Dim\big[\Ker\big(I_\gamma\vert_{T_\gamma\mathcal N\times T_\gamma\mathcal N}\big)\big].\]
We will establish in Lemma~\ref{thm:kernelIgamma} that $\mathrm n(\gamma)$ equals the dimension of
the space of periodic Jacobi fields along $\gamma$.

\subsection{The Palais--Smale condition}

Let us now assume that $(M,g)$ is a globally hyperbolic stationary Lorentzian manifold,
that admits a complete timelike Killing vector field $\mathcal Y$.
Let us recall that, in this situation, $(M,g)$ is a \emph{standard}
stationary manifold (see for instance \cite[Theorem~2.3]{CanFloSan}), i.e., denoting by $S$ a
smooth Cauchy surface of $M$, then $M$ is diffeomorphic to a product $S\times\R$,
and the Killing field $\mathcal Y$ is the vector field $\partial_t$ which is
tangent to the fibers $\{x\}\times\R$. One should observe that such product decomposition
of $M$ is \emph{not} canonical; however, all Cauchy surfaces of $M$ are homeomorphic.
In particular, $M$ is simply connected if and only if $S$ is, and
the inclusion of the free loop space $\Lambda S\hookrightarrow\Lambda M$
is a homotopy equivalence.

The projection onto the second factor $S\times\R\to\R$,
that will be denoted by $T$, is a smooth time function, that satisfies:
\begin{equation}\label{eq:nablaYcal}
\mathcal Y(T)=g\big(\nabla T,\mathcal Y\big)\equiv1
\end{equation}
on $M$. If $\mathbf L$ denotes the Lie derivative, from \eqref{eq:nablaYcal} it follows that $\mathbf L_{\mathcal Y}(\mathrm dT)$ vanishes identically.
For, given an arbitrary smooth vector field $X$ on $M$:
\[\mathbf L_{\mathcal Y}(\mathrm dT)(X)=\mathcal Y\big(X(T)\big)-\mathrm dt\big([\mathcal Y,X]\big)=\mathcal Y\big(X(T)\big)-\mathcal Y\big(X(T)\big)+X\big(\mathcal Y(T)\big)=0.\]
Since $\mathcal Y$ is Killing, then $\mathbf L_{\mathcal Y}(g)=0$, and \eqref{eq:nablaYcal} implies that the Lie bracket $[\mathcal Y,\nabla T]=
\mathbf L_{\mathcal Y}\big(g^{-1}\mathrm dT\big)$ also vanishes identically. It follows that the quantity $g(\nabla T,\nabla T)$ is constant
along the flow lines of $\mathcal Y$:
\[\mathcal Y\,g(\nabla T,\nabla T)=2g\big(\nabla_{\mathcal Y}\nabla T,\nabla T\big)=2g\big([\mathcal Y,\nabla T]-\nabla_{\nabla T}\mathcal Y,\nabla T\big)
=0.\]
\begin{lem}\label{thm:boundedbelow}
The restriction of the functional $f$ to $\mathcal N$ is bounded from below; more precisely, $f(\gamma)\ge0$ for
all $\gamma\in\mathcal N$, and $f(\gamma)=0$ only if $\gamma$ is a constant curve.
\end{lem}
\begin{proof}
Let $\gamma\in\mathcal N$ be fixed, and denote by $c_\gamma$ the value of the constant $g(\dot\gamma,\mathcal Y)$.
For almost all $t\in[0,1]$, the vector $\dot\gamma-g(\dot\gamma,\nabla T)\,\mathcal Y$ is (null or) spacelike, namely,
using \eqref{eq:nablaYcal}, one checks immediately that
it is orthogonal to the timelike vector $\nabla T$. Hence:
\begin{equation}\label{eq:primadisug}
0\le g\big(\dot\gamma-g(\dot\gamma,\nabla T)\,\mathcal Y,\dot\gamma-g(\dot\gamma,\nabla T)\,\mathcal Y\big)=g(\dot\gamma,\dot\gamma)
-2\,c_\gamma\,g\big(\dot\gamma,\nabla T\big)+g(\dot\gamma,\nabla T)^2g(\mathcal Y,\mathcal Y),\end{equation}
and thus:
\[g(\dot\gamma,\dot\gamma)\ge 2\,c_\gamma\,g\big(\dot\gamma,\nabla T\big)-g(\dot\gamma,\nabla T)^2g(\mathcal Y,\mathcal Y).\]
Integrating on $[0,1]$, and observing that since $\gamma$ is closed $\int_0^1g\big(\dot\gamma,\nabla T\big)\,\mathrm dt=0$, we get:
\begin{equation}\label{eq:secondadisug}2\,f(\gamma)\ge -\int_0^1g(\dot\gamma,\nabla T)^2g(\mathcal Y,\mathcal Y)\,\mathrm dt\ge0.\end{equation}
Equality in \eqref{eq:primadisug} holds only if $\dot\gamma-g(\dot\gamma,\nabla T)\,\mathcal Y=0$, while, in the last inequality
of \eqref{eq:secondadisug}, the equal sign holds only if $g(\dot\gamma,\nabla T)=0$ almost everywhere on $[0,1]$.
Hence, $f(\gamma)=0$ only if $\dot\gamma=0$ almost everywhere.
\end{proof}

We will assume in the sequel that the Cauchy surface $S$ is compact; recall that any two
Cauchy surfaces of a globally hyperbolic spacetime are homeomorphic.

\begin{lem}\label{thm:grcomplete}
The metric $\gR$ is complete, and thus $\Lambda M$ and $\mathcal N$ are  complete Hilbert manifolds when endowed
with the Riemannian structure \eqref{eq:metricaHilbert}.
\end{lem}
\begin{proof}
We will give an argument showing that, more generally, if a Riemannian manifold
admits a group of isometries all of whose orbits meet a given compact subset, then
the metric is complete.
Denote by $d_0$ the distance function induced by the Riemannian metric $\gR$.
It suffices to show that there exists $r>0$ such that, for every $p\in M$, the
closed $d_0$-ball $B[p;r]$ centered at $p$ and of radius $r$ is compact.
\marginpar{Esiste un altro\\ termine per\\ denotare il\\ ``radius of\\ compactness''??}
Clearly, for every single $p$ there exists $r(p)>0$ such that $B[p,r(p)]$ is compact;
we will call such $r(p)$ a \emph{radius of compactness} at $p$.\footnote{One can use,
alternatively, the function $\overline r'(p)=\sup\{r>0:\mathrm B[0;r]\subset\mathrm{Dom}(\exp_p)\big\}$, where
$\mathrm B[0;r]$ is the closed ball centered at $0$ and of radius $r$ in $T_pM$.
Clearly, $\overline r'\ge\overline r$} The map
$M\ni p\mapsto\overline r(p)=\sup\{r>0:r\ \text{is a radius of compactness at $p$}\big\}\in\left]0,+\infty\right]$
is lower semi-continuous.
Hence, given a compact Cauchy surface $S$ of $M$, $\overline r$ has a positive (possibly infinite)
minimum on $S$, and so there exists
$r>0$ which is a radius of compactness at all $p\in S$. Since the flow of $\mathcal Y$
preserves $\gR$, a radius of compactness at some $p\in S$ is a radius of compactness
at each point of the flow line of $\mathcal Y$ through $p$.
The conclusion follows easily from the fact that every flow line of $\mathcal Y$
has non empty intersection with $S$.
\end{proof}

The flow of the Killing vector field $\mathcal Y$ gives an isometric action of $\R$ in $\Lambda M$,
defined by $\R\times\Lambda M\ni(t,\gamma)\mapsto\mathcal F_t\circ\gamma\in\Lambda M$.
This action preserves $\mathcal N$, and the functional $f$ is invariant by this action;
the orbit of a critical point of $f$ consists of a collection of critical points
of $f$ with the same Morse index.
Such action is obviously free, and the quotient $\widetilde{\mathcal N}=\mathcal N/\R$
has the structure of a smooth manifold such that the product $\widetilde{\mathcal N}\times\R$ is diffeomorphic
to $\mathcal N$. For $\gamma\in\mathcal N$, we will denote by $[\gamma]$ its class in the quotient $\widetilde{\mathcal N}$;
the tangent space $T_{[\gamma]}\widetilde{\mathcal N}$ can be identified with:
\begin{equation}\label{eq:tansptildeN}
T_{[\gamma]}\widetilde{\mathcal N}\cong T_\gamma\mathcal N/E_\gamma,
\end{equation}
where $E_\gamma$ is the $1$-dimensional space of vector fields spanned by
the restriction of $\mathcal Y$ to $\gamma$. If $S$ is a Cauchy surface in $M$, then $\widetilde{\mathcal N}$ can
also be identified with the set:
\begin{equation}\label{eq:tildecalN}
\widetilde{\mathcal N}=\big\{\gamma\in\mathcal N:\gamma(0)\in S\big\};
\end{equation}
using this identification, for $\gamma\in\widetilde{\mathcal N}$ is given by:
\begin{equation}\label{eq:TtildecalN}
T_\gamma\widetilde{\mathcal N}=\big\{V\in T_\gamma\mathcal N:V(0)\in T_{\gamma(0)}S\big\}.
\end{equation}
Obviously, the quotient $\widetilde{\mathcal N}$ inherits an isometric action of $\mathrm O(2)$;
it should be observed that, if one uses the identification \eqref{eq:tildecalN}, then the action
of an element in $\mathrm O(2)$ is not simply a rotation in the parameter space, but a rotation
followed by a translation along the flow of $\mathcal Y$.

The function $f$ defines by quotient a smooth
function on $\widetilde{\mathcal N}$, that will still be denoted by $f$, and for which
the statement of Proposition~\ref{thm:constrainedvarprob} holds \emph{verbatim}.
In addition, $f$ satisfies the PS condition on $\widetilde{\mathcal N}$.

\begin{prop}\label{thm:PS}
$\widetilde{\mathcal N}$ is a complete Hilbert manifold, which is homotopically equivalent to $\mathcal N$ and
to $\Lambda M$. The critical points of the functional $f$ in $\widetilde{\mathcal N}$ correspond to orbits
\[[\gamma]=\{\mathcal F_t\circ\gamma\}_{t\in\R}\] where $\gamma$ is a closed geodesic in $M$;
the Morse index of each critical point $[\gamma]$ of $f$ equals the Morse index
of $\gamma$, while the nullity of $[\gamma]$ equals $\mathrm n(\gamma)-1$.
Moreover, $f$ satisfies the Palais--Smale condition in $\widetilde{\mathcal N}$.
\end{prop}
\begin{proof}
Most part of the statement is a direct consequence of the construction of $\widetilde{\mathcal N}$.
The statement on the Morse index and the nullity of a critical point $[\gamma]$ is obtained easily,
observing that the $1$-dimensional space $E_\gamma$ in formula \eqref{eq:tansptildeN} is contained in the kernel
of the index form $I_\gamma$ (see Lemma~\ref{thm:kernelIgamma} below).
The Palais--Smale condition is essentially the same as in \cite[Lemma~3.2]{Masiello2};
we will sketch here a more intrinsic proof along the lines of \cite{CanFloSan, CAG1}.
Using \cite[Section~5]{CAG1} and the compactness of $S$, for the PS condition it suffices to show that
$f$ is \emph{pseudo-coercive on $\widetilde{\mathcal N}$}, i.e., that given a sequence
$(\gamma_n)_{n\in\N}$ in $\widetilde{\mathcal N}$ such that $f(\gamma_n)$ is bounded, then
$\gamma_n$ admits a uniformly convergent subsequence.
Using the identification \eqref{eq:tildecalN}, let $(\gamma_n)_{n\in\N}$ be a sequence
in $\widetilde{\mathcal N}$ such that $f(\gamma_n)\le c$ for all $n$;
we claim that the real sequence $c_{\gamma_n}=g(\dot\gamma_,\mathcal Y)$ is bounded.
Namely, the vector field $\dot\gamma_n-c_{\gamma_n}\nabla T$ along $\gamma_n$ is a.e.\ spacelike or null
for all $n$, because it is a.e.\ orthogonal to the timelike vector field $\mathcal Y$. Hence,
\[\int_0^1g\big(\dot\gamma_n-c_{\gamma_n}\nabla T, \dot\gamma_n-c_{\gamma_n}\nabla T\big)\,\mathrm dt=2 f(\gamma_n)+c_{\gamma_n}^2
\int_0^1g(\nabla T,\nabla T)\,\mathrm dt\ge0,\]
that gives:
\[c_{\gamma_n}^2\le2f(\gamma_n)\left(\int_0^1-g(\nabla T,\nabla T)\,\mathrm dt\right)^{-1}.\]
Observe that the functions $g(\mathcal Y,\mathcal Y)$ and $g(\nabla T,\nabla T)$ admit minimum and maximum in $M$, because they are constant along
the flow lines of $\mathcal Y$, and because $S$ is compact. The claim on the boundedness
of $c_{\gamma_n}$ follows. From this, it follows that the sequence:
\[\int_0^1\gR(\dot\gamma_n,\dot\gamma_n)\,\mathrm dt=2f(\gamma_n)-2c_{\gamma_n}^2\int_0^1g(\mathcal Y,\mathcal Y)^{-1}\,\mathrm dt\]
is bounded. Since $\gR$ is complete and $S$ is compact, the theorem of Arzel\`a--Ascoli implies that,
up to subsequences, $\gamma_n$ is uniformly convergent  in $M$. This concludes the proof.
\end{proof}

From Lemma~\ref{thm:boundedbelow} and Proposition~\ref{thm:PS}, one obtains
the existence of one non trivial closed geodesic in $M$, as proved in \cite{Masiello2}.
Namely, using the theory of Ljusternik and Schnirelman, one shows the existence of
a sequence $\big([\gamma_r]\big)_{r\ge1}$ of critical points of $f\vert_{\widetilde{\mathcal N}}$
with $f(\gamma_r)\to\infty$. Thus, these critical points are not constant curves;
observe however that the Ljusternik--Schnirelman theory does not give information
on whether such curves are geometrically distinct. In the non simply connected case,
the following result follows immediately:
\begin{cor}
Let $(M,g)$ be a Lorentzian manifold that admits a complete timelike Killing vector field
and a compact Cauchy surface. Then, there is a closed geodesic in each free homotopy class of
$M$.\qed
\end{cor}

\begin{rem}\label{thm:remazioneO2}
The orthogonal group $\mathrm O(2)$ acts isometrically on $\Lambda M$ via the operation of
$\mathrm O(2)$ on the parameter circle $\mathds S^1$. It is easy to observe that the stabilizer
of each $\gamma\in\Lambda M$ with respect to this action is a finite cyclic subgroup of
$\mathrm{SO}(2)$ generated by the rotation of $\frac{2\pi}N$, for some
$N\ge1$. A closed $\gamma\in\Lambda M$ will be called
\emph{prime} if its stabilizer in $\mathrm O(2)$ is trivial, i.e., if $\gamma$ is not the
$N$-th iterate of some other curve in $\Lambda M$ with $N>1$.
The functional $f$ defined in $\Lambda M$ is invariant by the action of $\mathrm O(2)$; moreover, this
action leaves $\mathcal N$ invariant, and it commutes with the time translations $(\mathcal F_t\,\circ{})$.
We therefore get an equivariant and isometric action of $\mathrm O(2)$ on the manifold $\widetilde{\mathcal N}$
by $g\cdot[\gamma]=\big[g\cdot\gamma\big]$, $g\in\mathrm O(2)$. An element $[\gamma]\in\widetilde{\mathcal N}$
will be called prime if $\gamma$ is prime, in which case its orbit $\mathrm O(2)\cdot[\gamma]$ will
contain only prime curves. The existence of infinitely many geometrically distinct (in the sense of
the definition given in the Introduction) closed geodesics in $M$ is equivalent to the existence of infinitely many distinct
prime critical $\mathrm O(2)$-orbits of $f$ in $\widetilde{\mathcal N}$.
\end{rem}

It will be useful to prove the following two results:
\begin{lem}\label{thm:isolatedcriticalorbits}
If $(M,g)$ has only a finite number of geometrically distinct closed geodesics, then the critical
orbits of $f$ in $\widetilde{\mathcal N}$ are isolated.
\end{lem}
\begin{proof}
If $\gamma_1,\ldots,\gamma_r$ is a maximal set of pairwise geometrically distinct closed geodesics
in $M$ with $\gamma_j(0)\in S$ for all $j$, then the critical orbits of $f$ in $\widetilde{\mathcal N}$ is the countable set
formed by all the iterates $\mathrm O(2)\big[\gamma_j^{(N)}\big]$, $j=1,\ldots,r$, $N\ge1$; observe that $f(\gamma_j)>0$ for all
$j$. Any sequence $k\mapsto \gamma_{j_k}^{(N_k)}$ of pairwise distinct iterates of the $\gamma_j$'s
would necessarily have $N_k\to\infty$, hence $f(\gamma_{j_k}^{(N_k)})\to+\infty$.
In particular, no subsequence of such sequence could have a converging subsequence in $\widetilde{\mathcal N}$.
The group $\mathrm O(2)$ is compact, and the conclusion follows easily.
\end{proof}
Let $S$ be a Cauchy surface in $(M,g)$;
we will use the identification \eqref{eq:tildecalN} to prove the existence of a strong deformation
retract from the $\varepsilon$-sublevel of $f$ in $\widetilde{\mathcal N}$ to the set of constant curves in $S$.
\begin{lem}\label{thm:defretract}
For $\varepsilon>0$ small enough, the closed $\varepsilon$-sublevel of $f$ in $\widetilde{\mathcal N}$:
\[f^\varepsilon=\big\{[\gamma]\in\widetilde{\mathcal N}:f(\gamma)\le\varepsilon\big\}\]
is homotopically equivalent to (the set of constant curves in) $S$.
\end{lem}
\begin{proof}
Let us show that the map $f^\varepsilon\ni\gamma\mapsto\gamma(0)\in S$ is a deformation retract.
By part \eqref{itm:Ycomplete} of Proposition~\ref{thm:constrainedvarprob}, it suffices to show that there
exists a continuous map $\Phi:f^\varepsilon\times[0,1]\to \Lambda M$ with $\Phi(\gamma,0)=\gamma$ and $\Phi(\gamma,1)$ equal to
the constant curve $\gamma(0)$.
To this aim, consider the auxiliary Riemannian metric $h$ on $M$ defined by $h(v,w)=g(v,w)-2g(v,\nabla T)g(w,\nabla T)g(\nabla T,\nabla T)^{-1}$.
Recalling  that the functions $g(\mathcal Y,\mathcal Y)$ and $g(\nabla T,\nabla T)$ admit minimum in $M$,
set $a_0=\min\big[-g(\mathcal Y,\mathcal Y)\big]>0$ and $b_0=\min\big[-g(\nabla T,\nabla T)\big]>0$.
From \eqref{eq:secondadisug}, if $[\gamma]\in f^\varepsilon$, then:
\[\int_0^1g(\dot\gamma,\nabla T)^2\,\mathrm dt\le2\varepsilon a_0^{-1},\]
and thus:
\begin{equation}\label{eq:stimaenergiah}
\int_0^1h(\dot\gamma,\dot\gamma)\,\mathrm dt=\int_0^1\Big[g(\dot\gamma,\dot\gamma)-2g(\dot\gamma,\nabla T)^2g(\nabla T,\nabla T)^{-1}\Big]\,
\mathrm dt\le\varepsilon\left(1+\frac2{a_0b_0}\right).
\end{equation}
Using the Cauchy--Schwartz inequality, we get that the $h$-length of every curve $\gamma\in f^\varepsilon$ is less than
or equal to $\varepsilon\left(1+\frac2{a_0b_0}\right)$.
Let $\rho_0>0$ be the minimum on the compact manifold $S$ of the radius of injectivity of the Riemannian metric $h$;
choose a positive $\varepsilon<\rho_0\left(\frac{a_0b_0}{a_0b_0+2}\right)$; if $\gamma$ is a curve in $f^\varepsilon$ and
$t\in[0,1]$, then the $h$-distance between $\gamma(t)$ and $\gamma(0)$ is less than $\rho_0$.
The required deformation retract $\Phi$ is given by setting $\Phi(\gamma,s)(t)=c(s)$, where
$c:[0,1]\to M$ is the unique affinely parameterized minimal $h$-geodesic from $\gamma(0)$ to $\gamma(t)$.
\end{proof}

\end{section}

\begin{section}{The Morse index theorem}\label{sec:indexthm}
In this section we will prove an index theorem for closed  geodesics
in a stationary Lorentz\-ian manifold with arbitrary endpoints,
generalizing the result in \cite{MarPicTau}. The result is now obtained
as a corollary of Theorem~\ref{thm:magicbilin}, together with the semi-Riemannian
Morse index theorem for fixed endpoints geodesics proved in \cite{GiaPicPorCOMPTES}.
An earlier version of the theorem was proven in \cite{MarPicTau} for the nondegenerate case,
under the further assumption that the closed geodesic be orientation preserving.
The use of Theorem~\ref{thm:magicbilin} allows to get rid of both these extra assumptions
at the same time.
\subsection{The index theorem}
Let us consider a closed geodesic $\gamma$ in $M$;
it is easy to check that  $T_\gamma\mathcal N$ contains
the space of all Jacobi fields $J$ along $\gamma$ such that $J(0)=J(1)$.
The following lemma tells us that $\gamma$ is a nondegenerate critical point of $f$ if and only if
it is a nondegenerate critical point of $f\vert_{\mathcal N}$:
\begin{lem}\label{thm:kernelIgamma}
Let $\gamma$ be a critical point of $f\vert_{\mathcal N}$, i.e., a closed geodesic in $M$.
Then, the kernel of the index form $I_\gamma$ in $T_\gamma\mathcal N$ coincides with the Kernel of $I_\gamma$ in $T_\gamma\Lambda M$, and it
is given by the space of periodic Jacobi fields along $\gamma$:
\[\Ker\big(I_\gamma\vert_{T_\gamma\mathcal N\times T_\gamma\mathcal N}\big)=\Big\{J\ \text{Jacobi field along $\gamma$}: J(0)=J(1),\ \Ddt J(0)=\Ddt J(1)\Big\}.\]
Moreover, consider the following closed subspace $\mathcal W_\gamma\subset T_\gamma\mathcal N$:
\[\mathcal W_\gamma=\Big\{ V\in T_\gamma\mathcal N:V(0)=V(1)=0\Big\}.\]
Then, the $I_\gamma$-orthogonal space of $\mathcal W_\gamma$ in $T_\gamma\mathcal N$ is given by:
\[\mathcal S_\gamma=\Big\{J\ \text{Jacobi field along $\gamma$}: J(0)=J(1)\Big\}.\]
\end{lem}
\begin{proof}
The statement on the kernel of $I_\gamma$ is proved readily using the following two facts:
\begin{itemize}
\item[(a)] $T_\gamma\Lambda M=T_\gamma\mathcal N+\mathfrak Y$, where $\mathfrak Y$ is the space of vector fields in $T_\gamma\Lambda M$ that are pointwise
multiple of the Killing field $\mathcal Y$;
\item[(b)] $\mathfrak Y$ is contained in the $I_\gamma$-orthogonal complement of $T_\gamma\mathcal N$ in $T_\gamma\Lambda M$.
\end{itemize}
In order to prove (a), simply observe that, for any $W\in T_\gamma\Lambda M$, then the vector field $V$ along $\gamma$ defined below
belongs to $T_\gamma\mathcal N$:
\[V(t)=W(t)+\lambda_W(t)\cdot\mathcal Y\big(\gamma(t)\big),\quad t\in[0,1],\]
where
\[\lambda(t)=\int_0^t\frac{C_W+g\big(W,\Ddt\mathcal Y\big)-g\big(\Ddt W,\mathcal Y\big)}{g(\mathcal Y,\mathcal Y)}\,\mathrm ds,\]
and
\[C_W=\left[\int_0^1\frac{g\big(\Ddt W,\mathcal Y\big)-g\big(W,\Ddt\mathcal Y\big)}{g(\mathcal Y,\mathcal Y)}\,\mathrm ds\right]\cdot\left(
\int_0^1\frac{\mathrm ds}{g(\mathcal Y,\mathcal Y)}\right)^{-1}.\]
Part (b) is a simple partial integration calculation, which is omitted; similarly, the last part of the statement
is obtained by an immediate calculation using the fundamental Lemma of calculus of variations.
\end{proof}

\begin{rem}
In the case of a periodic geodesic $\gamma$ the index form $I_\gamma$ is always degenerate, being the tangent field
$\dot\gamma$ in its kernel. Moreover, also the restriction
of the Killing field $\mathcal Y$ to $\gamma$ is a non trivial Jacobi field in
$\Ker\big(I_\gamma\vert_{T_\gamma\mathcal N \times T_\gamma\mathcal N}\big)$.
Thus, $\Dim\big[\Ker\big(I_\gamma\vert_{T_\gamma\mathcal N \times T_\gamma\mathcal N}\big)\big]\ge2$.
\end{rem}
\begin{rem}\label{thm:remJacobiNtilde}
If $S$ is a Cauchy surface in $(M,g)$, using the identifications \eqref{eq:tildecalN} and \eqref{eq:TtildecalN},
the null space of the Hessian of $f\vert_{\widetilde{\mathcal N}}$ at $[\gamma]$ is given by the space of
periodic Jacobi fields $J$ along $\gamma$ such that $J(0)\in T_{\gamma(0)}S$.
The tangent space $T_{[\gamma]}\big(\mathrm O(2)[\gamma]\big)$ is given by the space of all
constant multiples of the periodic Jacobi vector field $J$ along $\gamma$ given by $J(t)=\dot\gamma(t)+\alpha\mathcal Y\big(\gamma(t)\big)$,
where $\alpha\in\R$ is such that $J(0)\in T_{\gamma(0)}S$.
\end{rem}

By Lemma~\ref{thm:kernelIgamma}, the nullity $\mathrm n(\gamma)$ is equal to the dimension of the space of
periodic Jacobi fields $J$ along $\gamma$.

\begin{teo}[Morse index theorem for closed geodesics with arbitrary endpoints]
\label{thm:MorseIndexTheorem}
Let $\gamma:[0,1]\to M$ be a closed geodesic in $M$. Then, the  Morse index
$\mu(\gamma)$ of $f\vert_{\mathcal N}$ at $\gamma$ is given by:
\begin{equation}\label{eq:MorseIndexTHM}
\mu(\gamma)=\iMaslov(\gamma)+1+\mathrm n_-(B_0)-n_1,
\end{equation}
where $B_0$ is the symmetric bilinear form on the finite dimensional vector space
$\mathcal S_\gamma$ given by:
\begin{equation}\label{eq:restrIgSg}B_0(J_1,J_2)=g\big(\Ddt J_1(0),J_2(0)\big),\end{equation}
and $n_1$ is the dimension of the vector space:
\[\mathcal W_\gamma\cap\Ker(I_\gamma)=\Big\{J\ \text{Jacobi field along $\gamma$}: J(0)=J(1)=0,\ \Ddt J(0)=\Ddt J(1)\Big\}.\]

\end{teo}
\begin{proof}
Formula \eqref{eq:MorseIndexTHM} follows  from Theorem~\ref{thm:magicbilin} applied to the index form $I_\gamma$ and the closed
spaces $\mathcal W_\gamma$ and $\mathcal S_\gamma$ introduced in Lemma~\ref{thm:kernelIgamma}.
One has:
 \[\mathrm n_-\big(I_\gamma\vert_{\mathcal W_\gamma\times\mathcal W_\gamma}\big)=\iMaslov(\gamma)+1-n_0,\]
where $n_0$ is the dimension of the vector space:
\[\mathcal W_\gamma\cap\mathcal S_\gamma=\big\{J\ \text{Jacobi field along $\gamma$}: J(0)=J(1)=0\big\}.\]
Such equality is given by the Morse index theorem for fixed endpoints geodesics, which is proved in \cite{asian} in
the nondegenerate case, and in \cite{GiaPicPorCOMPTES} for the general case\footnote{Recall also that the definition
of Maslov index $\iMaslov(\gamma)$ employed here differs by $1$ from the definition in \cite{GiaPicPorCOMPTES}.}. In order to apply the result
of \cite{GiaPicPorCOMPTES}, one needs to observe that the extended index (i.e., index plus nullity) of $I_\gamma$
in $\mathcal W_\gamma$ is equal to the \emph{spectral flow} of the path of Fredholm
symmetric bilinear forms $[0,1]\ni s\mapsto I_{\gamma\vert_{[0,s]}}$ defined on the space of fixed endpoints
variational vector fields along $\gamma\vert_{[0,s]}$. This follows easily from the fact that $I_\gamma$
is negative semi-definite on the space $\mathfrak Y$ defined in the proof of Lemma~\ref{thm:kernelIgamma}.
An immediate partial integration argument shows that the restriction of $I_\gamma$ to $\mathcal S_\gamma$ is given
by \eqref{eq:restrIgSg}, and equality \eqref{eq:MorseIndexTHM} follows readily.
\end{proof}
Observe that the following inequalities hold:
\begin{gather}\notag
0\le n_1\le\Dim(M),\quad 0\le\mathrm n_-(B_0)\le\Dim(M),\quad0\le\mathrm n(\gamma)\le\Dim(M),\\
\label{eq:ineqn0n1} 0\le \mathrm n_-(B_0)+n_0-n_1\le \Dim(M).
\end{gather}
Inequality \eqref{eq:ineqn0n1} is obtained easily using part~\eqref{thm:indcodindmonot} in Lemma~\ref{thm:KerBSSperp},
and observing that $\mathcal W_\gamma$ has codimension equal to $\Dim(M)$ in $T_\gamma\mathcal N$.
\subsection{Morse index of an iteration}
\label{sub:iterMorse}
Throughout this subsection, we will consider a fixed critical point $\gamma$ of $f\vert_{\mathcal N}$.
Given an integer $N\ge1$, let us denote by $\gamma^{(N)}$ the \emph{$N$-iterated of $\gamma$}, defined
by $\gamma^{(N)}(t)=\widetilde\gamma(Nt)$ for all $t\in[0,1]$, where $\widetilde\gamma:\R\to M$ is the periodic
extension of $\gamma$. Observe that $\gamma^{(N)}$ is a critical
point of $f\vert_{\mathcal N}$ for all $N\ge1$.
One of the central results of this paper will be to establish the growth of the sequence
$\mu\big(\gamma^{(N)}\big)$  (Proposition~\ref{thm:lineargrowth} and
Corollary~\ref{thm:corcrescitalinearearb}). The result will be first established for orientation
preserving closed geodesics, and then extended to the general case using Lemma~\ref{thm:nondecreasing}
below.
\smallskip

Although it is not clear at all whether the Morse index of a closed geodesic increases by iteration,
an argument using a finite codimensional restriction of the index form yields the following interesting
consequence:

\begin{lem}\label{thm:nondecreasing}
There exists a \emph{bounded} sequence of integers $(d_N)_{N\ge1}$ such that the sequence $N\mapsto\mu\big(\gamma^{(N)}\big)+d_N\in\Z$ is nondecreasing.
\end{lem}
\begin{proof}
Let us introduce the following space:
\begin{equation}\label{eq:deWgamma0}
\mathcal W_\gamma^o=\Big\{V\in\mathcal W_\gamma: g\big(\Ddt V,\mathcal Y\big)-g\big(V,\Ddt\mathcal Y\big)\equiv0\Big\}.
\end{equation}
Clearly, $\mathcal W_\gamma^o$ is a $1$-codimensional closed subspace of $\mathcal W_\gamma$, being the kernel of the
bounded linear functional $\mathcal W_\gamma\ni V\mapsto C_V\in\R$ (see \eqref{eq:defCV}).
Hence, recalling part~\eqref{thm:indcodindmonot} in Lemma~\ref{thm:KerBSSperp}:
\[\mathrm n_-\big(I_\gamma\vert_{\mathcal W_\gamma^o\times \mathcal W_\gamma^o}\big)\le \mathrm n_-\big(I_\gamma\vert_{\mathcal W_\gamma\times\mathcal W_\gamma}\big)\le
\mathrm n_-\big(I_\gamma\vert_{\mathcal W_\gamma^o\times\mathcal W_\gamma^o}\big)+1.\]
Thus, keeping in mind formula \eqref{eq:MorseIndexTHM} and inequality~\ref{eq:ineqn0n1}, in order to prove the Lemma it suffices to show that the sequence
$\bar\mu\big(\gamma^{(n)}\big)$ is nondecreasing, where
\begin{equation}\label{eq:defbarmun}
\bar\mu\big(\gamma\big)= \mathrm n_-\big(I_{\gamma}\vert_{\mathcal W_{\gamma}^o\times\mathcal W_{\gamma}^o}\big).
\end{equation}
To this aim,
let $1\le N\le M$ be given, and consider the map:
\[\mathcal E_{N,M}:\mathcal W_{\gamma^{(N)}}^o\longrightarrow \mathcal W_{\gamma^{(M)}}^o\]
defined by $\mathcal E_{N,M}(V)=\widetilde V$, where:
\[\widetilde V(t)=\begin{cases}V(tM/N),&\text{if $t\in[0, N/M]$;}\\ 0,&\text{if $t\in\left]N/M,1\right]$.} \end{cases}\]
Obviously, $\mathcal E_{N,M}$ is an injective bounded linear map; an immediate computation shows that
the following equality holds:
\begin{equation}\label{eq:pullbackindexform}
I_{\gamma^{(M)}}\big(\mathcal E_{N,M}(V),\mathcal E_{N,M}(W)\big)=\tfrac{M}N\,I_{\gamma^{(N)}}(V,W),\quad\forall\, V,W\in \mathcal W_{\gamma^{(N)}}^o.
\end{equation}
Hence, if $\mathcal V\subset \mathcal W_{\gamma^{(N)}}^o$ is a subspace such that
$\Dim(\mathcal V)=\mathrm n_-\big(I_\gamma\vert_{\mathcal W_{\gamma^{(N)}}^o\times\mathcal W_{\gamma^{(N)}}^o}\big)$ and such that
$I_{\gamma^{(N)}}$ is negative definite on $\mathcal V$, then $\Dim\big(\mathcal E_{N,M}(\mathcal V)\big)=\Dim(\mathcal V)$ and, by \eqref{eq:pullbackindexform},
$I_{\gamma^{(M)}}$ is negative definite on $\mathcal E_{N,M}(\mathcal V)$.
This shows that $\bar\mu\big(\gamma^{(N)}\big)\le\bar\mu\big(\gamma^{(M)}\big)$ and concludes the proof.
\end{proof}
\noindent It will be useful to record here the following relation between the Morse index $\mu(\gamma)$, the Maslov index $\iMaslov(\gamma)$ and
the \emph{restricted Morse index\/} $\bar\mu(\gamma)$ (see \eqref{eq:defbarmun}) of a closed geodesic~$\gamma$:
\[\bar\mu(\gamma)\le\iMaslov(\gamma)\le\mu(\gamma);\]
more precisely:
\begin{equation}\label{eq:relindices}
\begin{aligned}
&\mu(\gamma)=\iMaslov(\gamma)+A_\gamma,\quad 0\le A_\gamma\le\Dim(M)-1,\\
&\iMaslov(\gamma)=\bar\mu(\gamma)+B_\gamma,\quad 0\le B_\gamma\le1.
\end{aligned}
\end{equation}
Exploiting the same idea in Lemma~\ref{thm:nondecreasing}, one has the following result on the additivity of
the Morse index:
\begin{lem}\label{thm:lemadditivity}
There exists a bounded sequence $(e_N)_{N\ge1}$ of nonnegative integers such that for all $r,s>0$, the following inequality holds:
\begin{equation}\label{eq:additivita}
\mu\big(\gamma^{(r+s)}\big)\ge\mu(\gamma^{(r)}\big)+\mu(\gamma^{(s)}\big)-e_r-e_s.
\end{equation}
\end{lem}
\begin{proof}
As in the proof of Lemma~\ref{thm:nondecreasing}, the sequence $\bar\mu\big(\gamma^{(N)}\big)$ satisfies:
\[\bar\mu\big(\gamma^{(r+s)}\big)\ge\bar\mu(\gamma^{(r)}\big)+\bar\mu(\gamma^{(s)}\big);\]
the conclusion follows easily using \eqref{eq:relindices} and setting $e_N=A_{\gamma^{(N)}}+B_{\gamma^{(N)}}\le\Dim(M)$.
\end{proof}
Finally, we have our aimed results on the growth of the Maslov index:
\begin{prop}\label{thm:lineargrowth}
Given any  closed geodesic $\gamma$ in $M$, the sequence of Morse indices $N\mapsto\mu\big(\gamma^{(N)}\big)$ is either bounded
(by a constant depending only on the dimension of $M$),
or it has superlinear growth in $N$ for large $N$.
\end{prop}
\begin{proof}
Assume first that $\gamma$ is orientation preserving, and that $\mu\big(\gamma^{(N)}\big)$ is not bounded.
Let $k_*\in\N$ be the first positive integer such that:
\[\mu\big(\gamma^{(k_*)}\big)>8\,\Dim(M)+1.\]
Using Theorem~\ref{thm:MorseIndexTheorem}, the noncreasing property of the restricted Morse index proved in
Lemma~\ref{thm:nondecreasing} and formulas \eqref{eq:relindices}, for $m\ge k_*$, we compute as follows:
\begin{multline*}
\mu\big(\gamma^{(m)}\big)=\bar\mu\big(\gamma^{(m)}\big)+A_{\gamma^{(m)}}+B_{\gamma^{(m)}}
\stackrel{\text{Lemma~\ref{thm:nondecreasing}}}\ge \bar\mu\big(\gamma^{(\lfloor \frac m{k_*}\rfloor k_*)}\big)+A_{\gamma^{(m)}}+B_{\gamma^{(m)}}
\\ \stackrel{\text{by \eqref{eq:relindices}}}\ge\iMaslov\big(\gamma^{(\lfloor \frac m{k_*}\rfloor k_*)}\big)-1+A_{\gamma^{(m)}}+B_{\gamma^{(m)}}
\\ \stackrel{\text{Corollary~\ref{thm:iterateMaslovgeo}}}\ge\Big(\iMaslov\big(\gamma^{(k_*)}\big)-7\,\Dim(M)\Big)\cdot\left\lfloor\frac m{k_*}\right\rfloor-5\,\Dim(M)
-1+A_{\gamma^{(m)}}+B_{\gamma^{(m)}}\\
\stackrel{\text{by \eqref{eq:relindices}}}\ge\bigg(\frac{\mu\big(\gamma^{(k_*)}\big)-8\,\Dim(M)-1}{k_*}\bigg)\cdot m-\mu\big(\gamma^{(k_*)}\big)+2\,\Dim(M)
-1+A_{\gamma^{(m)}}+B_{\gamma^{(m)}}.
\end{multline*}
Here, $\lfloor\cdot\rfloor$ denotes the integer part function. The conclusion follows, recalling from formulas
\eqref{eq:relindices} that $A_{\gamma^{(m)}}$ and $B_{\gamma^{(m)}}$ are bounded sequences.

For the general case of possibly non orientation preserving closed geodesics, observe that the
double iterate $\gamma^{(2)}$ of any closed geodesic is orientation preserving. Observe also
that, by Lemma~\ref{thm:nondecreasing}, the sequence $\mu\big(\gamma^{(N)}\big)$ is bounded if
and only if $\mu\big(\gamma^{(2N)}\big)$ is bounded. Based on these observations and on Lemma~\ref{thm:nondecreasing},
establishing the superlinear growth of $\mu\big(\gamma^{(N)}\big)$ in the non orientable case is obtained
by elementary arithmetics from the previous case.
\end{proof}
We will need a slightly refined property on the growth of the Morse index, which
is some sort of uniform superlinear growth:
\begin{cor}\label{thm:corcrescitalinearearb}
Let $\gamma$ be a closed geodesic in $M$ such that $\mu\big(\gamma^{(N)}\big)$ is not bounded.
Then, there exist positive constants $\bar\alpha,\bar\beta\in\R$,
such that, for $s$ sufficiently large, the following inequalities hold:
\begin{equation}\label{eq:corcrescitalinearearb}
\mu\big(\gamma^{(r+s)}\big)\ge\mu\big(\gamma^{(r)}\big)+s\,\bar\alpha-\bar\beta,\qquad\forall\,r>0.
\end{equation}
\end{cor}
\begin{proof}
Let $k_*$ be as in the proof of Proposition~\ref{thm:lineargrowth}, and set \[\bar\alpha=\frac{\mu\big(\gamma^{(k_*)}\big)-8\,\Dim(M)-1}{k_*},
\qquad \bar\beta=\mu\big(\gamma^{(k_*)}\big)+1.\] For $s\ge k_*$, inequality \eqref{eq:corcrescitalinearearb} follows readily from
Lemma~\ref{thm:lemadditivity} and Proposition~\ref{thm:lineargrowth}.
\end{proof}
\subsection{Nullity of an iteration}\label{sub:nullity}
The nullity of an iterated closed geodesic $\gamma$ will be computed using the spectrum of the linearized Poincar\'e map
$\mathfrak P_\gamma$ defined below.
Given a closed geodesic $\gamma:[0,1]\to M$, denote by $\mathcal V$ the space $T_{\gamma(0)}M\oplus T_{\gamma(0)}M^*$,
endowed with its canonical symplectic structure, and let $\mathfrak P_\gamma:\mathcal V\to\mathcal V$ be the linear map defined
by:
\[\mathfrak P_\gamma\big(J(0),g\Ddt J(0)\big)=\big(J(1),g\Ddt J(1)\big),\]
where $J$ is a Jacobi field along $\gamma$. The map $\mathfrak P_\gamma$ is a symplectomorphism of $\mathcal V$; denote by
$\mathfrak s(\mathfrak P_\gamma)$ its spectrum.
It follows from Lemma~\ref{thm:kernelIgamma} that $\Ker\big(I_\gamma\vert_{T_\gamma\mathcal N\times T_\gamma\mathcal N}\big)$
consists of all Jacobi fields $J$ along $\gamma$ such that $\big(J(0),g\Ddt J(0)\big)$ belongs to the $1$-eigenspace
of $\mathfrak P_\gamma$.
The subspace of $\mathcal V$ spanned by $\big(\gamma'(0),0\big)$ and by $\big(\mathcal Y(\gamma(0)),
g\nabla_{\gamma'(0)}\mathcal Y\big)$ is a $2$-dimensional isotropic subspace of $\Ker(\mathfrak P_\gamma)$.
From Proposition~\ref{thm:PS}, it follows that $\mathrm O(2)[\gamma]$ is a nondegenerate critical
orbit of $f$ in $\widetilde{\mathcal N}$ when  $\Dim\big[\Ker\big(I_\gamma\vert_{T_\gamma\mathcal N\times T_\gamma\mathcal N}\big)\big]=2$.
We have a result on the nullity of an iteration, which is totally analogous to the Riemannian case (see \cite[Lemma~2]{GroMey2}
and \cite[Proposition~4.2.6]{Kli}); its proof, repeated here for the reader's convenience, is purely arithmetical.
\begin{lem}\label{thm:nuliteration}
Let $\gamma$ be a closed geodesic in $M$ and let $\gamma^{(N)}$ denote its $N$-th iterate, $N\ge1$. Then,
$\mathrm O(2)\big[\gamma^{(N)}\big]$ is a nondegenerate critical orbit of $f$ in $\widetilde{\mathcal N}$
if and only if:
\begin{itemize}
\item[(a)] $\mathrm O(2)[\gamma]$ is a nondegenerate critical orbit of $f$ in $\widetilde{\mathcal N}$;
\item[(b)] $\mathfrak s(\mathfrak P_\gamma)\setminus\{1\}$ does not contain any $N$-th root of unity.
\end{itemize}
Moreover, there exists a sequence $m_1,\ldots,m_s$ of positive integers, $s\le 2^{\Dim(M)}$, and, for each
$j\in\{1,\ldots,s\}$, a strictly increasing sequence $q_{j1}<q_{j2}<\ldots<q_{jm}<\ldots$ of positive
integers such that the sets $N_j=\big\{m_jq_{ji},i=1,2,\ldots\}$ form a partition of $\N\setminus\{0\}$, and such that
\begin{equation}\label{eq:nullitauguali}
\mathrm n\big(\gamma^{m_jq_{ji}}\big)=\mathrm n\big(\gamma^{m_j}),\quad\forall\,i\in\N.
\end{equation}
\end{lem}
\begin{proof}
The first statement is proved easily observing that $\mathfrak P_\gamma^N=\mathfrak P_{\gamma^{(N)}}$.

For the second statement, consider all the elements in $\mathfrak s(\mathfrak P_\gamma)$ of the form $e^{\pm2\pi\frac pqi}$, with
$p,q$ positive integers and relatively prime. Let $D$ the possibly empty set of all these denominators, and for all $E\subset D$
denote by $m(E)$ the least common multiple of all elements of $E$, setting $m(\emptyset)=1$. Denote by
$m_1,\ldots,m_s$ the set of all pairwise distinct numbers obtained as $m(E)$, for all subsets $E\subset D$, where
$m_1=1$.
Clearly, $s\le2^{\Dim(M)}$.  Finally, for all $j\in\{1,\ldots,s\}$, consider a maximal sequence
$\{q_{ji},i\ge1\}$ of positive integers such that none of the $m_k$, with $k\ne j$, divides $m_jq_{ji}$.
Then, \eqref{eq:nullitauguali} holds; furthermore, every $m\in\N\setminus\{0\}$ can be written as
the product $m_jq$, where $q$ is a positive integer, and $m_j$ is some divisor of $m$ among the
elements $m_1,\ldots,m_s$. If $m_j$ is the maximum of such divisors, then $q$ must be one of the $q_{ji}$'s,
for some $i\ge1$. This concludes the proof.
\end{proof}

\begin{rem}\label{thm:remfiniternul}
By Lemma~\ref{thm:nuliteration}, we have the following situation. Assuming that there is only a finite number of geometrically distinct
closed geodesics in $M$, it is possible to find a finite number of closed geodesics $\gamma_1,\ldots,\gamma_r$
(possibly not all geometrically distinct) such that any closed geodesics $\gamma$ in $M$ is geometrically
equivalent to some iterate $\gamma_{i_0}^{(N)}$ of one of the $\gamma_i$'s, and it has the same nullity as $\gamma_{i_0}$.
\end{rem}
\end{section}

\begin{section}{Equivariant Morse theory for the action functional}
\label{sec:proof}
\subsection{Abstract Morse relations}Given sequences $(\mu_k)_{k\ge0}$ and $(\beta_k)_{k\ge0}$ in $\N\bigcup\{+\infty\}$, we will say that the sequence of pairs
$(\mu_k,\beta_k)_{k\ge0}$ satisfies the \emph{Morse relations} if there exists a formal power series $Q(t)=\sum_{k\ge0}q_k\,t^k$ with
coefficients in $\N\bigcup\{+\infty\}$ such that:
\[\sum_{k\ge0}\mu_k\,t^k=\sum_{k\ge0}\beta_k\,t^k+(1+t)Q(t).\]
This condition implies (and, in fact it is equivalent to if all $\mu_k$'s are finite)
the familiar set of inequalities:
\begin{eqnarray*}
\mu_0&\ge&\beta_0,\\
\mu_1-\mu_0&\ge&\beta_1-\beta_0\\
\mu_2-\mu_1+\mu_0&\ge&\beta_2-\beta_1+\beta_0,\\
&\vdots&\\
\mu_k-\mu_{k-1}+\cdots+(-1)^k\mu_0&\ge&\beta_k-\beta_{k-1}+\cdots+(-1)^k\beta_0,
\\
&\ldots&
\end{eqnarray*}
that are called the \emph{strong Morse inequalities}. In turn, these inequalities imply the \emph{weak Morse inequalities}:
\begin{equation}\label{eq:weakMorserel}\phantom{\quad\forall\,k\ge0.}\mu_k\ge\beta_k,\quad\forall\,k\ge0.\end{equation}
Given a pair $Y\subset X$ of topological space and a coefficient field $\mathds K$, let us denote by
$H_k(X,Y;\mathds K)$ the $k$-th relative homology vector space with coefficients in $\mathds K$, and
by $\beta_k(X,Y;\mathds K)=\Dim\big(H_k(X,Y;\mathds K)\big)$ the $k$-th Betti number of the pair.
We set $H_k(X;\mathds K)=H_k(X,\emptyset;\mathds K)$ and $\beta_k(X;\mathds K)=\beta_k(X,\emptyset;\mathds K)$.
Using standard homological techniques, one proves the following:
\begin{prop}\label{thm:Morseinequalities}
Let $\mathds K$ be a field, and let $(X_n)_{n\ge0}$ be a filtration of a topological space $X$; assume that every
compact subset of $X$ is contained in some $X_n$. Setting:
\[\phantom{,\qquad k\ge0,}\mu_k=\sum_{n=0}^\infty\beta_k(X_{n+1},X_n;\mathds K),\qquad k\ge0,\]
and $\beta_k=\beta_k(X,X_0;\mathds K)$, then the sequence $(\mu_k,\beta_k)_{k\ge0}$ satisfies the
Morse relations.\qed
\end{prop}

\subsection{Homological invariants at isolated critical points and critical orbits}
\label{sub:hominvariants}
Let us recall here a few basic facts on the homological invariants associated to isolated critical points
and group orbits;
the basic references are \cite{Cha, GroMey1, GroMey2, GroTan, Wang}. Let $\mathcal M$ be a smooth Hilbert manifold and let
$\mathfrak f:\mathcal M\to\R$ be a smooth function; for $d\in\R$, denote by $\mathfrak f^d$ the closed sublevel $\big\{x\in\mathcal M:\mathfrak f(x)\le d\big\}$.
Let $p\in\mathcal M$ be a critical point of
$\mathfrak f$, and assume that the Hessian $\mathrm H^{\mathfrak f}(p)$ of $\mathfrak f$ at $p$ is represented
by a compact perturbation of the identity of $T_p\mathcal M$. A generalized Morse Lemma for this situation (\cite[Lemma~1]{GroMey1})
says that there exists a smooth local parametrization of $\mathcal M$ around $p$, $\Phi:U\to V$, where
$U$ is an open neigborhood of $0\in T_p\mathcal M\cong\Ker\big(\mathrm H^{\mathfrak f}(p)\big)^\perp\oplus \Ker\big(\mathrm H^{\mathfrak f}(p)\big)$,
$V$ is an open neighborhood of $p$, with $\Phi(0)=p$, and there exists an orthogonal projection $P$ on $\Ker\big(\mathrm H^{\mathfrak f}(p)\big)^\perp$ such that
$\mathfrak f\circ\Phi(x,y)=\Vert Px\Vert^2-\Vert(1-P)x\Vert^2+\mathfrak f_0(y)$, where $\mathfrak f_0:U\cap \Ker\big(\mathrm H^{\mathfrak f}(p)\big)\to\R$ is a smooth
function having $0$ as an isolated completely degenerate critical point.
Using this decomposition of $\mathfrak f$, a homological invariant $\mathfrak H(\mathfrak f,p;\mathds K)$ of $\mathfrak f$ at $p$ is defined by:
\[\mathfrak H(\mathfrak f,p;\mathds K)=H_*(W_p,W_p^-;\mathds K),\]
where $\mathds K$ is any coefficient field, and $(W,W^-)$ is a pair of topological spaces constructed in \cite{GroMey1} and
called \emph{admissible pair} (a \emph{GM-pair} in the language of \cite{Wang}). Let us describe briefly
such construction.
Denote by $\eta:\R\times\mathcal M\to\mathcal M$ the flow of $-\nabla\mathfrak f$ and set $\mathfrak f(p)=c$; an admissible pair $(W_p,W_p^-)$
is characterized by the following properties (see \cite[Definition~2.3]{Wang}):
\begin{enumerate}
\item $W_p$ is a closed neighborhood of $p$ that contains a unique critical point of $\mathfrak f$ and such that:
\begin{itemize}
\item[(a)] if $t_1<t_2$ and $\eta(t_i,x)\in W$ for $i=1,2$, then $\eta(t,x)\in W$ for all $t\in[t_1,t_2]$;
\item[(b)] there exists $\varepsilon>0$ such that $\mathfrak f$ has no critical value in $\left[c-\varepsilon,c\right[$
and such that $W\cap \mathfrak f^{c-\varepsilon}=\emptyset$;
\end{itemize}
\item $W^-=\big\{x\in W:\eta(x,t)\in W,\ \forall\,t>0\big\}$ is closed in $W$;
\item $W^-$ is a (piecewise smooth) hypersurface of $\mathcal M$ which is transversal to $\nabla\mathfrak f$.
\end{enumerate}
By \cite[Theorem~2.1]{Wang}, if $(W_p,W_p^-)$ is an admissible pair, then:
\[H_*(W_p,W_p^-;\mathds K)=H_*(\mathfrak f^c,\mathfrak f^c\setminus\{p\};\mathds K);\]
furthermore, by excision, if $U$ is any open subset of $\mathcal M$ containing $p$, then:
\[\mathfrak H_*(\mathfrak f,p;\mathds K)=H_*\big(U\cap\mathfrak f^c,U\cap(\mathfrak f^c)\setminus\{p\});\mathds K\big).\]
If $\mathcal M$ is complete, $c$ is the only critical value of $\mathfrak f$ in $[c-\varepsilon,c+\varepsilon]$, and
$p_1,\ldots,p_r$ are the critical points of $\mathfrak f$ in $\mathfrak f^{-1}(c)$, then the relative homology
$H_*\big(\mathfrak f^{c+\varepsilon},\mathfrak f^{c-\varepsilon};\mathds K\big)$ can be computed as:
\[H_*\big(\mathfrak f^{c+\varepsilon},\mathfrak f^{c-\varepsilon};\mathds K\big)=\bigoplus_{i=1}^r\mathfrak H_*(\mathfrak f,p_i;\mathds K).\]
Another homological invariant $\mathfrak H^o(\mathfrak f,p;\mathds K)$ is defined by setting:
\[\mathfrak H^o(\mathfrak f,p;\mathds K)=\mathfrak H(\mathfrak f_0,p;\mathds K),\]
where $\mathfrak f_0$ is the degenerate component of $\mathfrak f$ described above.
Among the main results of \cite{GroMey1}, the celebrated \emph{shifting theorem} gives a relation between $\mathfrak H(\mathfrak f,p;\mathds K)$
and $\mathfrak H^o(\mathfrak f,p;\mathds K)$. The shifting theorem states that if $\mu(p)$ is the Morse index of $\mathfrak f$ at $p$, then:
\begin{equation}\label{eq:uguaShTh}\phantom{,\qquad\forall\,k\in\mathds Z.}\mathfrak H_{k+\mu(p)}(\mathfrak f,p;\mathds K)=
\mathfrak H^o_k(\mathfrak f,p;\mathds K),\qquad\forall\,k\in\mathds Z.\end{equation}
The homological invariant $\mathfrak H$, as well as $\mathfrak H^o$, is of finite type, i.e., $\mathfrak H_k$ is finite dimensional
for all $k$ and $\mathfrak H_k=\{0\}$ except for a finite number of $k$'s.  Moreover, the homological invariant $\mathfrak H^o$ has the following localization property
\begin{lem}\label{thm:localizationproperty}
Let $\mathcal M$ be a smooth Hilbert manifold, $\mathfrak f:\mathcal M\to R$ be a smooth map, $p\in \mathcal M$ an isolated critical point
of $\mathfrak f$ such that the Hessian $\mathrm H^{\mathfrak f}(p)$ is represented by compact perturbation of the identity.
Let $\widehat{\mathcal M}$ be a smooth closed submanifold of $\mathcal M$ containing $p$ such that $\nabla \mathfrak f_q\in T_q\widehat{\mathcal M}$ for
all $q\in \widehat M$, and such that the null space of the Hessian $\mathrm H^{\mathfrak f}(p)$ is contained in $T_p\widehat M$.
Then, $\mathfrak H^o(\mathfrak f,p)=\mathfrak H^o\big(\mathfrak f\vert_{\widehat{\mathcal M}},p\big)$.
\end{lem}
\begin{proof}
See \cite[Lemma~7, p.\ 368--369]{GroMey1}.
\end{proof}

Consider now the case of a compact Lie group $G$ acting by isometries on $\mathcal M$, and let $\mathfrak f:\mathcal M\to\R$ be
a $G$-invariant smooth function satisfying the Palais--Smale condition. If $p$ is a critical point of $\mathfrak f$, denote
by $Gp$ its $G$-orbit, which consists of critical points of $\mathfrak f$. If such critical orbit is isolated, i.e.,
if there exists an open neighborhood of $Gp$ that does not contain critical points of $\mathfrak f$ oustide $Gp$, then
one defines a homological invariant at the critical orbit $Gp$ by setting:
\[\mathfrak H(\mathfrak f,Gp;\mathds K)=H_*(\mathfrak f^c,\mathfrak f^c\setminus Gp;\mathds K), \]
where $c=\mathfrak f(p)$. Again, by excision, if $U$ is any open subset of $\mathcal M$ containing $Gp$, then:
\[\mathfrak H_*(\mathfrak f,Gp;\mathds K)=H_*\big(U\cap\mathfrak f^c,U\cap(\mathfrak f^c)\setminus Gp);\mathds K\big).\]
If $\mathcal M$ is complete, $c$ is the unique critical value of $\mathfrak f$ in $[c-\varepsilon,
c+\varepsilon]$, and the critical set of $\mathfrak f$ at $c$ consists of a finite number of isolated critical
orbits $Gp_1$, \ldots, $Gp_r$, then by \cite[Theorem~2.1]{RobSal}, the relative homology $H_*(\mathfrak f^{c+\varepsilon},\mathfrak f^{c-\varepsilon};
\mathds K)$ can be computed as:
\begin{equation}\label{eq:homrelsublevcrorb}
H_*(\mathfrak f^{c+\varepsilon},\mathfrak f^{c-\varepsilon};
\mathds K)=\bigoplus_{i=1}^r\mathfrak H_*(\mathfrak f,Gp_i;\mathds K).
\end{equation}

\subsection{Local homological invariants at critical $\mathrm O(2)$-orbits in $\widetilde{\mathcal N}$}
Let us now consider the Hilbert manifold $\widetilde{\mathcal N}$ \eqref{eq:tildecalN}
and the geodesic action functional $f:\widetilde{\mathcal N}\to\R$.
Consider a non constant critical point $[\gamma]$ of $f$ and assume that the critical orbit $\mathrm O(2)[\gamma]$ is isolated.
Recalling that the Hessian of $f$ at each critical orbit is a Fredholm form which is a compact perturbation of
the identity (part \eqref{thm:constrainedvarprobitm5} of Proposition~\ref{thm:constrainedvarprob}), the completeness of $\widetilde{\mathcal N}$
and the Palais--Smale condition (Proposition~\ref{thm:PS}), the construction of the local homological invariant
at the critical orbit $\mathrm O(2)[\gamma]$ can be performed as follows.
Denote by $\Gamma\subset\mathrm{SO}(2)$ the stabilizer of $\gamma$, which is a finite cyclic group;
observe that the quotient $\mathrm O(2)/\Gamma\cong\mathrm O(2)[\gamma]$ is diffeomorphic to the union of two copies of the circle
and denote by $\nu\big(\mathrm O(2)[\gamma]\big)\subset T\widetilde{\mathcal N}$ the normal bundle of $\mathrm O(2)[\gamma]$ in
$\widetilde{\mathcal N}$. Denote by $\mathrm{EXP}$ the exponential map of $\widetilde{\mathcal N}$ relatively to
the metric \eqref{eq:metricaHilbert}, and let $r>0$ be chosen small enough so that $\mathrm{EXP}$ gives a diffeomorphism between:
\[\mathcal A_r=\Big\{v\in\nu\big(\mathrm O(2)[\gamma]\big):\Vert v\Vert<r\Big\}\]
and an open subset $\mathcal D$ of $\widetilde{\mathcal N}$ containing $\mathrm O(2)[\gamma]$.
For $u\in\mathrm O(2)[\gamma]$, set \[\mathcal D^u=\mathrm{EXP}_u\big(\mathcal A_r\cap T_u\widetilde{\mathcal N}\big);\]
$\mathcal D$ is a normal disc bundle over $\mathrm O(2)[\gamma]$ whose fiber at $u$ is $\mathcal D^u$.
Observe that, since $\mathrm O(2)$ acts by isometries on $\widetilde{\mathcal N}$, then for all
$g\in\mathrm O(2)$ and all $u\in\mathrm O(2)[\gamma]$, $g\mathcal D^u=\mathcal D^{gu}$.
In particular, the restriction of the $\mathrm O(2)$-action gives an action of $\Gamma$ on each fiber $\mathcal D^u$.

Consider the principal fiber bundle $\mathrm O(2)\mapsto\mathrm O(2)/\Gamma\cong\mathrm O(2)[\gamma]$;
we claim that the bundle $\mathcal D$ can be described as the fiberwise product:\footnote{%
Recall that given a $\mathcal G$-principal fiber bundle $\mathcal P\to \mathcal X$ over the manifold $\mathcal X$, and given a
topological space $\mathcal Y$ endowed with a left $\mathcal G$-action, the fiberwise product $\mathcal P\times_{\mathcal G}\mathcal Y$
is a fiber bundle over $\mathcal X$ whose fiber at $x\in\mathcal X$ is the quotient of the product $\mathcal P_x\times\mathcal Y$
by the left action of $\mathcal G$ given by: \[\mathcal G\times(\mathcal P_x\times\mathcal Y)\ni\big(g,(p,y)\big)=(pg^{-1},gy)\in \mathcal P_x\times\mathcal Y.\]
Since the right action of $\mathcal G$ on $\mathcal P_x$ is free and transitive, then each fiber of $\mathcal P\times_{\mathcal G}\mathcal Y$ is
homeomorphic to $\mathcal Y$. Fiberwise products are examples of associated bundles to principal bundles.}
\begin{equation}\label{eq:strufiberprod}
\mathcal D\cong\mathrm O(2)\times_\Gamma\mathcal D^{\gamma}.
\end{equation}
Namely, consider the local diffeomorphism $\psi:\mathrm O(2)\times\mathcal D^{\gamma}\to\mathcal D$:
\[\mathrm O(2)\times\mathcal D^\gamma\ni (g,\sigma)\longmapsto g\sigma\in \mathcal D; \]
assuming $\psi(g,\sigma)=\psi(g',\sigma')$ gives $g^{-1}g'=h\in\Gamma$, and $h\sigma=\sigma'$, thus $(g',\sigma')=(gh^{-1},h\sigma)$
and $\psi$ passes to the quotient giving a diffeomorphism $\bar\psi:\mathrm O(2)\times_\Gamma\mathcal D^\gamma\to\mathcal D$,
and \eqref{eq:strufiberprod} is proved.
As observed above, by excision, the local homological invariant $\mathfrak H\big(f,\mathrm O(2)[\gamma];\mathds K\big)$
can be computed as:
\[\mathfrak H_*\big(f,\mathrm O(2)[\gamma];\mathds K\big)=H_*\big(f^c\cap\mathcal D,(f^c\setminus\mathrm O(2)[\gamma])\cap\mathcal D;\mathds K\big),\]
where $c=f(\gamma)$. Since $f$ is $\mathrm O(2)$-invariant, with this construction, we have $g(f^c\cap\mathcal D^\gamma)=f^c\cap
\mathcal D^{g\gamma}$ for all $g\in\mathrm O(2)$; in particular, $f^c\cap\mathcal D^\gamma$ is $\Gamma$-invariant,
and we have two fiber bundles over $\mathrm O(2)$:
\begin{equation}\label{eq:pairofbundles}
f^c\cap\mathcal D=\mathrm O(2)\times_\Gamma(f^c\cap \mathcal D^\gamma),
\quad (f^c\cap\mathcal D)\setminus\mathrm O(2)[\gamma]=\mathrm O(2)\times_\Gamma\big((f^c\cap\mathcal D^\gamma)\setminus\{[\gamma]\}\big).
\end{equation}
 If $c$ is the only critical value of $f$ in $[c-\varepsilon,c+\varepsilon]$, and
$\mathrm O(2)\big[\gamma_1\big],\ldots,\mathrm O(2)\big[\gamma_1\big]$ are the critical orbits of $f$ in $f^{-1}(c)$, then by \eqref{eq:homrelsublevcrorb}
the relative homology $H_*\big(f^{c+\varepsilon},f^{c-\varepsilon};\mathds K\big)$ is given by:
\begin{equation}\label{eq:sommauguale}
H_*\big(f^{c+\varepsilon},f^{c-\varepsilon};\mathds K\big)=\bigoplus_{i=1}^r\mathfrak H_*\big(f,\mathrm O(2)\big[\gamma_i\big];\mathds K\big).
\end{equation}
\begin{rem}\label{thm:remhessianorestriction}
The restriction $f\vert_{\mathcal D^\gamma}$ of $f$ to the disc $\mathcal D^\gamma$ has an isolated critical
point at $[\gamma]$. By \eqref{thm:constrainedvarprobitm5} of Proposition~\ref{thm:constrainedvarprob},
the Hessian $\mathrm H^{f\vert_{\mathcal D^\gamma}}$ at $[\gamma]$ of the restriction $f\vert_{\mathcal D^\gamma}$ is
essentially positive (see part \eqref{itm:RCPPI4} in Lemma~\ref{thm:propRCPPI}).
We can therefore define the local homological invariant $\mathfrak H(f\vert_{\mathcal D^\gamma},[\gamma];\mathds K)$
as the relative homology \[\mathfrak H(f\vert_{\mathcal D^\gamma},[\gamma];\mathds K)=H_*(f^c\cap \mathcal D^\gamma, (f^c\cap\mathcal D^\gamma)\setminus\{[\gamma]\};\mathds K).\]
Observe also that the Morse index of $[\gamma]$ as a critical point of the restriction
$f\vert_{\mathcal D^\gamma}$ equals the Morse index of $f$ at $[\gamma]$; the dimension of the kernel of $\mathrm H^{f\vert_{\mathcal D^\gamma}}$
at $[\gamma]$ equals the dimension of the kernel of $\mathrm H^f$ at $[\gamma]$ minus one.
\end{rem}
For all $k\ge0$, set:
\[
B_k(\gamma;\mathds K)=\Dim\big[\mathfrak H_k(f,\mathrm O(2)[\gamma];\mathds K)\big],\]
\[C_k(\gamma;\mathds K)=\Dim\big[\mathfrak H_k(f\vert_{\mathcal D^\gamma},[\gamma];\mathds K)\big],\quad
\text{and}\quad C_k^o(\gamma,\mathds K)=\Dim\big[\mathfrak H_k^o(f\vert_{\mathcal D^\gamma},[\gamma];\mathds K)\big].
\]
Our construction of the local homological invariants does not clarify that,
in fact, the invariants $C_k$ and $C_k^o$ do \emph{not} depend on the metric structure of $\widetilde{\mathcal N}$;
observe that in Proposition~\ref{thm:stessiBk0} we will need to employ different Riemannian
structures on $\widetilde{\mathcal N}$.
In order to prove the independence on the metric, we will now establish that
$C_k(\gamma;\mathds K)$ and  $C_k^o(\gamma,\mathds K)$ can be computed by considering
restrictions of $f$ to \emph{any} hypersurface $\Sigma$ of $\widetilde{\mathcal N}$ through
$[\gamma]$ which is transversal to the orbit $\mathrm O(2)[\gamma]$:
\begin{lem}\label{thm:Ckrestrictions}
Let $\mathrm O(2)[\gamma]$ be an isolated critical orbit of $f$ in $\widetilde{\mathcal N}$,
with $f(\gamma)=c$, and
let $\Sigma$ be any smooth hypersurface of $\widetilde{\mathcal N}$ with $[\gamma]\in \Sigma$ and with
$T_{[\gamma]}\widetilde{\mathcal N}=T_{[\gamma]}\Sigma\oplus T_{[\gamma]}\big(\mathrm O(2)[\gamma]\big)$.
Then, $[\gamma]$ is an isolated critical point of $f\vert_\Sigma$, and
\[\mathfrak H_*(f\vert_{\mathcal D^\gamma},[\gamma];\mathds K)\cong
H_*\big(\Sigma\cap f^c,(\Sigma\cap f^c)\setminus\{[\gamma]\};\mathds K\big).
\]
Moreover, the Morse indexes and the nullities of $[\gamma]$ as a critical point of $f\vert_{\mathcal D^\gamma}$ and
of $f\vert_\Sigma$ coincide.
\end{lem}
\begin{proof}
Let $\Sigma$ be as above; the entire result will follow from the existence of an $f$-invariant diffeomorphism $\psi$ from (a small
neighborhood of $[\gamma]$ in) $\mathcal D^\gamma$ onto (a small neighborhood of $[\gamma]$ in) $\Sigma$ with $\psi\big([\gamma]\big)=[\gamma]$.
Consider the smooth map $\Sigma\times\mathrm O(2)\ni (u,g)\mapsto gu\in\widetilde{\mathcal N}$; the assumption of transversality
of $\Sigma$ to the orbit $\mathrm O(2)[\gamma]$ implies that the differential of this map at the point $\big([\gamma],1\big)$ is an
isomorphism, hence the map restricts to a diffeomorphism from a neighborhood of $\big([\gamma],1\big)$ to
a neighborhood of $[\gamma]$ in $\widetilde{\mathcal N}$. Since $\mathcal D^\gamma$ is also transversal to
$\mathrm O(2)[\gamma]$, a neighborhood of $[\gamma]$ in $\mathcal D^\gamma$ is diffeomorphic, via this map, to the
graph of a smooth function $\varphi:\widetilde\Sigma\to\mathrm O(2)$, where $\widetilde\Sigma$ is a neighborhood of $[\gamma]$ in $\Sigma$
and $\varphi\big([\gamma]\big)=1$.
The required $f$-invariant diffeomorphism $\psi$ is given by $\widetilde\Sigma\ni u\mapsto\varphi(u)u\in\mathcal D^\gamma$.
\end{proof}
\begin{cor}\label{thm:calchominv0restr}
Under the assumptions of Lemma~\ref{thm:Ckrestrictions}:
\[\mathfrak H_*^o\big(f\vert_{\mathcal D^\gamma},[\gamma];\mathds K\big)\cong\mathfrak H_*^o\big(f\vert_\Sigma,[\gamma];\mathds K\big).\]
\end{cor}
\begin{proof}
Follows immediately from Lemma~\ref{thm:Ckrestrictions} and the shifting theorem \eqref{eq:uguaShTh}.
\end{proof}
\begin{rem}\label{thm:remhessianorestrictionbis}
More generally, from the proof of  Lemma~\ref{thm:Ckrestrictions} we get that if $\Sigma$ is any
hypersurface of $\widetilde{\mathcal N}$ as in the statement, all the properties of $f\vert_{\mathcal D^\gamma}$
discussed in Remark~\ref{thm:remhessianorestriction} also hold for the restriction $f\vert_\Sigma$.
Under the circumstance that $\Sigma$ is a hypersurface through $[\gamma]$ in $\widetilde{\mathcal N}$ that
is orthogonal (relatively to an arbitrary Riemannian metric on $\widetilde{\mathcal N}$) to the critical
orbit $\mathrm O(2)[\gamma]$ at $[\gamma]$, then the null space of the Hessian of $f\vert_{\Sigma}$ at
$[\gamma]$ is the intersection of the null space of the Hessian of $f\vert_{\widetilde{\mathcal N}}$ at $\gamma$
and $T_{[\gamma]}\Sigma$. This follows easily from the observation that $T_{[\gamma]}\big(\mathrm O(2)[\gamma]\big)$
is contained in the kernel of the Hessian of $f\vert_{\widetilde{\mathcal N}}$ at $\gamma$.
\end{rem}
Finally, the key result of this subsection is to show that the local homological invariants at $[\gamma]$ coincide
with the invariants at the iterate $\big[\gamma^{(N)}\big]$ when $\gamma$ and $\gamma^{(N)}$ have
the same nullity (Proposition~\ref{thm:stessiBk0}). It will therefore be necessary to study the
\emph{$N$-times iteration map} $\mathfrak N:\widetilde{\mathcal N}\to\widetilde{\mathcal N}$, defined by
$\mathfrak N\big([\gamma]\big)=\big[\gamma^{(N)}\big]$.

\begin{lem}\label{thm:Nitermapembedding}
$\mathfrak N$ is a smooth embedding.
\end{lem}
\begin{proof}
We use the following criterion, which is proved easily. Let $A,B$ be Banach manifolds and let $A'\subset A$ be an embedded
submanifold. Let $\mathfrak g:A'\to B$ and $\mathfrak h:B\to A$ be smooth maps such that $\mathfrak h\circ\mathfrak g$ is
the inclusion of $A'$ into $A$. Then $\mathfrak g$ is a smooth embedding. In  order to prove the Lemma, the criterion
is used in the following setup. The manifolds $A$ and $B$ are the sets of all curves $\sigma:[0,1]\to M$
of Sobolev class $H^1$, with $\sigma(0)\in S$, and satisfying $g(\dot\sigma,\mathcal Y)$ constant almost everywhere
(the esistence of a Hilbert manifold structure of this set is proved exactly as for $\widetilde{\mathcal N}$).
The submanifold $A'$ is $\widetilde{\mathcal N}$, which corresponds to the subset of $A$ consisting of closed
curves. The map $\mathfrak g$ is the $N$-times iteration map $\mathfrak N$, and the map $\mathfrak h$ is defined by
$\mathfrak h(\sigma)=\tilde\sigma$, and $\tilde\sigma(t)=\sigma(t/N)$, for all $t\in[0,1]$.
\end{proof}

The differential $\mathrm d\mathfrak N_{[\gamma]}$ at $[\gamma]$ is the $N$-times iteration map for
vector fields along $\gamma$. Let us now prove the following central result:

\begin{prop}\label{thm:stessiBk0}
Let $\gamma$ be a closed geodesic in $M$, let $N\ge1$ be fixed, and assume that $\mathrm O(2)\big[\gamma^{(N)}\big]$
is an isolated critical orbit of $f$ in $\widetilde{\mathcal N}$. Then, $\mathrm O(2)[\gamma]$
is an isolated critical orbit of $f$ in $\widetilde{\mathcal N}$, and if $\mathrm n(\gamma)=\mathrm n\big(\gamma^{(N)}\big)$,
one has $C_k^o(\gamma; \mathds K)=C_k^o\big(\gamma^{(N)}; \mathds K\big)$ for all $k$.
\end{prop}
\begin{proof}
The idea of the proof is analogous to that of \cite[Theorem~3]{GroMey2} and  \cite[Proposition~3.6]{GroTan}; several
adaptations are needed due to the fact that we are dealing with different metric structures in the manifold $M$: the Lorentzian structure
$g$ and the Riemannian structure $\gR$ (recall \eqref{eq:defgr}) employed in the definition of the Hilbert structure of $\widetilde{\mathcal N}$.

Consider a modified Riemannian structure on $\widetilde{\mathcal N}$ induced by the inner product (compare with \eqref{eq:metricaHilbert})
on each tangent space $T_\gamma\mathcal N$ given by:
\begin{equation}\label{eq:metricaHilbertmodified}
\llangle V,W\rrangle_N=\int_0^1\Big[N^2\gR(V,W)+\gR\big(\DdtR V,\DdtR W\big)\Big]\,\mathrm dt.
\end{equation}
Consider the $N$-times iteration map $\mathfrak N:\big(\widetilde{\mathcal N},\llangle\cdot,\cdot\rrangle\big)
\to \big(\widetilde{\mathcal N},\llangle\cdot,\cdot\rrangle_N\big)$, which is an embedding
onto a smooth submanifold $\mathfrak N(\widetilde{\mathcal N})$ of $\widetilde{\mathcal N}$ by Lemma~\ref{thm:Nitermapembedding},
and it  preserves the metric up to a factor $N^2$.
We claim that, at the points in the image of the map $\mathfrak N$, the gradient $\nabla^N f$ of the functional $f\vert_{\widetilde{\mathcal N}}$ relatively to
the metric $\llangle\cdot,\cdot\rrangle_N$ is tangent to the image of $\mathfrak N$. The set of points in the image of $\mathfrak N$
where this situation occurs is closed, and so, by a density argument, it suffices to prove the
claim at those points $\sigma^{(N)}=\mathfrak N(\sigma)$ in the image of $\mathfrak N$ that are curves of class $C^2$.
Given one such point $\sigma^{(N)}$, using the fundamental theorem of Calculus of Variations, one sees that the
the gradient $\nabla^N f({\sigma^{(N)}})$  of $f$ at $\sigma^{(N)}$ is the unique periodic vector field $X$ along $\sigma^{(N)}$ that solves the
differential equation:
\begin{equation}\label{eq:secorderegradient}
\DdttR X-N^2X-2\frac{g(\DdttR X-N^2X,\mathcal Y)}{g(\mathcal Y,\mathcal Y)}\,\mathcal Y=\Ddt\ddt\sigma^{(N)}.
\end{equation}
Now, if $X_*$ is the vector field along $\sigma$ which is the unique periodic solution of:
\[\DdttR X_*-X_*-2\frac{g(\DdttR X_*-X_*,\mathcal Y)}{g(\mathcal Y,\mathcal Y)}\,\mathcal Y=\Ddt\ddt\sigma,\]
i.e., $X_*$ is the gradient of $f$ relatively to the metric $\llangle\cdot,\cdot\rrangle$ at $\sigma$, then
the iterate $X_*^{(N)}=\mathrm d\mathfrak N_\sigma(X_*)$ satisfies \eqref{eq:secorderegradient}, which proves the claim.\footnote{%
Observe that $\mathrm d\mathfrak N_\sigma\big(\mathcal Y\vert_\sigma\big)=\mathcal Y\vert_{\sigma^{(N)}}$.}

Let $\Gamma\subset\mathrm{SO}(2)$ be the stabilizer of $\gamma$; consider a normal disc bundle
$\mathcal D=\mathrm O(2)\mathcal D^\gamma\cong\mathrm O(2)\times_\Gamma\mathcal D^\gamma$ of the critical orbit $\mathrm O(2)[\gamma]$
as described in Subsection~\ref{sub:hominvariants}. The image $\mathfrak N(\mathcal D^\gamma)$ is a smooth embedded submanifold of $\widetilde{\mathcal N}$
containing $\gamma^{(N)}$; since $\mathfrak N:\mathcal D^\gamma\to\mathfrak N(\mathcal D^\gamma)$
is a diffeomorphism and $f\circ\mathfrak N=N^2f$, then:
\begin{equation}\label{eq:primaug}
\mathfrak H^o_*\big(f\vert_{\mathfrak N(\mathcal D^\gamma)},
\gamma^{(N)};\mathds K\big)=\mathfrak H^o_*\big(f\vert_{\mathcal D^\gamma},\gamma;\mathds K\big).
\end{equation}
In order to conclude the proof, we will now determine a hypersurface $\Sigma$ in $\widetilde{\mathcal N}$ through $\gamma^{(N)}$ which
is transversal at $\big[\gamma^{(N)}\big]$ to the orbit $\mathrm O(2)\big[\gamma^{(N)}\big]$
and satisfying the following two properties:
\begin{itemize}
\item[(a)] $\mathfrak N\big(\mathcal D^\gamma)\subset\Sigma$;
\item[(b)] the gradient $\nabla^N\big(f\vert_\Sigma\big)$ at the points of $\mathfrak N(\mathcal D^\gamma)$ is tangent to $\mathfrak N(\mathcal D^\gamma)$;
\item[(c)] the null space of the Hessian $\mathrm H^{f\vert_\Sigma}$ at $\big[\gamma^{(N)}\big]$ is contained in $T_{[\gamma^{(N)}]}\mathfrak N(\mathcal D^\gamma)$.
\end{itemize}
By Corollary~\ref{thm:calchominv0restr} it will follow that:
\begin{equation}\label{eq:secondaug}
C_k^o(\gamma^{(N)};\mathds K)=\Dim\Big[H_k\big(\Sigma\cap f^d,(\Sigma\cap f^d)\setminus\{[\gamma^{(N)}]\};\mathds K\big)\Big],\qquad\forall\,k\ge0,
\end{equation}
where $d=f\big(\gamma^{(N)}\big)=cN^2$ and $c=f(\gamma)$. Moreover, using Lemma~\ref{thm:localizationproperty}, properties (a), (b) and
(c) will imply that:
\begin{equation}\label{eq:terzaug}
H_*\big(\Sigma\cap f^d,(\Sigma\cap f^d)\setminus\{[\gamma^{(N)}]\};\mathds K\big)\cong \mathfrak H^o_*\big(f\vert_{\mathfrak N(\mathcal D^\gamma)},
\gamma^{(N)};\mathds K\big).
\end{equation}
The thesis will follow then from \eqref{eq:primaug}, \eqref{eq:secondaug} and \eqref{eq:terzaug}.

For the construction of the desired $\Sigma$, consider be the normal bundle $\nu\big(\mathfrak N(\mathcal D)\big)$  of the submanifold
$\mathfrak N(\mathcal D)$ in $\widetilde{\mathcal N}$ relatively to the metric $\llangle\cdot,\cdot\rrangle_N$. Let
$\widetilde{\mathrm{EXP}}$ be the exponential map of $\widetilde{\mathcal N}$ relatively to the metric $\llangle\cdot,\cdot\rrangle_N$;
define $\Sigma$ to be the image under $\widetilde{\mathrm{EXP}}$ of a small neighborhood $U$ of the zero section of
the bundle $\nu\big(\mathfrak N(\mathcal D)\big)\vert_{\mathfrak N(\mathcal D^\gamma)}$, i.e., the restriction to
$\mathfrak N(\mathcal D^\gamma)$ of the normal bundle of $\mathfrak N(\mathcal D)$.
Since $\mathcal D^\gamma$ is a hypersurface
in $\mathcal D$, if $U$ is sufficiently small, then $\Sigma$ is a hypersurface in $\widetilde{\mathcal N}$;
clearly, $\mathfrak N(\mathcal D^\gamma)\subset\Sigma$.

The image $\mathfrak N\big(\mathrm{SO}(2)[\gamma]\big)$ coincides with the orbit $\mathrm{SO}(2)\big[\gamma^{(N)}\big]$;
this is easily seen observing that the map $\mathrm{SO}(2)\ni g\mapsto g^N\in\mathrm{SO}(2)$ is surjective.
Since $\mathcal D^\gamma$ is orthogonal to $\mathrm O(2)[\gamma]$ and $\mathfrak N$ is metric preserving up to a constant factor,
it follows that $\Sigma$ is orthogonal to $\mathrm O(2)\big[\gamma^{(N)}\big]$ at $\big[\gamma^{(N)}\big]$
(observe that $\big[\gamma^{(N)}\big]$ belongs to the connected component $\mathrm{SO}(2)\big[\gamma^{(N)}\big]$
of $\mathrm O(2)\big[\gamma^{(N)}\big]$) relatively to the metric $\llangle\cdot,\cdot\rrangle_N$.

For $u\in\mathfrak N(\mathcal D^\gamma)$, the tangent space
$T_u\Sigma$ is given by the orthogonal direct sum (see Lemma~\ref{thm:tangspacimexp} below):
\begin{equation}\label{eq:tangspSigma}
T_u\Sigma=T_u\big(\mathfrak N(\mathcal D^\gamma)\big)\stackrel\perp\oplus T_u\big(\mathfrak N(\mathcal D)\big)^\perp.
\end{equation}
From the first part of the proof we know that at the points $u\in \mathfrak N(\mathcal D^\gamma)$, the gradient
$\nabla^N f(u)$ is tangent to $\mathfrak N(\mathcal D)$; from \eqref{eq:tangspSigma}, the orthogonal
projection of $\nabla^N f(u)$ onto $T_u\Sigma$, which is the gradient of $f\vert_\Sigma$ at $u$, must be tangent
$\mathfrak N(\mathcal D^\gamma)$. Property (b) is thus satisfied.

Finally, we claim that the differential $\mathrm d\mathfrak N_{[\gamma]}$ of $\mathfrak N$ at $[\gamma]$ carries
the null space of the Hessian of $f\vert_{\mathcal D^\gamma}$ at $[\gamma]$ (injectively) into the null
space of the Hessian of $f\vert_{\Sigma}$ at $[\gamma]$. Namely, recall from Remark~\ref{thm:remhessianorestrictionbis}
that the null space of the Hessian of $f\vert_{\mathcal D^\gamma}$ (resp., of $f\vert_\Sigma$) at
$[\gamma]$ (resp., at $\big[\gamma^{(N)}\big]$) consists of all periodic Jacobi fields that are
orthogonal to the critical orbit $\mathrm O(2)[\gamma]$ (resp., $\mathrm O(2)\big[\gamma^{(N)}\big]$).
Thus, the proof of the claim follows easily observing that the map $\mathrm d\mathfrak N_{[\gamma]}$:
\begin{itemize}
\item carries periodic Jacobi fields along $\gamma$ to periodic Jacobi fields along $\gamma^{(N)}$;
\item carries $T_{[\gamma]}\big(\mathrm O(2)[\gamma]\big)$ isomorphically onto $T_{[\gamma^{(N)}]}\big(\mathrm O(2)\big[\gamma^{(N)}\big]\big)$;
\item preserves orthogonality.
\end{itemize}
The null spaces of the two Hessians have the same dimension, because
of our assumption on the nullity of $[\gamma]$ and of $\big[\gamma^{(N)}\big]$ (recall from Proposition~\ref{thm:PS}
and Remarks~\ref{thm:remhessianorestriction},
\ref{thm:remhessianorestrictionbis} that these two spaces have dimensions $\mathrm n(\gamma)-2$ and $\mathrm n\big(\gamma^{(N)}\big)-2$ respectively).
This implies that the null space of $\mathrm H^{f\vert_\Sigma}$ at $\big[\gamma^{(N)}\big]$ is in the image
of $\mathrm d\mathfrak N_{[\gamma]}$, hence it is contained in $T_{[\gamma^{(N)}]}\mathfrak N(\mathcal D^\gamma)$, which gives property (c).
This concludes the proof.
\end{proof}
Lemma~\ref{thm:tangspacimexp} below has been used in the proof of Proposition~\ref{thm:stessiBk0} to the following setup:
$A=\widetilde{\mathcal N}$, $B=\mathfrak N(\mathcal D^\gamma)$ and $E=\nu\big(\mathfrak N(\mathcal D)\big)\vert_{\mathfrak N(\mathcal D^\gamma)}$.
\begin{lem}\label{thm:tangspacimexp}
Let $A$ be a Hilbert manifold and $B\subset A$ a submanifold. Let $\nu(B)\subset TA$ be the normal bundle of $B$ in $A$ and let
$E\subset\nu(B)$ be a subbundle. Let $U\subset\nu(B)$ be a small open subset containing the zero section and
set $\Sigma=\exp(U\cap E)$. Then, $B$ is a submanifold of $\Sigma$, and for all $b\in B$, the tangent space
$T_b\Sigma$ is the orthogonal direct sum $T_bB\oplus E_b$.
\end{lem}
\begin{proof}
$B$ is the image of the zero section $\mathbf 0$ of $E$.
At each point $0_b\in\mathbf 0$, $b\in B$, there is a canonical isomorphism $T_{0_b}E\cong T_bB\oplus E_b$, where
$T_bB$ is identified with the tangent space at $0_b$ of $\mathbf 0$. Using this identification, the differential
$\mathrm d(\exp\vert_{\mathbf 0})(0_b):T_bB\oplus\{0\}\to T_bB$ of the restriction of $\exp$ to $\mathbf 0$ at $0_b$ is the identity. Moreover,
the restriction of $\mathrm d\exp(0_b)$ to $\{0\}\oplus E_b$ coincides with the differential
$\mathrm d\exp_b(0_b)$, which is the identity. Thus, $\mathrm d\exp(0_b)$ carries $T_{0_b}E$ isomorphically onto
$T_bB\oplus E_b$, and the conclusion follows.
\end{proof}

\subsection{Equivariant Morse theory for closed geodesics}
As observed in Remark~\ref{thm:remazioneO2}, in order to prove the theorem we need to show the existence
of infinitely many distinct prime critical $\mathrm O(2)$-orbits of the functional $f$ in $\widetilde{\mathcal N}$;
this will be obtained by contradiction, showing that assuming the existence of only a finite number of geometrically distinct
closed geodesics will yield a uniform upper bound on the Betti numbers of $\Lambda M$.

Let us assume that there is only a finite number of geometrically distinct critical orbits, hence, by Lemma~\ref{thm:isolatedcriticalorbits},
the critical orbits of $f$ in $\widetilde{\mathcal N}$ are isolated.
If $0\le a<b$ are regular values of $f$, and if $\mathrm O(2)\big[\gamma_1\big]$, \ldots, $\mathrm O(2)\big[\gamma_r\big]$ are all
the critical orbits of $f$ in $f^{-1}\big([a,b]\big)$, then, using \eqref{eq:sommauguale} and the fact that the $\beta_k$'s
are subadditive functions,  one has the Morse inequalities:
\begin{equation}\label{eq:Morseinequalities}
\beta_k(f^b,f^a;\mathds K)\le\sum_{j=1}^r B_k(\gamma_j;\mathds K).
\end{equation}
In particular, since $\mathfrak H$ is of finite type, i.e., $B_k(\gamma;\mathds K)$ is finite for all $k$ and $B_k(\gamma;\mathds K)=0$
except for a finite number of $k$'s, then $\beta_k(f^b,f^a;\mathds K)<+\infty$ for all $a,b$ and $k$.

Using the relative Mayer--Vietoris sequence to the pair of bundles \eqref{eq:pairofbundles}
over $\mathrm O(2)$, which is homeomorphic to the disjoint union
of two copies of the circle,
one proves that the following inequality:
\begin{equation}\label{eq:disugHkorbHkpto}
B_k(\gamma;\mathds K)\le2\big(C_{k}(\gamma,\mathds K)+C_{k-1}(\gamma,\mathds K)\big),
\end{equation}
holds for all $k\ge1$. The details of this computation will be given in Appendix~\ref{sec:homfiberS1}; it should be observed that
in \cite{GroMey2, GroTan} the inequality is stated only in the case of a field $\mathds K$ of characteristic zero.

By the shifting theorem (see \eqref{eq:uguaShTh}), inequalities \eqref{eq:disugHkorbHkpto} become:
\begin{equation}\label{eq:stimaKkB0k}
B_k(\gamma;\mathds K)\le 2\big(C^o_{k-\mu(\gamma)}(\gamma;\mathds K)+C^o_{k-\mu(\gamma)-1}(\gamma;\mathds K)\big).
\end{equation}

\begin{prop}\label{thm:B0k(gN)unifbounded}
Let $\gamma$ be a closed geodesic in $M$. If all the critical orbits $\mathrm O(2)\big[\gamma^{(N)}\big]$
of $f$ in $\widetilde{\mathcal N}$ are isolated, then the double sequence $(k,N)\mapsto C^o_k\big(\gamma^{(N)}; \mathds K\big)$
is uniformly bounded:
\begin{equation}\label{eq:B0k(gN)unifbounded}
\phantom{,\quad\forall\,k,N\in\N\setminus\{0\}.}C_k^o\big(\gamma^{(N)}; \mathds K\big)\le B,\quad\forall\,k,N\in\N\setminus\{0\}.
\end{equation}
Moreover, there exists $k_0$ such that $C_k^o\big(\gamma^{(N)}; \mathds K\big)=0$ for all $k>k_0$ and all $N\ge1$.
\end{prop}
\begin{proof}
Inequality \eqref{eq:B0k(gN)unifbounded} follows readily from Lemma~\ref{thm:nuliteration} (see Remark~\ref{thm:remfiniternul})
and Proposition~\ref{thm:stessiBk0}. For a fixed $N$, the existence of $k_0$ as above is guaranteed by the fact
that the invariant $\mathfrak H^o$ is of finite type. Again, independence on $N$ is obtained easily
from Lemma~\ref{thm:nuliteration} and Proposition~\ref{thm:stessiBk0}.
\end{proof}
\begin{cor}\label{thm:cor2}
Under the assumptions of Proposition~\ref{thm:B0k(gN)unifbounded}, the following inequality holds:
\begin{equation}\label{eq:stimaBk}
\phantom{,\quad\forall\,N\ge1}B_k\big(\gamma^{(N)}; \mathds K\big)\le4B,\quad\forall\,N\ge1.
\end{equation}
Moreover, for $k>k_0+8\,\Dim(M)+2$, the number of iterates $\gamma^{(N)}$ of $\gamma$ such that $B_k\big(\gamma^{(N)}; \mathds K\big)\ne0$ is bounded
by a constant $C$ which does not depend on $k$.
\end{cor}
\begin{proof}
Inequality \eqref{eq:stimaBk} follows from \eqref{eq:stimaKkB0k} and \eqref{eq:B0k(gN)unifbounded}.
Moreover, using \eqref{eq:stimaKkB0k} and Proposition~\ref{thm:B0k(gN)unifbounded} we get that $B_k\big(\gamma^{(N)}; \mathds K\big)\ne0$
only if:
\begin{equation}\label{eq:condizioneBknonzero}
k-k_0-1\le\mu\big(\gamma^{(N)}\big)\le k.
\end{equation}
If the sequence $\mu\big(\gamma^{(N)}\big)$ is bounded, then by our assumption on $k$ and
Proposition~\ref{thm:lineargrowth}, no iterate $\gamma^{(N)}$ of $\gamma$ satisfies \eqref{eq:condizioneBknonzero}.
Assume that $\mu\big(\gamma^{(N)}\big)$ is not bounded, and let $k_*$ be as in the
proof of Proposition~\ref{thm:lineargrowth}. Let $\bar k\ge k_*$ be the smallest integer
for which $\mu\big(\gamma^{(\bar k)}\big)\ge k-k_0-1$; we need to estimate the numbers of positive
integers $s$ such that $\mu\big(\gamma^{(\bar k+s)}\big)\le k$. If $s\ge k_*$, then by
Corollary~\ref{thm:corcrescitalinearearb}:
\[k_0+1\ge\mu\big(\gamma^{(\bar k+s)}\big)-\mu\big(\gamma^{(\bar k)}\big)\ge\bar\alpha\,s-\bar\beta,\]
where $ \bar\alpha,\bar\beta>0$. Thus, the number of iterates $\gamma^{(N)}$ such that $B_k\big(\gamma^{(N)}; \mathds K\big)\ne0$ is bounded
by the constant:
\[\max\left\{k_*,\frac{k_0+1+\bar\beta}{\bar\alpha}\right\}.\qedhere\]
\end{proof}
\begin{prop}\label{thm:boundedbettinumbers}
Let $(M,g)$ be a Lorentzian manifold that has a complete timelike Killing vector field
and a compact Cauchy surface. If there is only a finite number of geometrically distinct
non trivial closed geodesics in $M$, then the Betti numbers $\beta_k(\Lambda M;\mathds K)$ form a
bounded sequence for $k$ large enough.
\end{prop}
\begin{proof}
Since $\Lambda M$ is homotopically equivalent to $\widetilde{\mathcal N}$,
$\beta_k(\Lambda M;\mathds K)=\beta_k(\widetilde{\mathcal N};\mathds K)$ for all $k\ge0$.
Denote by $\gamma_1,\ldots,\gamma_r$ a maximal family of pairwise geometrically distinct non trivial
closed geodesic in $M$, and let $0=c_0<c_1<\ldots<c_n<\ldots$ be the critical values of $f$ in
$\widetilde{\mathcal N}$ corresponding to the critical orbits $\mathrm O(2)\big[\gamma_i^{(N)}\big]$,
$N\ge1$, $i=1,\ldots,r$. By Lemma~\ref{thm:isolatedcriticalorbits}, these critical orbits are isolated,
and each fixed sublevel $f^b$ of $f$ in $\widetilde{\mathcal N}$ contains only a finite number of
them. By Corollary~\ref{thm:cor2}, the sequence $(i,k,N)\mapsto B_k\big(\gamma_i^{(N)}; \mathds K\big)$
takes a finite number of values, and we can define $\widehat B=\max\limits_{i,k,N}B_k\big(\gamma_i^{(N)}; \mathds K\big)$.

For each geodesic $\gamma_i$, choose numbers $k_0^{(i)}$ and $C^{(i)}$ as in
Proposition~\ref{thm:B0k(gN)unifbounded} and Corollary~\ref{thm:cor2}; set $\hat k_0=\max\limits_ik_0^{(i)}$ and
$\widehat C=\max\limits_iC^{(i)}$.

By Corollary~\ref{thm:cor2}, for all $k>\hat k_0+8\,\Dim(M)+2$, the
constant $\widehat C$ is an upper bound for the number of orbits $\mathrm O(2)\big[\gamma_i^{(N)}\big]$
with $B_k\big(\gamma_i^{(N)}; \mathds K\big)\ne0$. Using the Morse inequalities \eqref{eq:Morseinequalities},
we have that for all regular values $a,b$ of $f$ in $\widetilde{\mathcal N}$, with $0<a<b$,  and for all
$k>\hat k_0+8\,\Dim(M)+2$ the following inequality holds:
\begin{equation}
\label{eq:disugBettirel}
\beta_k(f^b,f^a;\mathds K)\le 4\widehat B\widehat C.
\end{equation}
By Lemma~\ref{thm:defretract}, there exists $\varepsilon\in\left]0,c_1\right[$  such that the
sublevel $f^\varepsilon$ is homotopically equivalent to a Cauchy surface $S$ of $M$.
For all $n\ge1$, set $d_n=\frac12(c_n+c_{n+1})$, $d_0=\varepsilon$, and for all $n\ge0$ set $X_n=f^{d_n}$; each $d_n$ is a regular value of $f$ in $\widetilde{\mathcal N}$,
and the $X_n$'s form a filtration of $\widetilde{\mathcal N}$ as in Proposition~\ref{thm:Morseinequalities}.
Since $X_0$ is homotopically equivalent to $S$,  which is a finite dimensional compact manifold,
for $k$ large enough, $\beta_k(\widetilde{\mathcal N},X_0;\mathds K)= \beta_k(\widetilde{\mathcal N};\mathds K)$.
We claim that, for $k>\hat k_0+8\,\Dim(M)+2$, the number of indices $n$ such that $\beta_k(X_{n+1},X_n;\mathds K)\ne0$
is bounded by a constant $N_0$ that does not depend on $k$. Namely, arguing as in the proof of Corollary~\ref{thm:cor2},
one proves easily that such constant $N_0$ can be taken equal to $\sum_{i=1}^rC^{(i)}$.

Now, using \eqref{eq:disugBettirel}, it follows that:
\[
\sum_{n=0}^\infty\beta_k(X_{n+1},X_n;\mathds K)\le 4\widehat B\widehat CN_0,\]
for all $k>\hat k_0+8\,\Dim(M)+2$. Using Proposition~\ref{thm:Morseinequalities} (and the weak Morse inequalities \eqref{eq:weakMorserel}), we get:
\[\beta_k(\Lambda M;\mathds K)=\beta_k(\widetilde{\mathcal N};\mathds K)=\beta_k(\widetilde{\mathcal N}, X_0;\mathds K)\le 4\widehat B\widehat CN_0\]
for $k$ large enough, which concludes the proof.
\end{proof}

We are now in the position of finalizing the proof of our main result.

\begin{proof}[Proof of the main theorem]
Assume that $(M,g)$ is a simply connected stationary globally hyperbolic spacetime, having a
compact Cauchy surface $S$ and a complete timelike Killing vector field $\mathcal Y$.
\marginpar{\textbf{I risultati di Serre\\ valgono con coeffi\-cien\-ti in qualsiasi\\ corpo?}}
Then, $S$ is simply connected and, by \cite{Ser}, the Betti numbers of the free loop space of $S$
(or, equivalently, of $M$) are finite.
Then, by Proposition~\ref{thm:boundedbettinumbers}, the finiteness of the number of geometrically distinct closed
geodesics in $M$ implies that the Betti numbers of $\Lambda M$ form a bounded sequence. The thesis follows.
\end{proof}
\end{section}

\begin{section}{Final remarks}
A few observations on the result presented in the paper and its proof are in order.

\begin{rem}\label{thm:remindependence}
As to the notion of geometric equivalence for closed geodesics given in the Introduction,
\marginpar{\textbf{C'\`e un argomento\\ pi\`u semplice?}}
and based on the choice of some complete timelike Killing vector field, we observe
that the property of existence of infinitely many geometrically distinct closed geodesic
in independent on such choice. This can be seen using the following construction.
Assume that $S\subset M$ is a Cauchy surface of $(M,g)$; given a complete timelike
Killing vector field $\mathcal Y$, one can define a diffeomorphism
$\mathbb P_{\mathcal Y}:\widetilde{\mathcal N}\to\Lambda S$ by considering \emph{projections}
onto $S$ along the flow lines of $\mathcal Y$ (note that also the definition of $\widetilde{\mathcal N}$
employs the given vector field $\mathcal Y$). More precisely, given $\gamma\in\widetilde{\mathcal N}$,
the curve $x=\mathbb P_{\mathcal Y}(\gamma)$ is defined by $x(t)=\mathcal F_{h_\gamma(t)}\big(\gamma(t)\big)$,
where $\mathcal F$ is the flow of $\mathcal Y$ and $h_\gamma:[0,1]\to\R$ is uniquely
defined by the property that $\mathcal F_{h_\gamma(t)}\big(\gamma(t)\big)\in S$.
By an elementary ODE argument, it is easy to see that $\mathbb P_{\mathcal Y}$ is indeed a bijection, by proving
that, given $x\in\Lambda S$, there exists a unique closed curve $\gamma$ with
$\gamma(0)=x(0)$ such that $\mathbb P_{\mathcal Y}(\gamma)=x$ and such that
$g(\dot\gamma,\mathcal Y)$ is constant. The smoothness of $\mathbb P_{\mathcal Y}$
is obtained by standard smooth dependence results for ODE's.
The map $\mathbb P_{\mathcal Y}$ is $\mathrm O(2)$-equivariant; thus,  geometrically distinct closed
geodesics in $M$ correspond to  distinct critical
$\mathrm O(2)$-orbits of the functional $\mathfrak f_{\mathcal Y}=f\circ\mathbb P_{\mathcal Y}^{-1}:\Lambda S\to\R$
(this is precisely the variational problem considered in \cite{Masiello2}).
Given two complete timelike Killing vector field $\mathcal Y_1$ and $\mathcal Y_2$ in $M$,
the number of critical $\mathrm O(2)$-orbits of the functionals $\mathfrak f_1$ and
$\mathfrak f_2=\mathfrak f_1\circ \mathbb P_{\mathcal Y_1}\circ\mathbb P_{\mathcal Y_2}^{-1}$
on $\Lambda S$ coincide, which proves that the number of geometrically distinct closed
geodesic in $(M,g)$ is an intrinsic notion.
\end{rem}

\begin{rem}
Under the assumptions of our main result, if in addition the Killing vector field $\mathcal Y$ is irrotational,
i.e., if the orthogonal distribution $\mathcal Y^\perp$ is integrable, then the proof of our result is immediate.
Namely, in this situation, a maximal integrable submanifold $S$ of $\mathcal Y^\perp$ is a compact \emph{totally geodesic}
Cauchy surface in $(M,g)$. Thus, infinitely many closed geodesics in $M$ can be obtained applying the classical Gromoll
and Meyer result to the Riemannian manifold $(S,g\vert_S)$.
\end{rem}

\begin{rem} It must be emphasized that the estimates on the Conley--Zehnder index
and the Maslov index discussed in Section~\ref{sec:iteformMaslov} are very far from
being sharp, and they only serve the purposes of the present paper. An intense
literature on the iteration formulas in the context of periodic solutions of Hamiltonians on
symplectic manifolds has been produced in the last decade (see for instance \cite{CusDui, Lon1, Lon2}
and the references therein). On the other hand, the naive approach discussed in
Section~\ref{sec:iteformMaslov} seems to simplify significantly the approach using
Bott's deep results in \cite{Bo4} on the Morse index of an iteration, even in the Riemannian case.
\end{rem}

\begin{rem}\label{thm:remnonsimplyconnected}
As to the assumption that $M$ be simply connected, one should note that the central
result in Proposition~\ref{thm:boundedbettinumbers} does not use this. The simple
connectedness hypotheses is used in the final argument to guarantee
the finiteness of the dimensions of \emph{all} the homology spaces of the free loop space
of $M$, by a result on spectral sequences due to Serre \cite{Ser}.
Observe that Proposition~\ref{thm:boundedbettinumbers} does not give any information
on the dimension of the homology spaces $\beta_k(\Lambda M;\mathds K)$ for
$k=0,\ldots,\hat k_0+8\,\Dim(M)+2$.
As already observed in \cite{GroTan}, if $M$ is not simply connected, then $\Lambda M$
(and $\widetilde{\mathcal N}$) is not connected, and it might be the case that
$\beta_k(\Lambda M;\mathds K)=+\infty$ for small values of $k$ even if $(M,g)$ has
only a finite number of geometrically distinct non trivial closed geodesics. This
might happen when there is a non trivial closed geodesic whose iterates have bounded
Morse indexes.  Thus, one can state the main result of the paper in the following
slightly more general form:
\begin{teon}
Let $(M,g)$ be a  globally hyperbolic stationary Lorentzian manifold having
a complete timelike Killing vector field, and having a compact Cauchy surface.
Assume that the free loop space $\Lambda M$ has  Betti numbers $\beta_k$
with respect to some coefficient field that satisfy:
\[\limsup_{k\to\infty}\beta_k=+\infty.\]
Then, there are infinitely many geometrically distinct non trivial closed geodesics in $M$.
\end{teon}
\end{rem}

\begin{rem}
Although it is clear how to
produce examples of non trivial closed geodesics all of whose iterates have null
Morse index (any minimum of $f$ in a nontrivial free homotopy class of $M$),
it would be extremely interesting to produce Lorentzian examples having bounded,
but non zero, Morse indexes. The homology generated by the iterates of such closed
geodesics might be richer than the homology of the free loop space, as described
in \cite{BanKli} for the Riemannian case.
\end{rem}

\begin{rem}
Extensions of the result of existence of multiple closed geodesics in Lorentzian geometry are possible,
and indeed desirable, in more general classes of manifolds. The non simply connected case can be studied
following the lines of the corresponding results in Riemannian geometry, as in \cite{BalThoZil, BanHin}.
Finally, we observe that, in view to applications to General Relativity, it would be interesting
to establish multiplicity results for (causal) geodesics satisfying more general boundary conditions.
A particularly interesting case is that of causal geodesics whose spatial component is periodic.
In the stationary case such geodesics have endpoints related by a global isometry of the spacetime,
and an analysis of this case might be based on a variational setup as in \cite{GroHalVig, GroHal, GroTan, Tan}.
\end{rem}
\end{section}

\appendix
\begin{section}{An estimate on the relative homology of fiber bundles over $\mathds S^1$}
\label{sec:homfiberS1}
In this short appendix we will prove a result on the relative homology of fiber bundles
over the circle with coefficients in an arbitrary field $\mathds K$, that will allow a slight generalization
of the result of Gromoll and Mayer.

\begin{prop}\label{thm:homologybundles}
Let $\mathds K$ be a field, and let $\pi:E\to\mathds S^1$ be a fiber bundle with typical fiber $E_0$. Let
$E'\subset E$ and $E_0'\subset E_0$ be subsets such that for all $p\in\mathds S^1$ there exists a trivialization
$\phi_p:\pi^{-1}\big(\mathds S^1\setminus\{p\}\big)\to \big(\mathds S^1\setminus\{p\}\big)\times E_0$
whose restriction to $\pi^{-1}\big(\mathds S^1\setminus\{p\}\big)\cap E'$ gives a homeomorphism with $\big(\mathds S^1\setminus\{p\}\big)\times E_0'$.
Then, for all $k\ge0$, the following inequality holds:
\[\dim\big(H_k(E,E';\mathds K)\big)\le \dim\big(H_k(E_0,E_0';\mathds K)\big)+\dim\big(H_{k-1}(E_0,E_0';\mathds K)\big).\]
\end{prop}
\begin{proof}
Consider two distinct points $p_1,p_2\in\mathds S^1$ and set:
\[X_i=\pi^{-1}\big(\mathds S^1\setminus\{p_i\}\big),\quad A_i=X_i\cap E',\quad i=1,2,\]
so that $E=X_1\bigcup X_2$ and $E'=A_1\bigcup A_2$. The pairs $(X_1,X_2)$ and $(A_1,A_2)$ are
excisive couples for $E$ and $E'$ respectively, since $X_i$ is open in $X$ and $A_i$ is open in $A$, $i=1,2$.
Hence, there is an exact sequence (Mayer--Vietoris, see for instance \cite[\S~8.1]{MawWil}):
\begin{multline*}
\cdots\to H_k(X_1\cap X_2,A_1\cap A_2;\mathds K)\stackrel{\alpha^k_1\oplus\alpha^k_2}{\longrightarrow} H_k(X_1,A_1;\mathds K)\oplus H_k(X_2,A_2;\mathds K)\\ \to H_k(E,E';\mathds K)\to
 H_{k-1}(X_1\cap X_2,A_1\cap A_2;\mathds K)\stackrel{\alpha^{k-1}_1\oplus\alpha^{k-1}_2}{\longrightarrow} \cdots
 \end{multline*}
 Clearly,
 \[X_1\cap X_2=\pi^{-1}\big(\mathds S^1\setminus\{p_1,p_2\}\big),\quad A_1\cap A_2=\pi^{-1}\big(\mathds S^1\setminus\{p_1,p_2\}\big)\cap E'.\]
 We will determine an estimate on the size of the image and the kernel of the map:
 \[\alpha_1^j:H_j(X_1\cap X_2,A_1\cap A_2;\mathds K)\longrightarrow H_j(X_1,A_1;\mathds K), \]
 $j\ge0$, that is induced by the inclusion   $\mathfrak i_1:(X_1\cap X_2,A_1\cap A_2)\to (X_1,A_1)$.
 Choose a trivialization  $\phi:\pi^{-1}\big(\mathds S^1\setminus\{p_1\}\big)\to \big(\mathds S^1\setminus\{p_1\}\big)\times E_0$
 compatible with $E'$ as in the assumptions, and denote by $\widetilde\phi$ the restriction of $\phi$ to
 $\pi^{-1}\big(\mathds S^1\setminus\{p_1,p_2\}\big)$. We have induced isomorphisms:
 \[\xymatrix@C-30pt{H_j\Big(\pi^{-1}\big(\mathds S^1\setminus\{\!p_1\!\}\big),\pi^{-1}\big(\mathds S^1\setminus\{\!p_1\!\}\big)\!\cap\! E';\mathds K\Big)
\ar[d]_{\phi_*} \\ H_j\big((\mathds S^1\setminus\{\!p_i\})\!\times\! E_0,
 (\mathds S^1\setminus\{\!p_1\!\})\!\times\! E_0';\mathds K\big)&\cong\! H_j(E_0,E_0';\mathds K),}\]
 \[\xymatrix@C-30pt{H_j\Big(\pi^{-1}\big(\mathds S^1\setminus\{\!p_1,p_2\!\}\big),\pi^{-1}\big(\mathds S^1\setminus\{\!p_1,p_2\!\}\big)\!\cap\! E';\mathds K\Big)
\ar[d]_{\widetilde\phi_*} \\H_j\big((\mathds S^1\setminus\{\!p_1,p_2\!\})\!\times\! E_0,
 (\mathds S^1\setminus\{\!p_1,p_2\!\})\!\times\! E_0';\mathds K\big)&\cong\! H_j(E_0,E_0';\mathds K)\!\oplus\! H_j(E_0,E_0';\mathds K).}\]
It is immediate to verify that the map \[\phi_*\circ\alpha_1^j\circ\widetilde\phi_*^{-1}:H_j(E_0,E_0';\mathds K)\oplus H_j(E_0,E_0';\mathds K)\to H_j(E_0,E_0';\mathds K)\]
is the sum $(x,y)\mapsto x+y$, which is surjective.
It follows that the dimension of the  image of the map $\alpha_1^j\oplus\alpha_2^j$ is greater than or equal to $\Dim\big(H_j(E_0,E_0';\mathds K)\big)$,
while the kernel of $\alpha_1^j\oplus\alpha_2^j$ has dimension less than or equal to $\Dim\big(H_j(E_0,E_0';\mathds K)\big)$.
From the Mayer--Vietoris sequence, we now pass to the short exact sequence:
\[0\to V_k\to H_k(E,E';\mathds K)\to \mathrm{Ker}(\alpha_1^{k-1}\oplus\alpha_2^{k-1})\to 0,\]
where \[V_k=\big(H_k(E_0,E_0';\mathds K)\oplus H_k(E_0,E_0';\mathds K)\big)/\mathrm{Im}(\alpha_1^k\oplus\alpha_2^k),\]
obtaining:
\begin{multline*}
\Dim\big(H_k(E,E';\mathds K)=\Dim(V_k)+\Dim\big(\mathrm{Ker}(\alpha_1^{k-1}\oplus\alpha_2^{k-1})\big)\\ \le
\dim\big(H_k(E_0,E_0';\mathds K)\big)+\dim\big(H_{k-1}(E_0,E_0';\mathds K)\big).\qedhere\end{multline*}

\end{proof}
An example where Proposition~\ref{thm:homologybundles} applies
is given by considering fiber bundles $E$ that are \emph{associated bundles} $P\times_GE_0$
of a $G$-principal fiber bundle $P$ over $\mathds S^1$, where $E_0$ is a $G$-space (i.e., a topological space endowed with a continuous left $G$-action),
$E_0'\subset E_0$ is a $G$-subspace of $E_0$, and $E'=P\times_GE_0'$ (see \cite[Ch.~1]{PicTauGstruct}).
This is the situation in which Proposition~\ref{thm:homologybundles} is used in the present paper
(recall the definitions of the pair of bundles \eqref{eq:pairofbundles}).
\end{section}

\end{document}